\documentclass[red,11pt,a4paper]{article}

\usepackage{cite}

\usepackage{amsmath}
\usepackage{amscd}
\usepackage{amssymb}
\usepackage{latexsym}
\usepackage{color}

\newcommand{\erl}{\mbox{e}^{r_l}}
\newcommand{\erlN}{\mbox{e}^{r_{l,N}}}

\newcommand{\erk}{\mbox{e}^{r_k}}
\newcommand{\ezk}{\mbox{e}^{z_k}}

\newcommand{\erkN}{\mbox{e}^{r_{k,N}}}

\newcommand{\erki}[1]{\mbox{e}^{r_{k,#1}}}
\newcommand{\ezki}[1]{\mbox{e}^{z_{k,#1}}}

\newcommand{\lap}{\Delta}

\newcommand{\Div}{\operatorname{div}}

\newcommand{\pf}{{\noindent\it Proof.~}}
\newcommand{\Ov}[1]{\overline{#1}}

\newcommand{\Un}[1]{\underline{#1}}

\newcommand{\vC}{{C}}
\newcommand{\vy}{Y_{k}}
\newcommand{\vf}{{\vc F}}

\newcommand{\vQ}{\vc{Q}}
\newcommand{\vS}{{\bf S}}
\newcommand{\vw}{\omega}

\newcommand{\vr}{\varrho}

\newcommand{\vt}{\vartheta}

\newcommand{\vu}{\vc{u}}

\newcommand{\vd}{\vc{d}_{k}}
\newcommand{\vc}[1]{{\bf #1}}
\newcommand{\vcg}[1]{{\pmb #1}}

\newcommand{\Grad}{\nabla}

\newcommand{\pt}{\partial_{t}}
\newcommand{\ptb}[1]{\partial_{t}(#1)}
\newcommand{\Dt}{\frac{ d}{dt}}

\newcommand{\dx}{{\rm d} {x}}
\newcommand{\dt}{{\rm d} t }

\newcommand{\dxdt}{\dx \ \dt}
\newcommand{\lr}[1]{\left( #1 \right)}
\newcommand{\intO}[1]{\int_{\Omega} #1 \ \dx}

\newcommand{\intOB}[1]{\int_{\Omega} \left( #1 \right) \ \dx}

\newcommand{\intT}[1]{\int_0^T #1 \ \dt}
\newcommand{\intTO}[1]{\int_0^T\!\!\!\! \int_{\Omega} #1 \ \dxdt}

\newcommand{\inttO}[1]{\int_0^t\!\!\!\! \int_{\Omega} #1 \ \dxdt}
\newcommand{\inttauO}[1]{\int_0^\tau\!\!\!\! \int_{\Omega} #1 \ \dxdt}

\newcommand{\intTOB}[1]{ \int_0^T\!\!\!\! \int_{\Omega} \left( #1 \right) \ \dxdt}
\newcommand{\inttOB}[1]{ \int_0^t\!\!\!\! \int_{\Omega} \left( #1 \right) \ \dxdt}
\newcommand{\sumkN}[1]{\sum_{k=1}^n #1}
\newcommand{\sumlN}[1]{\sum_{l=1}^n #1}
\newcommand{\eq}[1]{\begin{equation}
\begin{split}
#1
\end{split}
\end{equation}}
\newcommand{\D}{{\vc{D}}}
\newcommand{\ep}{\varepsilon}

\newcommand{\R}{\mathbb{R}}
\newtheorem{thm}{Theorem}
\newtheorem{lemma}[thm]{Lemma}

\newtheorem{rmk}{Remark}

\title{Heat-conducting, compressible mixtures with multicomponent diffusion: construction of a weak solution}
\author{P.B. Mucha$^1$, M. Pokorn\'y$^2$, E. Zatorska$^{1,3}$}
\date{\today}

\begin{document}
\maketitle
\normalsize

\medskip

\begin{center}

1. Institute of Applied Mathematics and Mechanics, University of Warsaw, \\
ul Banacha 2, 02-097 Warszawa, Poland.

2. Charles University in Prague, Faculty of Mathematics and Physics, \\
Mathematical Institute of Charles University, \\
Sokolovsk\'a 83, 186 75 Praha, Czech Republic.

3.  Centre de Math\'{e}matiques Appliqu\'{e}es,
      \'{E}cole Polytechnique,\\
          91128 Palaiseau Cedex, France.

\bigskip

E-mails: p.mucha@mimuw.edu.pl, pokorny@karlin.mff.cuni.cz,\\ e.zatorska@mimuw.edu.pl

\end{center}

\maketitle
\normalsize

\medskip

\noindent{\bf Abstract:} We investigate a coupling between the compressible Navier--Stokes--Fourier system and the full Maxwell-Stefan equations. This model describes the motion of chemically reacting heat-conducting gaseous mixture. The viscosity coefficients are density-dependent functions vanishing on vacuum and the internal pressure depends on species concentrations. 
By several levels of approximation we prove the global-in-time existence of weak solutions on the three-dimensional torus. 

\medskip


\section{Introduction}

Mathematical modelling of mixtures encounters various problems due to wide range of chemical and mechanical aspects that can become important for a particular phenomena.  It turns out that even for one of the simplest systems $H_2-O_2,$ as many as 20 different reactions can occur involving 8 different species \cite{T2000}. Moreover, under certain circumstances, all of these reactions can become reversible. Therefore, it is important to build a rigorous mathematical theory that does not impose any specific bounds on the size of initial data, range of temperatures, form of the pressure and the extent of reaction.

The  purpose of this paper is to present  the first complete existence result for such models, in case when the pressure depends on the mixture composition, and when no smallness assumption is postulated. Our particular concern is to perform a detailed construction leading to global in time weak solution.  

Recently, this area of mathematical modelling has attracted a lot of attention. Especially  the issue of so-called {\it multicomponent diffusion} has become much better understood, mostly due to several works devoted to the Maxwell--Stefan equations. These equations describe in an implicit way the relation between the diffusion velocities $\vc{V}_i$ and the gradients of the molar concentrations of the species $X_i$,
\begin{equation}\label{SM}
\underbrace{\Grad X_i-{\left(Y_i-X_i\right)\Grad \log \pi}}_{:=\vc{d_i}}=\sum_{{j=1} \atop {j\neq i}}^{n}\left(\frac{X_i X_j}{{\mathcal D}_{ij}}\right)\left(\vc{V}_j-\vc{V}_i\right).
\end{equation}
Here $Y_i$ is the mass concentration, $\pi$ is the internal pressure of the mixture and  ${\mathcal D}_{ij}> 0$  denotes the binary diffusion coefficient, ${\mathcal D}_{ij}={\mathcal D}_{ji}$. The first rigorous mathematical treatment of this system is due to Giovangigli \cite{Gio90, Gio91}, where various iterative methods for solving the linear system were presented as well as multicomponent diffusion coefficients for gaseous mixtures were provided. Later on, the proof of local in time existence result in the isobaric ($\pi=const$), isothermal case was proven by Bothe \cite{B2010}, see also \cite{HMPW13}.  Then, J\"ungel and Stelzer generalized this result and combined it with the  entropy dissipation method to prove the global in time existence of weak solutions \cite{JS13}, still in the case of constant pressure and temperature. Another interesting result devoted to the global-in-time existence of regular solutions for a special ternary gaseous mixture is presented in \cite{BGS12}. Recently, the coupling between the multi-species system and the incompressible system describing the fluid motion was investigated by Marion and Temam \cite{MT13} and Chen and J\"{u}ngel \cite{CJ13} . 

For the case of full systems describing the compressible mixtures, see \eqref{1.1} below, not much is known. The Navier-Stokes-Fourier system coupled with the set of reaction-diffusion equations for arbitrary large number of reacting species was treated by Feiresl, Petzeltov\'a and Trivisa \cite{FTP}. They proved the global-in-time existence of weak entropy variational solutions for  the Fick diffusion law and the state equation which does not depend on the species concentrations. Several extensions of this result and  asymptotic studies are also available, see \cite{DT1,DT2,CHT,KT,EZ}. 

Concerning general system with physically justified assumptions on the form of fluxes and transport coefficients the only global-in-time result, to the best of our knowledge, is due to Giovangigli \cite{VG}. His result is however restricted to data sufficiently close to the equilibrium state. The main intention of our paper is to present a possible extension of this result to the large data case, however under  less general assumptions.

It should be emphasized that the dependence of  the state equation on the species concentrations results in much more complex form of diffusion. The Fick approximation fails in providing a relevant description
 in this case, since it does not take into account strong cross-diffusion that is well known to play an important role \cite{BGS12, Wal62,T2000, BE}. From the mathematical point of view, it interferes with proving that the system possesses a Lyapunov functional, based on the notion of physical entropy. The integral form of entropy (in)equality is a source of majority of a-priori estimates being the corner stone of global-in-time analysis \cite{FN, JS13,MPZ1}.

Although including more general (multi-component) diffusion in the non-isobaric systems leads to a thermodynamically consistent model, it makes the analysis much more complex. The associated reaction-diffusion equations have to take into account the variation of total pressure, leading to a hyperbolic deviation. In other words, in the associated Maxwell--Stefan system \eqref{SM} the second term from the l.h.s. cannot be neglected. Note that this additional term is not in the conservative form.
To handle it one needs better regularity of the density than it follows from the hyperbolic continuity equation. Recently, Mucha, Pokorn\'y and Zatorska \cite{MPZ1} studied the system of $n$ reaction-diffusion equations for the chemically reacting species 
 \begin{equation*}\label{a}
		\pt{\vr Y_k}+\Div (\vr Y_{k} \vu)+ \Div (\vf_k)  =  \vr\vw_{k},\quad k=1,...n,
\end{equation*}
where $\vf_k$ is the $k$-th species diffusion flux satisfying
$$\vf_{k}=-C_0(\vr,\vt)\sum_{l=1}^{n} {C}_{kl}\vc{d}_l,$$
with the diffusion deriving forces $\vc{d}_l$ defined as in \eqref{SM}. They proved the global-in-time existence of weak solutions provided the additional regularity information about the total density $\vr$ and velocity $\vu$ is available.
This result holds for a particular choice of matrix $C$, which corresponds to $m_i m_j {\cal D}_{ij}=const$, in \eqref{SM}. However, the main argument presented there does not relay on this assumption. Extension of this result to the case of more general diffusion matrixes $C$, (see \cite{VG}, Chapter 7) is the work in progress.

Such a result creates a possibility of coupling the reaction-diffusion system to the model of fluid-motion, provided the assumptions on the regularity 
of $\vr$ and $\vu$ can be fulfilled. It is well known, that in the usual case of Navier-Stokes(-Fourier) equations with constant viscosity coefficients \cite{PLL, FN, NS} such regularity is not expected, except some special situations, f.i. 1D domains \cite{HS01,LM,M1D}. 

 Similar problem appears in the works devoted to analysis of the continuum mechanics mixture model derived in \cite{RT}. This model assigns different velocity fields (and temperatures) to each of the species. However, usually, a part of the interaction term (the momentum source) associated with the difference of gradients of species densities $\vr_k=\vr Y_k$ is neglected \cite{FGM05a, FGM05b}. This simplification leads to the family of homogeneous (interpenetrating) compressible multi-fluids models, for which a relevant existence analysis was recently performed by Kucher, Mamontov and Prokudin \cite{KMP12}.
 
In the present paper, this problem is solved by assuming that the viscosity coefficients are density-dependent functions. They are subject to a condition introduced by Bresch and Desjardins for the Saint-Venant system \cite{BD2,BD3}. In \cite{MV07}, Mellet and Vasseur proved that the weak solutions to the barotropic Navier-Stokes equations
\begin{equation*}
\left.
\begin{array}{r}
\vspace{0.2cm}
\pt\vr+\Div (\vr \vu) = 0\\
\vspace{0.2cm}
\ptb{\vr\vu}+\Div (\vr \vu \otimes \vu) - \Div (2\mu(\vr) \D(\vu)) -\Grad(\nu(\vr)\Div\vu)+ \Grad \vr^\gamma =\vc{0}\\
\end{array}\right\}\quad\mbox{{in}}\ (0,T)\times \Omega,
\end{equation*}
with
\eq{\nu(\vr)=2\vr\mu'(\vr)-2\mu(\vr),\label{BD-MV}}
are weakly sequentially stable  in the periodic domain $\Omega=\mathbb{T}^3$ and in the whole space $\Omega=\mathbb{R}^N$, $N=2,3$. 
For possible extension of this result to the case of heat-conducting flow and other boundary conditions we refer to \cite{BD, BDG}. The analogous result for the model of 2-component mixture was proven by Zatorska \cite{EZ2}.
It turns out, however, that construction of an approximate solution to such kind of systems is still an open problem, some special cases were studied in {\cite{BDC, MV08}}. The main problem here is the possibility of appearance of vacuum regions. Indeed, even though relation \eqref{BD-MV} provides particular structure necessary to improve the regularity of $\vr$, the uniform bound from below for the density is only known to be valid in 1D \cite{MV08}.
As a consequence, one may have a problem with defining the velocity vector field $\vu$ (the concentrations of species and the temperature). So, in contrary to the Navier-Stokes equations with $\mu,\nu=const$ \cite{PLL}, where the strong convergence of the density does not follow from the a-priori estimates,  here the biggest problem is to prove the compactness of the sequences approximating the velocity.

In case of isothermal model of two-species mixture with multicomponent diffusion, this problem was solved by including in the state equation the stabilizing term in the form of {\it cold pressure} \cite{EZ2, MPZ}. This singular pressure prevents the appearance of vacuum on the sets of non-zero measure, {as proven in  \cite{BD} for single-component heat-conducting fluids}. But in case of heat-conductiong mixtures, this term may not be sufficient.  Roughly speaking, the problem lies in coupling all the components of the system due to very strong nonlinearity in the energy equation and comparatively low regularity of solutions to all "sub-systems". Even under very special assumptions, the notion of weak solution introduced in \cite{EZ3} for which the weak sequential stability was established does not allow to construct  a desired approximation, at least not immediately. Indeed, the transfer of energy due to species molecular diffusion and more complex form of the entropy  makes it more difficult to construct the solution within the variational entropy formulation than within the usual energy formulation.  Note that it is unlike the case of single gas flow modelled by the Navier-Stokes-Fourier system \cite{NP1,NP2,MP3}. There, it is the passage to the limit in the total energy balance which requires more restrictive assumptions and thus makes the weak energy solutions harder to get.

Here we focus on presenting the detailed approximation scheme for the full system with several assumptions on the constitutive relation that will be specified in the next section.  Moreover, we rather concentrate on the model inside the domain than on realistic modelling of  processes at the boundary. Hence, we choose the most simple case of boundary conditions, and the particular form of thermodynamic functions, which maybe make our paper less general, but easier to follow.

To summarize,
we assume that $\Omega$ is a periodic box in $\R^3$, i.e. $\Omega=\mathbb{T}^3$, and we consider the following system of equations:
        \begin{equation}\label{1.1}
     \left.  
      \begin{array}{r}
            \vspace{0.2cm}
        \pt\vr+\Div (\vr \vu) = 0\\
            \vspace{0.2cm}
        \ptb{\vr\vu}+\Div (\vr \vu \otimes \vu) - \Div \vS+ \Grad \pi =\vc{0}\\
            \vspace{0.2cm}
        \pt{E}+\Div (E\vu )+\Div(\pi\vu) +\Div{\vQ}-\Div (\vS\vu)=0\\
            \vspace{0.2cm}
        \pt{\vr_k}+\Div (\vr_{k} \vu)+ \Div (\vf_k)  =  \vr\vt\vw_{k},\quad k\in\{1,...,n\}
            \end{array}\right\}\quad\mbox{in}\ (0,T)\times\Omega
    \end{equation}
with the corresponding set of initial conditions. Above, $\vr$ denotes the density of the mixture, $\vu$ is the mean velocity of the mixture, $\vS$ is the stress tensor, $\pi$ the total pressure, $E$ the total energy, 
$\vQ$ the heat flux, $\vr_k$ the density of the $k$-th constituent, $\vf_k$ the multicomponent diffusion matrix, $\vt$ is the temperature of the mixture and $\omega_k$ the chemical source term.
The system is supplemented by initial data on density  -- $\vr^0$, momentum -- ${\bf m}^0=\vr^0 \vu^0$, temperature -- $\vt^0$ and densities of species -- $\vr_k^0$.

Recall that the first equation, usually called the continuity equation,  describes the balance of the mass, the second equation expresses the balance of the momentum and the third one the balance of the 
total energy. The last set of $n$ equations describes the balance of separate constituents. Note that the system of equations cannot be independent, the last $n$ equations must sum into the continuity 
equation \eqref{1.1}$_1$. Thus, here we meet a serious  mathematical obstacle,  system $(\ref{1.1})_4$ is  degenerate parabolic in terms of $\vr_k$.


The second equation is the balance of momentum, in which the temperature and the species concentrations appear only in the form of the internal pressure $\pi$. The relation between the density dependent viscosities appearing in the form of the stress tensor $\vS$, enables to prove better regularity of the density.

The third equation of the above system may be rewritten as the internal energy balance, since the total energy $E$ is a sum of kinetic and internal energies:
	\begin{equation*}
	E=\frac{1}{2}\vr|\vu|^{2}+\vr e.
	\end{equation*}
The kinetic energy balance is nothing but a consequence of the momentum equation multiplied by $\vc{u}$ and integrated over $\Omega$, thus we may write the balance of the internal energy in the form   	
        \begin{equation}\label{1.1b}
        \ptb{\vr e}+\Div (\vr e\vu ) +\Div{\vQ}+\pi\Div\vu-\vS:\Grad\vu=0.
            \vspace{0.2cm}
    \end{equation}
However, the balance of internal energy together with the momentum equation is equivalent to the balance of the total energy and the momentum just in case, when the solutions are sufficiently regular. 
Hence it is true for classical or strong solutions, but it might be not true for weak solution which we introduce below.

%

The outline of the paper is the following. In Section \ref{CR} we discuss the constitutive relations and their consequences. We also introduce a notion of a weak solution and state our main result -- Theorem \ref{main}. The rest of the paper is devoted to its proof, it consists of several levels of approximation and the subsequent limit passages.
Section \ref{Sec:app} provides a description of two most important levels of approximation. First we present the system including only the main regularizing terms, which is marked by a presence of three main approximation parameters -- $\ep$, $\lambda$ and $\delta$. The existence result for this system is stated in Theorem \ref{t 2}. Then we rewrite the approximate system so that the momentum equation is replaced by its Faedo--Galerkin approximation (denoted by $N$), the total energy equation is replaced by the approximate thermal energy equation and the finite dimensional projections of the species equations with new entropy variables are introduced. Finally, in Section \ref{Sec:Basic} the basic level of approximation is considered with all the necessary regularizations needed to prove the existence of regular solutions as stated in Theorem \ref{th1}.
After this, in Section  \ref{Sec:Gal}, we let $N\to \infty$ in the equations of approximate system proving Theorem \ref{t 2} and come back to the first weak formulation from Section \ref{Sec:app}. Then, in Section \ref{Sec:delta} we perform the limit passage $\delta\to 0$, so that we are left with only two parameters of approximation: $\ep$ and $\lambda$. At this level, derivation of the B-D estimate becomes possible, we present it in Section \ref{Sec:BD}. In Section \ref{Sec:ep_lambda}
 we first present the new uniform bounds arising from the B-D estimate and then we let the last two approximation parameters to $0$, which finishes the proof of Theorem \ref{main}.

\section{Constitutive relations and their consequences. Main result} \label{CR}

Before introducing the definition of the weak solution and stating the main result of this paper, we have to specify the constitutive relations that we are going to use in our paper. We try to use such relations which are close to models of the real processes; however, in some cases we have to simplify them in order to be able to prove our main result.

\subsection{Pressure and internal energy} \label{PINE}

In the above system we use
\begin{equation} \label{pre}
\pi=\pi_c(\vr)+\frac{\beta}{3}\vt^4+\pi_m(\vr_k,\vt), \quad \beta>0,
\end{equation}
where the latter denotes the internal pressure of the mixture which is determined through the Boyle law
\begin{equation} \label{intpre} 
\pi_m(\vt,\vr_{1},\ldots,\vr_{n})=\sum_{k=1}^{n}p_{k}(\vt,\vr_{k})=\sum_{k=1}^{n}\frac{\vt\vr_{k}}{m_{k}};
\end{equation}
above, $m_{k}$ is the molar mass of the species $k$ and for simplicity, we set the gaseous constant  equal to 1.
We further assume that $\pi_c$ is a continuously differentiable function on $(0,\infty)$ satisfying the following growth conditions
	\begin{equation}\label{coldp}
	\pi'_{c}(\vr)=\left\{
	\begin{array}{cl}
	c_1\vr^{-\gamma^- -1}&\mbox{for} \ \vr\leq 1,\\
	c_2\vr^{\gamma^+-1}&\mbox{for} \ \vr>1,
	\end{array}\right.
	\end{equation}
for positive constants $c_1,c_2$ and $\gamma^->5$, $\gamma^+>3$ specified below. Note that the lower bounds for $\gamma^+$ and $\gamma^-$ are rather mathematical than physical. We shall keep in mind that interpretation of
the form of $\pi_c(\cdot)$ for $\vr \leq 1$ implies that vacuum regions of the fluid are not admissible, {as observed in \cite{BD}}.  The term $\frac{\beta}{3} \vt^4$ models the radiative pressure.


The internal energy can be expressed as $ e=e_m+\beta\frac{\vt^4}{\vr}+e_c,$  where
	$$\vr e_m(\vt,\vr_k)= \sumkN\vr_ke_k=\sumkN \vr_k(e^{st}_{k}+(\vt-\vt^{st}) c_{vk}),$$
	$$\vr^2\frac{de_c(\vr)}{d\vr}=\pi_c(\vr).$$

It is convenient to define the internal energy $e_k^0$ of the $k$-th species at zero temperature
	$$e_k^0=e_k^{st}-\vt^{st}c_{vk},$$
where $e_k^{st}$ denotes the formation energy of the $k$-th species at the positive standard temperature $\vt^{st}$, see Giovangigli \cite{VG} (2.3.4). Here we take without loss of generality $\vt^{st}=1$ and  assume that  this energy is equal for all species, for simplicity $e_k^0=0$. 
In the above formulas $c_{vk}$ denotes a constant-volume specific heat for the $k$-th species and it is related to the constant-pressure  specific heat ($c_{pk}$) by the formula
	\begin{equation}\label{cpcv}
	c_{pk}=c_{vk}+\frac{1}{m_k}.
	\end{equation}
For the sake of simplicity we take
	$$c_{vk}=c_v=1.$$
Hence, since $\sumkN\vr_k=\vr$, the molecular part of internal energy can be reduced to
	\begin{equation*}
	\vr e_m=\vr\vt.
	\end{equation*}
This simplification leads to the following reduction in the form of specific enthalpy of the $k$-th species with respect to the general form (2.3.6) in Giovangigli \cite{VG}
	$$h_k=e_k+\frac{\vt}{m_k}=c_{pk}\vt.$$

\subsection{Flux diffusion matrix} \label{FDM}

A key element of the presented model is the structure of laws governing chemical reactions. We first define the flux diffusion.
We consider the following special case
	\begin{equation}\label{eq:diff}
	\vf_{k}=-\vC_{0}\sum_{l=1}^{n} \vC_{kl}\vc{d}_l, \quad k=1,...n,
	\end{equation}
where $\vC_{0}, \ \vC_{kl}$ are multicomponent flux diffusion coefficients and $\vd=(d_{k}^{1},d_{k}^{2},d_{k}^{3})$ is the species $k$ diffusion force 
	\begin{equation}\label{eq:}
	d_{k}^{i}=\Grad_{x_{i}}\left({p_{k}\over \pi_m}\right)+\left({p_{k}\over \pi_m}-{\vr_{k}\over \vr}\right)\Grad_{x_{i}} \log{\pi_m}.
	\end{equation}
To fix the idea, we shall concentrate on the following explicit form of $\vC$
	\begin{equation}\label{Cform}
	\vC =\left(
		\begin{array}{cccc}
			Z_{1} & -Y_{1} & \ldots & -Y_{1}\\
			-Y_{2} & Z_{2}  & \ldots & -Y_{2}\\
			\vdots & \vdots & \ddots & \vdots \\
			-Y_{n} & -Y_{n}& \ldots & Z_{n}
		\end{array} \right),
	\end{equation}
where $Y_k=\frac{\vr_k}{\sumkN{\vr_k}}$ and $Z_{k}=\sum_{{i=1} \atop {i\neq k}}^{n} Y_{i}$. We also assume that $\vC_0=\vC_0(\vr,\vt)$ is a continuous function in both variables such that
	\begin{equation}\label{assC}
	\vC_0(\vr,\vt)=c_0(\vr )\tilde C_0(\vt),\quad \Un{\vC_0}\vr(1+\vt)\leq\vC_0(\vr,\vt)\leq\Ov{\vC_0}\vr(1+\vt),
	\end{equation}
for some positive constants $\Un{\vC_0},\Ov{\vC_0}$. Matrix (\ref{Cform}) can be examined in more general form. However. it is fixed to reduce a number of technical computations, which would 
make our proof more difficult to follow. \\

	\begin{rmk}\label{rem1}
	Using expressions for the diffusion forces \eqref{eq:} and the properties of $\vC$ one can rewrite \eqref{eq:diff} into the following form
		\begin{equation}\label{difp}
		\vf_{k}=-\frac{C_0}{\pi_m}\left(\Grad p_{k}-Y_{k}\Grad \pi_m\right)=-\frac{C_0}{\pi_m}\sum_{l=1}^{n}\vC_{kl}\Grad p_{l},
		\end{equation}
	moreover, $\sumkN\vf_k=\vc{0}$, pointwisely.
	\end{rmk}
	\begin{rmk}
	The matrix $D_{kl}=\frac{C_{kl}}{Y_k}$ is symmetric and positive semi-definite.
	\end{rmk}

\subsection{Heat flux} \label{HF}

Next, we look at the energy equation $(\ref{1.1})_3$.
The heat flux is given by the Fourier law
	\begin{equation}\label{cotoQ}
        \vQ=\sumkN h_k \vf_{k}-\kappa\Grad\vt,
    	\end{equation}
    where $h_k$ stands for partial enthalpies $h_k(\vt)=c_{pk}\vt$
and the thermal conductivity coefficient $\kappa=\kappa(\vr,\vt)=k(\vr)\tilde\kappa(\vt)$ 
is a smooth function such that
	\begin{equation}\label{cotok}
	\kappa(\vr,\vt)={\kappa_0}+\vr+\vr\vt^2+\beta\vt^B,
	\end{equation}
where  ${\kappa_0}=const.>0$, $B\geq8$. Again the last limitation is a consequence of mathematical needs.

\subsection{Stress tensor (viscous part)} \label{ST}

The viscous part of the stress tensor obeys the {\it Newton rheological law}
	\begin{equation}\label{chF:Stokes}
	\vS(\vr,\vu)= 2\mu(\vr)\D(\vu)+\nu(\vr)\Div \vu{\bf{I}},
	\end{equation}
where $\D(\vu)=\frac{1}{2}\left(\Grad \vu+(\Grad \vu)^{T}\right)$ and the nonnegative viscosity coefficients $\mu(\vr)$, $\nu(\vr)$ satisfy the Bresch-Desjardins relation
	$$\nu(\vr)=2\mu'(\vr)-2\mu(\vr).$$
In this work we consider only the simplest example of such functions, namely $\mu=\vr$, $\nu=0$.

\subsection {Species production rates} \label{SPR}

We assume that the species production rates are Lipschitz continuous of $\vr_1,\ldots,\vr_n$ and that there exist positive constants
 $\Un{\vw}$ and $\Ov{\vw}$ such that
	\begin{equation}\label{wform}
	-\Un{\vw}\leq\vw_{k}(\vr_{1},\ldots,\vr_{n})\leq\Ov{\vw},\quad  \quad \mbox{ for\ all}\quad  0\leq Y_{k}\leq 1,\quad k=1,\ldots,n;
	\end{equation}
moreover, we suppose that
	\begin{equation}\label{wform1}
	\vw_{k}(\vr_{1},\ldots,\vr_{n})\geq 0\quad\mbox{ whenever}\ \ Y_{k}=0.
	\end{equation}
We also anticipate the mass constraint between the chemical source terms
	\begin{equation}\label{wform0}
	\sum_{k=1}^{n}\vw_{k}=0.
	\end{equation}
Another restriction that we postulate for chemical sources is dictated by the second law of thermodynamics; it asserts that the entropy production associated with any admissible chemical reaction is nonnegative. In particular, $\vw_{k}$ must enjoy the following  condition
	\begin{equation}\label{admiss}
	-{\vr\sumkN{g_k\vw_k}}\geq 0,
	\end{equation}
where $g_k$ are the Gibbs functions specified below.

\subsection{Entropy}

In accordance with the second law of thermodynamics we postulate the existence of a state function called the entropy. It is defined (up to an additive constant) in terms of differentials of energy, total density, and species mass fractions by the Gibbs relation:
	\begin{equation}\label{entropdiff}
	\vt \mathrm{D} s=\mathrm{D} e+\pi\mathrm{D}\left({1 \over \vr}\right)-\sumkN g_{k}\mathrm{D} \vy,
	\end{equation}
where $\mathrm{D}$ denotes the total derivative with respect to the state variables $\vt,\vr,Y_1,\ldots,Y_n$. The Gibbs function of the $k$-th species is defined by the formula
	$$g_k=h_k-\vt s_k,$$
	$$g_k=h_k-\vt s_k=c_{pk}\vt-\vt\log\vt+\frac{\vt}{m_k}\log{ {\vr_{k}}\over{m_k}}.$$
 Here
$s_k=s_k(\vt,\vr_{k})$ are the specific entropies of species and their general form (relation (2.3.20) from \cite{VG}, Giovangigli) is the following
	$$s_k=s_k^{st}+c_{vk}(\log\vt-\log\vt^{st})-\frac{1}{m_k}\log\frac{\vt^{st}\vr_k}{p^{st} m_k}+\frac{4\beta\vt^3}{3\vr_k}.$$

In comparison to the general framework we reduce the form of entropy  by an assumption that $s_k^0=s_k^{st}-c_{vk}\log\vt^{st}=0$, for all $k$. Moreover, we set the standard  pressure $p^{st}$ equal to 1.
Therefore
	\begin{equation}\label{entropycv}
	s_k=\log{\vt}-{1\over m_{k}}\log{ {\vr_{k}}\over{m_k}}+\frac{4\beta\vt^3}{3\vr_k}
	\end{equation}
and
	\begin{equation}\label{entropycv1}
	\vr s=\sumkN \vr_k s_{k}=\vr\log{\vt}-\sumkN{\frac{\vr_k}{m_{k}}}\log{ {\vr_{k}}\over{m_k}}+\frac{4a}{3}\vt^3,
	\end{equation}
and this expression will be understood as our definition of the entropy.

  The evolution of the total entropy of the mixture 
can be  described by the following equation
	\begin{equation}\label{entropy}
	\ptb{\vr s}+\Div(\vr s\vu)+\Div\left( \frac{\vQ}{\vt}+\sumkN \frac{g_{k}}{\vt}\vf_{k}\right)=\sigma,
	\end{equation}
where $\sigma$ is the entropy production rate
	\begin{multline} \label{sigma}
	\sigma=\frac{\vS:\Grad\vu}{\vt}-{\vQ\cdot\Grad\vt\over{\vt^{2}}}-\sumkN\vf_{k}\cdot\Grad\left(\frac{g_{k}}{\vt}\right)-{{\sumkN g_{k}\vw_{k}}}\\
	=\frac{2\vr\D\vu:\D\vu}{\vt}+{\kappa\left|\Grad\vt\right|^2\over{\vt^{2}}}-\sumkN\frac{\vf_{k}}{m_k}\cdot\Grad\left(\log p_k\right)-{{\sumkN \vr g_{k}\vw_{k}}}\geq 0.
	\end{multline}
The equation  follows from the internal energy balance \eqref{1.1b} and the assumptions we posed above.

	\begin{rmk}
	Note that the above manipulations work only for smooth solutions and if we know that $\sumkN\vr_k=\vr$, otherwise the Gibbs formula and the Gibbs functions must be modified.
	\end{rmk}

\subsection{Weak formulation and the main result}
In what follows we introduce the notion of the weak solution to system \eqref{1.1}. It is a set of space-periodic functions $(\vr,\vu,\vt,\{\vr_k\}_{k=1}^n)$ such that
$\vr >0$, $\vt >0$, $\vr_k \geq 0$, $\vr = \sum_{k=1}^n \vr_k$ a.e. in $(0,T)\times \Omega$,
\begin{equation} \label{reg_sol}
\begin{array}{c}
\displaystyle \vr \in L^\infty(0,T;L^{\gamma^+}(\Omega)), {\vr}^{-1} \in L^\infty(0,T;L^{\gamma^-}(\Omega)) \\
\displaystyle \sqrt{\vr}\vu \in L^\infty(0,T;L^{2}(\Omega)), \sqrt{\vr} \nabla \vu \in L^2(0,T; L^{2}(\Omega)) \\
\displaystyle \vt \in L^\infty(0,T;L^4(\Omega)), \vt \in L^2(0,T;W^{1,2}(\Omega)) \\
\displaystyle \vt^{\frac B2} \in L^2(0,T;W^{1,2}(\Omega)) \\
\displaystyle \sqrt{\vr_k} \in L^2(0,T;W^{1,2}(\Omega)),
\end{array}
\end{equation}  
and the following identities are fulfilled:
\begin{itemize}
\item[-]
the continuity equation
	\begin{equation*}\label{WS1}
	\pt\vr+\Div (\vr \vu)= 0,
	\end{equation*}
is satisfied pointwisely on $[0,T]\times\Omega$;

\item[-]
the momentum equation
	\begin{equation}\label{WS2}
	\begin{split}
	-&\intTO{\vr\vu\cdot\pt\vcg{\phi}}-\intO{\vc{m}^{0}\cdot\vcg{\phi}(0)}\\
	-&\intTO{(\vr{\vu}\otimes{\vu}):\Grad \vcg{\phi}}
	+\intTO{2\vr\D\vu:\D \vcg{\phi}}
	-\intTO{\pi\Div \vcg{\phi}}=0
	\end{split}
	\end{equation}
holds for any test function $\vcg{\phi}$ smooth function such that $\vcg{\phi}(\cdot,T)=\vc{0}$.

\item[-]
the species equations
	\begin{multline*}\label{WS3}
 	 \intTO{\vr_k\pt\phi}+\intO{\vr_k^0\phi(0)}\\+\intTO{ \vr_{k} \vu\cdot\Grad\phi}+\intTO{\vf_k\cdot\Grad\phi}  =  -\intTO{\vr\vt\vw_{k}\phi},\quad k\in\{1,...,n\}
	\end{multline*}
are fulfilled for any smooth function $\phi$ such that $\phi(\cdot,T)=0$;

\item[-]
the total energy equation
	\begin{equation}\label{WS4}
	\begin{split}
	&\intTO{\lr{\vr e+\frac{1}{2}\vr|\vu|^2}\pt\phi} +\intO{\lr{\vr e+\frac{1}{2}\vr|\vu|^2}(0)\phi(0)}\\
	&+\intTO{\lr{\vr e\vu+\frac{1}{2}\vr|\vu|^2\vu}\cdot\Grad\phi}-\intTO{\kappa\Grad\vt\cdot\Grad\phi}+\sumkN \intTO{ \vt\frac{\vf_{k}}{m_k}\cdot\Grad\phi}\\
	&+\intTO{\pi\vu\cdot\Grad\phi}-\intTO{(2\vr\D(\vu)\vu)\cdot\D\phi}= 0
	\end{split}
	\end{equation}
holds for any smooth function $\phi$ such that $\phi(\cdot,T)=0$, where the heat flux term is to be understood as in the distributional sense.
\end{itemize}


The main result of this paper reads
\begin{thm}\label{main}
Let $\Omega = {\mathbb T}^3$ be a periodic box, $\gamma^+ >3$, $\gamma^- >5$, $\gamma^- > \frac{5\gamma^+ -3}{\gamma^+ -3}$, $B\geq 8$, $\mu(\vr) = \vr$, $\nu(\vr)=0$. Let $\vr^0 \in L^{5\gamma^+/3}(\Omega)$, $1/\vr^0 \in L^{(5\gamma^- -1)/3}(\Omega)$, $\vc{m}^0 \in L^{1}(\Omega)$ such that $\frac{(\vc{m}^0)^2}{\vr^0} \in L^1(\Omega)$, $\vt^0 \in L^4(\Omega)$, $\vr_k^0 \in L^1(\Omega)$. Let $T>0$ be arbitrary. Then there exists a weak solution to \eqref{1.1} in the weak sense specified above. Moreover, the density $\vr>0$ and the temperature $\vt>0$ a.e. in $(0,T)\times \Omega$.  
\end{thm}

\section{Approximation}\label{Sec:app}

The aim of this section is to present two levels of approximation (in fact, there will be more of them as intermediate steps, however, we will not mention all of them explicitly). It is well-known that one of the main  problems with the so-called Bresch-Desjardins relation is the question of a good approximation. This has to allow to construct sequence of solutions which are compatible with the better estimate; e.g. for the isentropic Navier--Stokes system (i.e. $p(\vr) \sim \vr^\gamma$) the problem is totally open. In case of  no vacuum regions (as it is here), there is a chance to show existence of such approximations, {see f. i.  \cite{BDC}}, however, the problem is far from being trivial (in fact, it is more complex than the problem for the full Navier--Stokes--Fourier system presented in \cite{FN}).

First, we  take $\ep$, $\delta$ and $\lambda >0$ and fix $s$ a sufficiently large positive integer. 
{Our aim is to consider the regularized problem given below, in which $\ep$ is the rate of dissipation in the continuity equation, $\delta$  introduces additional relaxation into the species mass balance equations. By $\lambda$  we insert  to the momentum equation the artificial smoothing operator $\lambda\vr\Grad\lap^{2s+1}\vr$ with $s$ sufficiently large, inspired by the works of Bresch and Desjardins \cite{BDC, BD3}, we also introduce another regularization of the momentum $\lambda\vr\lap^{2s+1}(\vr\vu)$. }
Note that at the end, after passing subsequently with $\delta$, $\ep$ and $\lambda \to 0^+$, we recover our original problem \eqref{1.1}.

We look for space periodic functions $(\vr, \vu, \vt, \{\vr_k\}_{k=1}^n)$ such that
\begin{equation} \label{RG_Lev1}
\begin{array}{c}
\vr \in L^2(0,T; W^{2s+2}(\Omega)), \partial_t \vr \in L^2(0,T;L^2(\Omega)) \\
\vu \in L^2(0,T; W^{2s+1}(\Omega)) \\
\vt \in L^2(0,T; W^{1}(\Omega)) \cap L^B(0,T; L^{3B}(\Omega))\\
\vr_k \in L^q(0,T; L^q(\Omega)), \quad  q>\frac 53
\end{array}
\end{equation} 
solving the following problem:

\begin{itemize}

\item the approximate continuity equation
	\begin{equation}\label{apcontN}
	\begin{array}{c}
	\pt\vr+\Div (\vr \vu) -\ep\lap\vr= 0 \\
	\vr(0,x) = \vr_{\lambda}^0(x)
	\end{array}
	\end{equation}
is satisfied pointwisely on $[0,T]\times\Omega$ and the initial condition holds in the strong $L^2$ sense; here $\vr_{\lambda}^0 \in  C^\infty(\overline{\Omega})$ is a regularized initial condition such that $\vr_{\lambda}^0 \to \vr^0$ in $L^{\gamma^+}(\Omega)$ for $\lambda \to 0^+$, such that $\lambda \|\nabla^{2s+1}\vr^0_\lambda\|_{2}^2 \to 0$ for $\lambda\to 0^+$, and 
\begin{equation}\label{apcontreg}
	\quad \inf_{x\in\Omega}\vr^{0}_{\lambda}(x)>0;
\end{equation}	

\item the weak formulation of the approximate momentum equation
	\begin{equation}\label{FG1}
	\begin{split}
	&\intTO{\vr\vu\cdot\pt\vcg{\phi}}
	-\intTO{\lambda\lap^{s}\Grad(\vr\vu) :\lap^{s}\Grad (\vr\vcg{\phi})}\\
	+&\intTO{(\vr{\vu}\otimes{\vu}):\Grad \vcg{\phi}}
	-\intTO{2\vr^n\D\vu:\D \vcg{\phi}}
	+\intTO{\pi\Div \vcg{\phi}}\\
	-&\lambda\intTO{\lap^s\Div\left({\vr\vcg{\phi}}\right)\lap^{s+1}\vr}-\ep\intTO{(\Grad\vr\cdot\Grad)\vu\cdot \vcg{\phi}}=-\intO{\vc{m}^{0}\cdot\vcg{\phi}(0)}
	\end{split}
	\end{equation}
holds for any test function $\vcg{\phi}\in L^2(0,T;W^{2s+1}(\Omega))\cap W^{1,2}(0,T;W^{1,2}(\Omega))$ such that $\vcg{\phi}(\cdot,T)=\vc{0}$;

\item the weak formulation of the total energy equality
	\begin{equation}\label{totEN}
	\begin{split}
	&\intTO{\lr{\vr e+\frac{1}{2}\vr|\vu|^2+\frac{\lambda}{2}|\Grad^{2s+1}\vr|^2}\pt\phi} +\intTO{\lr{\vr e\vu +\frac{1}{2}\vr|\vu|^2\vu}\cdot\Grad\phi}\\
	-&\intTO{\kappa_{\ep}\Grad\vt\cdot\Grad\phi}+\sumkN \intTO{ \left(\vt\frac{\vf_{k}}{m_k}-\delta \vt\frac{\Grad \vr_{k}}{\vr_k}-\ep\vt\frac{\Grad \vr_k}{m_k} \right)\cdot\Grad\phi}\\
	+&\intTO{\pi\vu\cdot\Grad\phi}
	-\intTO{(2\vr^n\D(\vu)\vu)\cdot\Grad\phi}=-\intTO{\lr{\frac{\ep}{\vt^2}-\ep\vt^{5}}\phi}\\
	+&\intTO{R_{\ep,\lambda}(\vr,\vt,\vu,\phi)}-\intO{\lr{\vr_\lambda^0 e^0+\frac{1}{2}\vr^0_\lambda|\vu_\lambda^0|^2+\frac{1}{2}|\Grad^{2s+1}\vr^0_\lambda|^2}\phi(0)},
	\end{split}
	\end{equation}
with 
\begin{equation}\label{weak_kappa}
\intTO{\kappa_{\ep}\Grad\vt\cdot\Grad\phi}=\intTO{\lr{\kappa_0+\ep\vr^n+\vr+\vr\vt^2}\Grad\vt\cdot\Grad\phi}+\frac{\beta}{B+1}\intTO{\vt^{B+1}\lap\phi}
\end{equation}
and
\begin{equation}\label{RN} 
\begin{split}
{R_{\ep,\lambda}(\vr,\vt,\vu,\phi)} = & \lambda\left[\Delta^s (\Div (\vr \vu \phi)) \Delta^{s+1} \vr 
-\lap^s\Div(\vr\vu)\lap^{s+1}\vr\phi\right]\\
&-\lambda\lap^s\Div(\vr\vu)\Grad\lap^s\vr\cdot\Grad\phi
- \lambda\left[|\Delta^s (\nabla (\vr \vu))|^2 \phi- \Delta^s \nabla (\vr \vu): \Delta^s \nabla (\vr \vu \phi)\right] \\
&+ 
\lambda \ep\lap^{s+1}\vr\Grad \lap^{s}\vr \cdot\Grad\phi + \frac{\ep}{2} |\vu|^2 \nabla \vr \cdot \nabla\phi +\ep \nabla \vr \cdot \nabla \phi \left(e_c(\vr) + \frac{\pi_c(\vr)}{\vr}\right)
\end{split}
	\end{equation}
is satisfied
for any $\phi \in C^\infty([0,T]\times \overline{\Omega})$ with $\phi(T,\cdot) = 0$;
here $e^0 = e(\vr^0_\lambda, \vt^0_\lambda)$, $\vu^0_\lambda = \frac{\vc{m}^0}{\vr^0_\lambda}$ and $\vt^0_\lambda \in C^\infty(\overline{\Omega})$, $\vt^0_\lambda \to \vt^0$ for $\lambda \to 0^+$ in $L^{4}(\Omega)$, and 
\begin{equation} \label{ic_temp}
\quad 0<\inf_{x\in\Omega}\vt^{0}_{\lambda}(x) = \Un{\vt^0}\leq\vt^0(x)\leq\sup_{x\in\Omega}\vt^{0}_{\lambda}(x)= \Ov{\vt^0} <\infty;
\end{equation}

\item the weak formulation of the species equations

	\begin{multline}\label{rk_del}
 	 \intTO{\lr{\delta \ln \left(\frac{\vr_k}{m_k}\right)+\frac{\vr_k}{m_k}}\pt\phi}+\intTO{ \frac{\vr_{k}}{m_k} \vu\cdot\Grad\phi} \\+\intTO{\frac{\vf_k}{m_k}\cdot\Grad\phi}  =  -\intTO{\frac{\vr\vt\vw_{k}}{m_k}\phi} - \intO{\lr{\delta \ln\left(\frac{\vr_k^0}{m_k}\right)  +\frac{\vr_k^0}{m_k}}\phi(0)}, \quad k\in\{1,...,n\}
	\end{multline}
is fulfilled for any function $\phi \in C^\infty([0,T];\Omega)$ such that $\phi(\cdot,T)=0$.

Above, we denoted
\begin{equation} \label{notation_x}
\begin{split}
\vr^n = \sum_{k=1}^n \vr_k, \quad
\kappa_\ep = \frac{\ep}{\Un{m}} \vr^n + \kappa (\vr,\vt),
\end{split}
\end{equation}
where $\Un{m}=\min\{m_1,\ldots,m_n\}$.
\end{itemize}
We prove the following result
\begin{thm} \label{t 2}
Under the assumptions of Theorem \ref{main} and the assumptions specified in this section, for any $T>0$, $\ep$, $\delta$ and $\lambda>0$, there exists a solution to problem \eqref{RG_Lev1}--\eqref{notation_x} in the sense defined above.
\end{thm} 

Indeed, the proof of this result is far from being obvious. To prove Theorem \ref{t 2} we have to introduce another level of approximation, based on regularization of certain quantities and finite dimensional projection (Faedo--Galerkin approximation) of the momentum equation and the species equation as well as on replacing the species densities by their logarithms. More precisely, we additionally take $n \in \mathbb{N}$ and look for functions 
$(\vr,\vu,\vt,\{r_k\}_{k=1}^n)$ such that (for definition of $X_N$ and $Y_N$ see below)
\begin{equation} \label{RG_Lev2}
\begin{array}{c}
\vr \in L^2(0,T; W^{2s+2}(\Omega)), \partial_t \vr \in L^2(0,T;L^2(\Omega)) \\
\vu \in C([0,T]; X_N) \\
\vt \in L^\infty((0,T)\times \Omega) \cap L^2(0,T; W^{1,2}(\Omega))\\
r_k \in C([0,T]; Y_N)
\end{array}
\end{equation} 
solving the following problem:

\begin{itemize}

\item the approximate continuity equation:
	\begin{equation}\label{apcontN_2}
	\begin{array}{c}
	\pt\vr+\Div (\vr \vu) -\ep\lap\vr= 0 \\
	\vr(0,x) = \vr_{\lambda}^0(x)
	\end{array}
	\end{equation}
it is satisfied pointwisely on $[0,T]\times\Omega$ and the initial condition holds in the strong $L^2$ sense; here $\vr^0_\lambda$ is as above;

\item the Faedo-Galerkin approximation for the weak formulation of the momentum balance: we look for $\vu\in C([0,T]; X_N)$ such that 	
	\begin{equation}\label{FG}
	\begin{split}
	&\intO{\vr\vu(t)\cdot\vcg{\phi}}-\intO{\vc{m}^{0}\cdot\vcg{\phi}}
	-\lambda\inttO{\vr\lap^{2s+1}(\vr\vu) \cdot\vcg{\phi}}\\
	&-\inttO{(\vr{\vu}\otimes{\vu}):\Grad \vcg{\phi}}
	+\inttO{\vS:\Grad \vcg{\phi}}
	-\inttO{\pi\Div \vcg{\phi}}\\
	&-\lambda\inttO{\vr\Grad\lap^{2s+1}\vr\cdot \vcg{\phi}}+\ep\inttO{(\Grad\vr\cdot\Grad)\vu\cdot \vcg{\phi}}=0
	\end{split}
	\end{equation}
is satisfied  for all  $t\in[0,T]$ for any test function $\vcg{\phi}\in X_{N}$, and $X_{N}=\operatorname{span}\{\vcg{\phi}_{i}\}_{i=1}^{N}$,   where $\{\vcg{\phi}_{i}\}_{i=1}^{\infty}$ is an orthonormal basis in $L^2(\Omega)$, such that $\vcg{\phi}_i\in C^\infty(\Omega)$ for all $i\in\mathbb{N}$;

\item the approximate thermal energy equation:
		 \begin{multline}\label{app_vt}
	\ptb{\vr \vt+\beta\vt^4} +\Div\lr{\vu\lr{\vr\vt+\beta\vt^4}}-\Div\left(\kappa_{\ep}\Grad\vt\right)
	+\sumkN \Div \left(\vt\frac{\vf_{k}}{m_k}-\delta \vt\Grad r_{k}-\ep\vt\erk\Grad r_k\right)= \\
	\frac{\ep}{\vt^2}-\ep\vt^5-\lr{\pi_m + \frac{\beta}{3} \vt^4}\Div\vu+2\vr^n|\D(\vu)|^2
	+\lambda|\lap^{s}\Grad(\vr\vu)|^{2}+\lambda\ep|\lap^{s+1}\vr|^{2}+\ep\frac{1}{\vr}\frac{\partial\pi_c(\vr)}{\partial\vr}|\Grad\vr|^2,
	\end{multline}
it is fulfilled pointwisely on $(0,T)\times \Omega$, the initial condition $\vt^0_{\lambda}$ is as above;

\item the Faedo-Galerkin approximation for the weak formulation of the species mass balance equations: we look for $r=\{r_k\}_{1}^n\in C([0,T];Y_N)^n$ such that 
	\begin{multline}\label{rk_gal}
	\begin{split}
	&\intO{(\delta r_{k}+\erk)(t)\phi}-\intO{(\delta r^0_{k}+{\rm e}^{r_k^0})\phi(0)}=\\
	&\inttO{\left(\erk\vu-(\delta+\ep\erk)\Grad r_{k}\right)\cdot\Grad\phi}
	-\inttO{\left({C_{0}}(\vr,\vt)\sumlN\hat{D}_{kl}(r)\Grad r_{l}\right)\cdot\Grad\phi}\\
	+&\inttO{\left({C_{0}(\vr,\vt)}\Grad\log \vt\sumlN\hat D_{kl}(r)\right)\cdot\Grad\phi}
	+\inttO{\frac{\vr^{n}\vt\vw_k}{m_k}\phi},
	\end{split}
	\end{multline}
is satisfied  for each $k=1,\ldots,n$ and for all  $t\in[0,T]$ for any test function ${\phi}\in Y_{N}$, and $Y_{N}=\operatorname{span}\{{\phi}_{i}\}_{i=1}^{N}$,   where $\{{\phi}_{i}\}_{i=1}^{\infty}$ is an orthonormal basis in $L^2(\Omega)$ such that ${\phi}_i\in C^\infty(\Omega)$ for all $i\in\mathbb{N}$ and the initial condition $r_k^0 = \ln \frac{\vr^0_k}{m_k}$.
\end{itemize}

 
In the system above  we introduced the following notation
	\begin{equation} \label{21a}
\vr^{n}=\sumkN m_k\erk,
	\end{equation}
and
	$$\vf_{k}=\vf_{k}(\vr,\vt,r),\quad {\vf_{k}}(\vr,\vt,r)=-{C_{0}}(\vr,\vt)m_k\left(\sumlN\hat{D}_{kl}(r)\Grad r_{l}+\Grad\log \vt\sumlN\hat{D}_{kl}(r)\right)$$
	$$\hat{D}_{kl}(r)=\frac{\vt C_{kl}(r)}{\pi_{m}\erk  m_k}\erk \erl,\quad\pi_{m}=\pi_m(\vt,r),\quad \pi_m(\vt,r)=\sumkN\vt\erk,$$
	$$C_{kl}(r)=C_{kl}(Y_{1},\ldots,Y_{n}), \quad \mbox{where} \ Y_{k}=\frac{m_k\erk}{\vr^{n}}$$
	which formally corresponds to the definition of $\vr_k$ above.
From the definition it follows that the matrix $\hat D_{kl}$ is symmetric, positive semi-definite and its $L^\infty$ norm is bounded by a constant dependent only on $m_1,\ldots, m_k$. In particular
	\begin{equation}\label{propD}
	0\leq \hat D_{kl}(r)\leq c.
	\end{equation}
Note that passing with $N \to \infty$ in \eqref{apcontN_2}--\eqref{rk_gal} and setting $r_k = \ln \frac{\vr_k}{m_k}$ we formally recover our previous approximate system.

	\begin{rmk}
		The role of parameter $\beta$ is significant in obtaining the weak formulation for the total energy balance. We need it to pass to the limit in $\vQ$ especially in the term of the form $\vr\vt^2\Grad\vt$. This in turn is needed to close the B-D estimate.
	\end{rmk}

\begin{thm}[Existence of regular solution]\label{th1}
	Let $N \in\mathbb{N}$, $\ep$, $\delta$ and $\lambda >0$. Let $\vc{m}^0$, $\vr^0_\lambda$, $\vt^0_\lambda$, $\vr_k^0$ be as above. Under the assumptions of Theorem \ref{main} and assumptions stated above in this section, for any $T>0$, there exists a solution of system \eqref{apcontN_2}--\eqref{rk_gal} in the sense specified above.
\end{thm}

\section{Basic level of approximation}\label{Sec:Basic}
This section is dedicated to the proof of Theorem \ref{th1}. The strategy of the proof can be summarized as follows:
\begin{itemize}
\item We fix $\vu(t,x)$ in the space $C([0,T];X_{{N}})$ and $r_{k}$ in $C([0,T];Y_{{N}})$ and use it to find a unique smooth solution to (\ref{apcontN_2}) $\vr=\vr(\vu)$.
\item We smoothen functions $r_k(t,x)$, $k=1,\ldots,n$ with respect to time by a convolution with smoothing kernel denoted by $\xi$ and we find $\vt(\vu,r_{k},r_{k,\xi})$ a unique solution to a regularized version of equation (\ref{app_vt}).
\item We find the unique local-in-time solution to the momentum equation and the species equations by a fixed point argument.
\item We extend the local-in-time solution for the whole time interval using  uniform estimates.
\item We pass to the limit  $\xi\to 0^+$ and thus prove Theorem \ref{th1}.
\end{itemize}

\subsection{Continuity equation}\label{ex:cont}
We first prove the existence of a smooth, unique solution to the approximate continuity equation
in the situation when the vector field $\vu(x,t)$ is given and belongs to $C([0,T];X_{{N}})$.

The following result can be proven by the Galerkin approximation and the well known statements about the regularity of linear parabolic systems (for the details of the proof see \cite{FN}, Lemma 3.1).
	\begin{lemma}\label{kropa}
	Let $\vu\in C([0,T];X_{ N})$ for $ N$ fixed and let $\vr^{0}_{\lambda}$ be as above.
	Then there exists the unique classical solution to (\ref{apcontN_2}), i.e. $\vr\in V^\vr_{[0,T]}$, where
		\begin{equation}\label{regvr}
		V^\vr_{[0,T]}=\left\{\begin{array}{c}
		\vr\in C\left([0,T];C^{2+\nu}(\Omega)\right),\\
		\pt\vr\in C\left([0,T];C^{0,\nu}(\Omega)\right).
		\end{array}\right\}
		\end{equation}
	Moreover, the mapping $\vu\mapsto\vr(\vu)$ maps bounded sets in $C([0,T];X_{ N})$ into bounded sets in $V^\vr_{[0,T]}$ and is continuous  with values in $C([0,T];C^{2+\nu'}(\Omega))$, $0<\nu'<\nu<1$,\\
		\begin{equation}\label{below0}
		\Un{\vr^{0}}e^{-\int_{0}^{\tau}{\|\Div\vu\|_{\infty}}\dt}\leq\vr(\tau,x)\leq\Ov{\vr^{0}}e^{\int_{0}^{\tau}{\|\Div\vu\|_{\infty}}\dt}\quad for\ all\ \tau\in[0,T],\ x\in\Omega.
		\end{equation}
Finally, for fixed $N \in \mathbb{N}$, the function $\vr$ is smooth ($C^\infty$) in the space variable.		
	\end{lemma}

\subsection{Temperature equation}
The existence of unique solution to \eqref{app_vt} can be proven as in \cite{FN} with the necessary modifications to accommodate the extra terms due to the dependence of the  species concentrations. First, however, we need to smoothen the coefficients of \eqref{app_vt} with respect to time. We will consider the following system
	\begin{multline}\label{app_vt1}
	\ptb{\vr \vt+\beta\vt^4} +\Div\lr{\vu\lr{\vr \vt+\beta\vt^4}}-\Div\left(\kappa_{\ep,\xi}\Grad\vt\right)
	+\sumkN \Div \left(\vt\frac{\vf_{k,\xi}}{m_k}-\delta \vt\Grad r_{k}-\ep\vt\erki{\xi}\Grad r_k\right)= \\
	\frac{\ep}{\vt^2}-\ep\vt^5-\lr{\pi_m + \frac{\beta}{3} \vt^4}\Div\vu+2\vr^n|\D(\vu)|^2
	+\lambda|\lap^{s}\Grad(\vr\vu)|^{2}+\lambda\ep|\lap^{s+1}\vr|^{2}+\ep\frac{1}{\vr}\frac{\partial\pi_c(\vr)}{\partial\vr}|\Grad\vr|^2.
	\end{multline}
Here $\xi$ denotes the convolution with the standard regularizing kernel $\omega_\xi$ applied at several places to the time variable of function $r(t,x)=\{r_1(t,x),\ldots, r_n(t,x)\}$ (extended by the initial value to  negative times); moreover we denoted
	\begin{equation*}
	\kappa_{\ep,\xi}=\kappa_{\ep}(\vr,\vt,r_\xi)
	\end{equation*}
and
	$$\vf_{k,\xi}=\vf_{k}(\vr,\vt,r,r_\xi)
	=-{C_{0}}(\vr,\vt)m_k\left(\sumlN\hat{D}_{kl}(r_\xi)\Grad r_{l}+\Grad\log \vt\sumlN\hat{D}_{kl}(r_\xi)\right).$$

We have the following lemma.

	\begin{lemma}\label{lemma_vt}
	Let $\vu\in C([0,T];X_{ N})$ be a given vector field, $r_{k}\in C([0,T];Y_N)$, $k=1,\ldots,n$ be given functions and let $\vr(\vu)$ be the unique solution of (\ref{apcontN_2}). 

	Then, equation (\ref{app_vt1}) with the initial condition $\vt^0_\lambda$ defined above  admits a unique strong solution $\vt=\vt(\vu,r)$ which belongs to
		\begin{equation}\label{regvt}
		V_{[0,T]}^\vt=\left\{\begin{array}{c}
			\pt\vt\in L^2((0,T)\times\Omega),\ \lap\vt\in L^2((0,T)\times\Omega),\\
			\vt\in L^\infty(0,T;W^{1,2}(\Omega)),\quad \vt,\vt^{-1}\in L^{\infty}((0,T)\times\Omega).
			\end{array}\right\}
		\end{equation}
	Moreover, the mapping $(\vu,r)\mapsto\vt(\vu,r)$ maps bounded sets in $C([0,T];X_{ N})\times C([0,T];Y_N)$ into bounded sets in $V^\vt_{[0,T]}$ and is continuous  with values in $L^2([0,T];W^{1,2}(\Omega))$.
	\end{lemma}

\pf The rough idea of the proof is to transform and to regularize equation \eqref{app_vt1} in such a way that the classical theory for quasilinear parabolic equations could be applied.
Lest us consider the following approximate equation 
	 \begin{equation}\label{eqK}
	 \begin{split}
	&\lr{\vr_\eta+\frac{4\beta\vt^4}{\sqrt{\vt^2+\eta^2}}} \pt\vt+\Div((\vu\vr)_\eta\vt+\beta\vu\theta_\eta^4)-\lap K_{\eta,\xi}\\
	+&\frac{\partial}{\partial x_i}\int_{1}^\vt \frac{\partial}{\partial x_i}\left([\kappa_{\ep}(\vr,s, r_\xi)]^\eta+\sum_{k,l=1}^n[\hat D_{kl}(r_\xi)C_0(\vr,s)]^\eta\right)~ds\\
	-&\sumkN \Div \left({\langle{C_0\rangle}_\eta(\vr_\eta,\theta_\eta)}\theta_\eta\sumlN{\hat{D}_{kl}(r_\xi)}\Grad{ r_{l,\eta}}+\delta\vt\Grad {r_{k,\eta}}+\ep \vt\erki{\xi}\Grad r_{k,\eta}\right)= \\
	-&\vt\pt\vr_\eta+\frac{\ep}{\vt^2+\eta^2}-\theta_\eta^5-\lr{\langle{\pi_{m}\rangle}_\eta(\theta_\eta,r_\eta)+\frac{\beta}{3}\theta_\eta^4}\Div\vu_\eta
	+G_\eta(t,x).
	\end{split}
	\end{equation}
By $\eta$ we denoted mollification of functions in the following sense.  For functions $\vr(t,x)$, $\vu(t,x)$, $r(t,x)$, $G(t,x)$ it is the mollification in the time variable, in case of energy, pressure and the transport coefficients it is the mollification with respect to all independent variables; the functions are assumed to be extended to the whole space in the following way
	$$\langle{a\rangle}(z)=\left\{
	\begin{array}{c}
	a(z)\ \mbox{if}\ z\in(0,\infty)^M\\
	\max\{\inf_{z\in(0,\infty)^M}a(z),0\}
	\end{array}
	\right. , \ M=1,2.$$
Moreover we set
	$$\theta_\eta=\frac{\sqrt{\vt^2+\eta^2}}{1+\eta\sqrt{\vt^2+\eta^2}}$$
and
	$$ K_{\xi,\eta}=K_{\eta}(\vr,\vt,r_\xi)=\int_1^{\vt}\left([\kappa_{\ep,\xi}]^\eta+\sum_{k,l=1}^n[\hat D_{kl}(r_\xi)C_0(\vr,s)]^\eta\right)~ds,$$
	$$[\kappa_{\ep,\xi}]^\eta=\ep\vr^n_\xi+\langle{\kappa\rangle}_\eta(\vr_\eta,\theta_\eta)+\beta\theta_\eta^B, \quad [\hat D_{kl}(r_\xi)C_0(\vr,s)]^\eta=\hat D_{kl}(r_{\xi})\langle{C_0\rangle}_\eta(\vr_\eta,\theta_\eta),$$
	$$G(t,x)=2\vr^n|\D(\vu)|^2+\lambda|\lap^{s}\Grad(\vr\vu)|^{2}+\lambda\ep|\lap^{s+1}\vr|^{2}+\ep\frac{1}{\vr}\frac{\partial\pi_c(\vr)}{\partial\vr}|\Grad\vr|^2.$$

Now, for $\xi,\eta$, $N$ being fixed, we  can combine the Ladyzhenskaya-Solonnikov-Uralceva theorem about the existence and uniqueness of solutions to quasilinear parabolic problem \cite{LSU} together with standard estimates for parabolic equations in order to find unique $\vt=\vt^\eta$ solving \eqref{eqK}, such that 
	$$\vt^\eta\in C([0,T];C^{2,l}(\Omega))\cap C^{1}([0,T]\times\Omega),\quad\pt\vt^\eta\in C^{0,l/2}([0,T]; C(\Omega)).$$
The following bounds are satisfied uniformly with respect to $\eta$
	\begin{equation}\label{un_vt_eta}
	\mbox{ess}\sup_{t\in(0,T)}\|\vt^\eta\|_{W^{1,2}}^2+\intT{\lr{\|\pt\vt^\eta\|_{L^2(\Omega)^2}+\|\lap\vt^\eta\|_{L^2(\Omega)}^2}}\leq c
	\end{equation}
with the constant $c$ which depends only on the following quantities: $\|\vu\|_{ C([0,T];X_N)}$, $\|\vr\|_{C^{1}([0,T]\times\Omega)}$, $\|\lap^{s+1}\vr\|_{C([0,T]\times\Omega)}$, $\|\vr^{-1}\|_{L^\infty((0,T)\times\Omega)}$, $\|r\|_{C([0,T],Y_N)^n}$, $\|r_\xi\|_{C^1([0,T],Y_N)^n}$ and the initial value $\|\vt_0\|_{W^{1,2}(\Omega)}$.

In addition, analogously as in \cite{FN}, Chapter 3, one can prove a comparison principle in the class of strong super and subsolutions to \eqref{eqK}, which thanks to a presence of singular terms  $\frac{\ep}{\vt^2+\eta^2}$, $\theta_\eta^5$ causes that the temperature $\vt^\eta$ stays away from $0$ and is bounded from above
	$$\left\|\vt^\eta\right\|_{L^{\infty}((0,T)\times\Omega)}+\left\|(\vt^\eta)^{-1} \right\|_{L^{\infty}((0,T)\times\Omega)}\leq c.$$
With these bounds at hand we may let $\eta\to 0^+$ in equation \eqref{eqK} to obtain the unique solution of \eqref{app_vt1} belonging to the class \eqref{regvt}. Moreover, the uniform estimates from \eqref{un_vt_eta} imply compactness of solutions in $L^2([0,T];W^{1,2}(\Omega))$. The continuity of the mapping $(\vu,r)\mapsto\vt(\vu,r)$ then follows from uniqueness of solutions established above. $\Box$

\subsection{Fixed point argument}\label{ex:mom}

Having all the necessary elements prepared, we are ready to apply the fixed point argument. We use the Schauder fixed point theorem to find a solution to the momentum and the species mass balance equations.

More precisely, we prove that there exists $\tau=\tau(N)$ such that $\vu$ solves  the approximate momentum equation \eqref{FG} and $r_{k}$ solves the species balance equation \eqref{rk_gal} with regularized coefficients,
	\begin{multline}\label{rk_gal1}
	\begin{split}
	&\intO{(\delta r_{k}+\erk)(t)\phi}-\intO{(\delta r_{k}+\erk)(0)\phi}
	=\\
	&\inttO{\left(\erk\vu-(\delta+\ep\erki{\xi})\Grad r_{k}\right)\cdot\Grad\phi}
	-\inttO{\left({C_{0}}(\vr,\vt)\sumlN\hat{D}_{kl}(r_\xi)\Grad r_{l}\right)\cdot\Grad\phi}\\
	-&\inttO{\left({C_{0}(\vr,\vt)}\Grad\log \vt\sumlN\hat D_{kl}(r_\xi)\right)\cdot\Grad\phi}
	+\inttO{\frac{\vr^{n}\vt\vw_k}{m_k}\phi},
	\end{split}
	\end{multline}
for $\phi \in C^\infty(\overline{\Omega})$ and $t \in (0,\tau]$. 
To this purpose we consider the following mapping
	\begin{equation}\label{T}
		\begin{array}{c}
		\vspace{0.5cm}
		{\cal{T}}:C([0,\tau];X_{N})\times C([0,\tau];Y_{N})\rightarrow C([0,\tau];X_{N})\times C([0,\tau];Y_{N}),\\
		{\cal{T}}(\vc{v}, z)=(\vu, r),
		\end{array}
	\end{equation}
which attains a solution to the following problem
	\begin{equation}\label{T1}
		\begin{array}{c}
		\vspace{0.5cm}
		\vu(t)={\cal{M}}_{\vr(\vc{v})(t)}\left[\vc{m}^0+\int_{0}^{t}{P_{X_N}{\cal {N}}(\vc{v}, z)(s) {\rm d}s}\right],\\
		r_{k}(t)={\cal K}_{z_{k}(t)}\left[(\delta z_{k}+\mbox{e}^{z_{k}})^0+\int_0^t{P_{Y_N}{\cal {L}}_k(\vc{v}, z)(s) {\rm d}s}\right],
		\end{array}
	\end{equation}
where 
	\begin{multline}\label{uN_vt}
	\langle{{\cal {N}} (\vc{v},z),\vcg{\phi}\rangle} =\intO{(\vr \vc{v} \otimes \vc{v}):\Grad\vcg{\phi}} -\intO{ 2\vr^{n}\D(\vc{v}):\Grad\vcg{\phi}}
	+\intO{ {\pi}(\vt,\vr, z_k)\Div\vcg{\phi}}\\
	+\lambda\intO{\vr\Grad\lap^{2s+1}\vr \cdot\vcg{\phi}}
	+\lambda\intO{\vr\lap^{s}\Div(\lap^{s}\Grad(\vc{v}\vr))\cdot\vcg{\phi}}-\ep\intO{(\Grad\vr\cdot\Grad)\vc{v}\cdot\vcg{\phi}},
	\end{multline}
and
	$${\cal{M}}_{\vr}\left[\cdot\right]:X_{N}\rightarrow X_{N},\quad\intO{\vr{\cal{M}}_{\vr}\left[\vc{w}\right]\vcg{\phi}}=\langle{\vc{w},\vcg{\phi}\rangle}, \quad\vc{w},\vcg{\phi}\in X_{N}.$$
Similarly,
	\begin{multline}
	\langle{{\cal {L}} _{k}(\vc{v},z),\phi\rangle}= \intO{\frac{\vr^{n}\vt\vw_k}{m_k}\phi}\\
	+ \intO{\left(\ezk\vc{v}-(\delta+\ep\ezki{\xi})\Grad z_{k}
	-{C_{0}}(\vr,\vt)\lr{\sumlN\hat{D}_{kl}(z_\xi)\Grad z_{l}+\Grad\log\vt\sumlN\hat{D}_{kl}(z_\xi)}\right)\cdot\Grad\phi},
	\end{multline}
and
	$${\cal{K}}_{z_{k}}\left[\cdot\right]:Y_{N}\rightarrow Y_{N},\quad\intO{(\delta+\mbox{e}^{z_{k}}){\cal{K}}_{z_{k}}\left[y\right]\phi}=\langle{y,\phi\rangle}, \quad y,\phi\in Y_{N}.$$
Next, we consider a ball ${\cal B}$ in the space $C([0,\tau];X_{N})\times C([0,\tau];Y_{N})$:
	$${\cal B}_{R,\tau}=\left\{(\vc{v},z)\in C([0,\tau];X_{N})\times C([0,\tau];Y_{N}):\|\vc{v}\|_{C([0,\tau];X_N)}+\sumkN\|z_{k}\|_{C([0,\tau];Y_N)}\leq R\right\}$$

We need to show that the operator ${\cal T}$ is continuous and maps ${\cal B}_{R,\tau}$ into itself, provided $\tau$ is sufficiently small. First observe that we have
	\begin{multline}\label{uN}
	\|{\cal {N}}(\vu, r)\|_{X_N}\leq c\left[\lr{\|\vr^{n}\|_{Y_N}+\|\vr\|_{L^\infty(\Omega)}}\left(\|\vu\|_{X_n}^2+\|\vu\|_{X_N}\right)
	+\|\vr\|_{L^\infty(\Omega)}^\gamma+ \|\vt\|^4_{L^\infty(\Omega)}\right.\\
	\left.+\|\vr^{n}\|_{Y_N}\|\vt\|_{L^\infty(\Omega)}+\|\vr\|_{L^\infty(\Omega)}\lr{\|\vr\|_{W^{4s+3,\infty}(\Omega)}+\|\vr\|_{W^{4s+2,\infty}(\Omega)}\|\vu\|_{X_{N}}}\right],
	\end{multline}
and	
	\begin{multline}\label{rN}
	\|{\cal {L}}_k(\vu, r)\|_{Y_N}\leq c[\lr{1+\|\erki{\xi}\|_{L^\infty(\Omega)}}   \lr{\|\vu\|_{X_N}+\|r_{k}\|_{Y_N}}\\
	+\|\vr\|_{L^\infty(\Omega)}\lr{\|\vt\|_{L^\infty(\Omega)}\|r\|_{Y_N}+\|\vt\|_{W^{1,2}(\Omega)}}+\|\vr^{n}\|_{Y_N}\|\vt\|_{L^\infty(\Omega)}].
	\end{multline}
To justify that the higher order gradients of the density  are bounded, one needs to recall that the unique solution $\vr$ to the approximate continuity equation \eqref{apcontN_2} is smooth in the space variable.  Therefore we can put the term $\Div(\vr\vu)$ to the r.h.s. of \eqref{apcontN_2} and then bootstrap the procedure leading to regularity \eqref{regvr}, see e.g. \cite{LSU}, Chapter IV. By this argument, the term $\vr\Grad\lap^{2s+1}\vr$ in the approximate momentum equation makes sense, i.e. it is bounded in $L^{\infty}(0,\tau; C^\infty(\Omega))$.

From estimates \eqref{uN}, $\eqref{rN}$, the estimates established in Lemmas \ref{kropa}, \ref{lemma_vt} and  from \eqref{T1} it follows that for sufficiently small $\tau$, operator ${\cal T}$ maps the ball ${\cal B}_{R,\tau}$ into itself. Moreover, ${\cal T}$ is a continuous mapping and its image consists of Lipschitz functions, thus it is compact in ${\cal B}_{R,\tau}$. It allows us to apply the theory of topological degree to infer that there exists at least one fixed point $(\vu, r)$ solving \eqref{FG} and \eqref{rk_gal1} on $[0,\tau]$.

\subsection{Uniform estimates and global in time solvability}

In order to extend this solution for the whole time interval $[0,T]$, we need a uniform bound of the solution. It follows from \eqref{T1} that $\vu$ and $r_k$ are continuously differentiable function, therefore, system \eqref{FG}, \eqref{rk_gal1} may be transformed to the following one
	\begin{equation}\label{FG_diff}
	\begin{split}
	&\intO{\pt(\vr\vu)\cdot\vcg{\phi}}
	-\lambda\inttO{\vr\lap^{2s+1}(\vr\vu) \cdot\vcg{\phi}}
	-\intO{(\vr{\vu}\otimes{\vu}):\Grad \vcg{\phi}}
	+\intO{\vS:\Grad \vcg{\phi}}\\
	-&\intO{\pi\Div \vcg{\phi}}
	-\lambda\inttO{\vr\Grad\lap^{2s+1}\vr\cdot \vcg{\phi}}+\ep\inttO{(\Grad\vr\cdot\Grad)\vu\cdot \vcg{\phi}}=0,
	\end{split}
	\end{equation}
satisfied for any $\vcg{\phi}\in X_N$ on $(0,\tau)$, $\vr\vu(0)=P_{X_N}\vc{m}^0$, 
 and
	\begin{multline}\label{rk_gal_diff}
	\begin{split}
	\intO{\pt(\delta r_{k}+\erk)\phi}
	=&\intO{\left(\erk\vu-(\delta+\ep\erki{\xi})\Grad r_{k}\right)\cdot\Grad\phi}
	+\intO{\left({C_{0}}(\vr,\vt)\sumlN\hat{D}_{kl}(r_\xi)\Grad r_{l}\right)\cdot\Grad\phi}\\
	&+\intO{\left({C_{0}(\vr,\vt)}\Grad\log \vt\sumlN\hat D_{kl}(r_\xi)\right)\cdot\Grad\phi}
	+\intO{\frac{\vr^{n}\vt\vw_k}{m_k}\phi}
	\end{split}
	\end{multline}
holds for any $\phi\in Y_N$ on $(0,\tau)$ and $(\delta r_{k}+\erk)(0)=P_{Y_N}[(\delta r^0_{k}+{\rm e}^{r_k^0})]$.
Therefore we can test \eqref{FG_diff} and \eqref{rk_gal_diff} by $\vu(t)$ and $r_{k}(t)$ respectively. For the approximate momentum equation, using the continuity equation, we obtain the kinetic energy balance
	\begin{multline}\label{firstep}
	\Dt\intO{\lr{{1\over2}\vr|\vu|^2+\frac{\lambda}{2}|\Grad^{2s+1}\vr|^{2}+\vr e_c(\vr)}}+\ep\intO{\frac{1}{\vr}\frac{\partial\pi_c}{\partial\vr}|\Grad\vr|^2}\\
	+\intOB{2\vr^n|\D(\vu)|^2
	+\lambda|\lap^{s}\Grad(\vr\vu)|^{2}+\lambda\ep|\lap^{s+1}\vr|^{2}}=
	\intO{\pi_\beta\Div\vu}.
	\end{multline}
Adding to this equality \eqref{app_vt1} integrated with respect to space and integrating the resulting sum with respect to time we obtain
	\begin{multline}\label{main_ener}
	\intO{\lr{{1\over2}\vr|\vu|^2(t)+\frac{\lambda}{2}|\Grad^{2s+1}\vr|^{2}(t)+\vr e(t)}}+\ep\inttO{\vt^{5}}\\
	=\inttO{\frac{\ep}{\vt^2}}+\intO{\lr{{1\over2}\vr|\vu|^2(0)+\frac{\lambda}{2}|\Grad^{2s+1}\vr|^{2}(0)+\vr e(0)}}.
	\end{multline}
Next taking $\phi=r_{k}(t)$ in \eqref{rk_gal_diff} 
and then summing with respect to $k=1,\ldots,n$ we get
	\begin{multline}\label{firstrk}
	\sumkN\Dt\intOB{\delta\frac{r_{k}^2}{2}+\erk r_{k}-\erk}
	+\sumkN\intOB{(\delta+\ep\erki{\xi})|\Grad r_{k}|^2+{C_{0}({\vr},\vt)}\sumlN\hat{D}_{kl}(r_\xi)\Grad r_{l}\Grad r_{k}}\\
	=-\sumkN \intO{\erk\Div\vu}+\sumkN\intO{\frac{\vr^{n}\vt\vw_k}{m_k}r_{k}}-\intO{{C_{0}(\vr,\vt)}\Grad\log\vt\sum_{l,k=1}^n\hat{D}_{kl}(r_\xi)\Grad r_{k}}.
	\end{multline}
Since $\hat D_{kl}(r)$ is a symmetric, positive semi-definite matrix, the last term on the l.h.s. is nonnegative.
On the other hand, for $\phi=1$ we get from \eqref{rk_gal_diff} 
	\begin{equation}\label{single:mas}
	\intO{\pt(\delta r_{k}+\erk)}
	=\intO{\frac{\vr^{n}\vt\vw_k}{m_k}},
	\end{equation}
multiplying these equations by $m_k$, summing them and integrating with respect to time yields
	\begin{equation}\label{cons:mas}
	\delta\sumkN m_k \intO{r_k(t)}+\sumkN\intO{{m_k\erk}(t)}= \delta\sumkN m_k \intO{r_k(0)}+\sumkN\intO{{m_k\erk}(0)},
	\end{equation}
for $t\in(0,\tau)$.

Note that these equalities do not provide any uniform bounds for the velocity yet. Indeed, proving boundedness of r.h.s. of \eqref{main_ener} and \eqref{firstrk} requires an estimate of the temperature. The next step will be hence devoted to derivation of such estimate from the so called entropy balance equation.

\subsubsection{Entropy estimate}
Our aim now is to derive a fundamental estimate for our system. It can be viewed as a total global entropy balance. Indeed, using \eqref{entropycv1} as a definition of the entropy, we will rearrange the internal energy equation, the continuity equation and the species mass balance equations in order to get a relevant approximation of \eqref{entropy}. However, due to the very low regularity of $\erk$, we can only hope to have an integral equality rather then a pointwise one.

From Lemma \ref{lemma_vt} it follows in particular that $\vt=\vt(\vu,r)$ is bounded from below by a constant. Therefore, dividing internal energy equation \eqref{app_vt1} by $\vt$ is possible and the equation
	 \begin{equation}\label{previous}
	 \begin{split}
	&{\frac{4\beta}{3}\pt\vt^3} +\pt\left({\vr\log\vt}\right)+\Div\lr{\vu\vr\log\vt}+\Div\lr{\frac{4\beta}{3}\vu\vt^3}+\ep\lap\vr(1-\log\vt)\\
	-&\Div\left(\kappa_{\ep}(\vr,\vt,r_\xi)\Grad\log\vt\right)-\frac{\kappa_{\ep}(\vr,\vt,r_\xi)|\Grad\vt|^2}{\vt^2}
	+\sumkN  \Div \left(\frac{\vf_{k,\xi}}{m_k}-\delta \Grad r_{k}-\ep\erki{\xi}\Grad r_{k}\right)\\
	+&\sumkN  \left(\frac{\vf_{k,\xi}}{m_k}-\delta \Grad r_{k}-\ep\erki{\xi}\Grad r_k\right)\cdot\Grad\log\vt= 
	\frac{\ep}{\vt^3}-\ep\vt^{4}-\sumkN\erk\Div\vu+\ep\frac{1}{\vr\vt}\frac{\partial\pi_c}{\partial\vr}|\Grad\vr|^2\\
	+&\frac{2\vr^n|\D\vu|^2+\lambda|\lap^{s}\Grad(\vr\vu)|^2+\lambda\ep(\lap^{s+1}\vr)^{2}}{\vt}
	\end{split}
	\end{equation}
is satisfied on $(0,\tau)$. In the next step we sum \eqref{firstrk} (without summing over $k$) with \eqref{single:mas} and obtain
	\begin{multline*}
	\Dt\intOB{\erk r_{k}+\delta r_{k}\lr{\frac{r_{k}}{2}+1}}
	=\intO{\left(\erk\vu-(\delta+\ep\erki{\xi})\Grad r_{k}-{C_{0}}(\vt,\vr)\sumlN\hat{D}_{kl}(r_\xi)\Grad r_{l}\right)\cdot\Grad r_{k}}\\
	-\intO{\left({C_{0}(\vt,\vr)}\sumlN\hat D_{kl}(r_\xi)\Grad\log \vt\right)\cdot\Grad r_{k}}
	+\intO{\frac{\vr^n\vt\vw_k}{m_k}(r_{k}+1)}.
	\end{multline*}	

Summing the above equations with respect to $k=1,\ldots,n$ and subtracting the obtained sum from \eqref{previous} integrated over $\Omega$ we get
	 \begin{equation}\label{main_entropy}
	 \begin{split}
	&\Dt\intOB{{\frac{4\beta}{3}\vt^3} +{\vr\log\vt}-\sumkN \erk r_{k}-\delta\sumkN r_{k}\lr{\frac{r_{k}}{2}+1}}\\
	+&\intO{\Div\left(\frac{4\beta}{3}\vu\vt^3+\vu\vr\log\vt+\vu\sumkN\erk\right)}\\
	-&\intOB{\Div\left(\kappa_{\ep}(\vr,\vt,r_\xi)\Grad\log\vt\right)
	-\sumkN  \Div \left(\frac{\vf_{k,\xi}}{m_k}-\delta \Grad r_{k}-\ep\erki{\xi}\Grad r_k\right)}\\
	+&\ep\intO{\vt^{4}}+\ep\intO{\Grad\vr\cdot\Grad\log\vt}-\sumkN\intO{(\delta+\ep\erki{\xi})\Grad r_{k}\cdot\Grad\log\vt}\\
	=&\intO{\vt^{-1}\lr{2\vr^n|\D\vu|^2
	+\lambda|\lap^{s}\Grad(\vr\vu)|^2+\lambda\ep(\lap^{s+1}\vr)^{2}}}\\
	+&\sum_{k,l=1}^n\intO{C_0(\vr,\vt)\hat{D}_{k,l}(r_\xi)\left(\Grad r_{k}+\Grad\log \vt\right)\left(\Grad r_{k}+\Grad\log \vt\right)}+\intO{\frac{\kappa_{\ep}(\vr,\vt, r_\xi)|\Grad\vt|^2}{\vt^2}}\\
	-&\sumkN\intO{\vr^n g_k\vw_k}
	+\sumkN\intO{(\delta+\ep\erki{\xi})|\Grad r_{k}|^2}+\ep\intO{\frac{1}{\vr\vt}\frac{\partial\pi_c}{\partial\vr}|\Grad\vr|^2}+\intO{\frac{\ep}{\vt^3}},
	\end{split}
	\end{equation}
since $\sumkN\vw_k=0$.
In order to get uniform estimates we need to control the l.h.s.; the simplest way to do it is to integrate the above inequality with respect to time and then to subtract it from the energy balance \eqref{main_ener}. Hence we get
	\begin{equation}\label{eq:entr}
	\begin{split}
	&\inttauO{\vt^{-1}\lr{2\vr^n|\D\vu|^2
	+\lambda|\lap^{s}\Grad(\vr\vu)|^2+\lambda\ep(\lap^{s+1}\vr)^{2}}}\\
	&+\sum_{k,l=1}^n\inttauO{{C_{0}(\vr,\vt)}\hat{D}_{kl}(r_\xi)\left(\Grad r_{l}+\Grad\log \vt\right)\left(\Grad r_{k}+\Grad\log \vt\right)}\\
	&+\inttauO{\frac{\kappa_{\ep}(\vr,\vt,r_\xi)|\Grad\vt|^2}{\vt^2}}
	-\sumkN\inttauO{\vr^n g_k\vw_k}\\
	&+\sumkN\inttauO{(\delta+\ep\erki{\xi})|\Grad r_{k}|^2} +\ep\inttauO{\frac{1}{\vr\vt}\frac{\partial\pi_c}{\partial\vr}|\Grad\vr|^2}
	+\inttauO{\frac{\ep}{\vt^3}}\\
	&+\ep\inttauO{\vt^{5}}+\intOB{\frac{4\beta}{3}\vt^3(0)+{\vr\log\vt}(0)-\sumkN \erk r_{k}(0)-\delta\sumkN r_{k}\lr{\frac{r_{k}}{2}+1}(0)}\\
	&+\intOB{{1\over2}\vr|\vu|^2(\tau)+\frac{\lambda}{2}|\Grad^{2s+1}\vr|^{2}(\tau)+\vr e(\tau)}\\
	=&\intOB{{1\over2}\vr|\vu|^2(0)+\frac{\lambda}{2}|\Grad^{2s+1}\vr|^{2}(0)+\vr e(0)}\\
	&+\intOB{\frac{4\beta}{3}\vt^3(\tau)+{\vr\log\vt}(\tau)-\sumkN \erk r_{k}(\tau)-\delta\sumkN \lr{\frac{r_{k}}{2}+1}(\tau)}+\inttauO{\frac{\ep}{\vt^2}}\\
	&+\ep\inttauO{\vt^{4}}+\ep\inttauO{\Grad\vr\cdot\Grad\log\vt}-\sumkN\inttauO{(\delta+\ep\erki{\xi})\Grad r_{k}\cdot\Grad\log\vt}.
	\end{split}
	\end{equation}
To control the r.h.s. we take adventage of the fact that the heat conductivity coefficient depends on the partial densities. We  write
	\begin{multline*}
	\ep\intO{\Grad\vr\cdot\Grad\log\vt}-\sumkN\intO{(\delta+\ep\erki\xi)\Grad r_{k}\cdot\Grad\log\vt}
	\leq \ep\|\sqrt{\lr{\vr\vt}^{-1}}\Grad\vr\|_{L^2(\Omega)}\|\sqrt{\vr\vt}\Grad\log\vt\|_{L^2(\Omega)}\\
	 +\delta \sumkN\|\Grad r_{k}\|_{L^2(\Omega)}\|\Grad\log\vt\|_{L^2(\Omega)}
	+\ep\sumkN\|\sqrt{\erki\xi}\Grad r_{k}\|_{L^2(\Omega)}\|\sqrt{\erki\xi}\Grad\log\vt\|_{L^2(\Omega)}.
	\end{multline*}
Obviously,
\begin{equation*}
\ep\sumkN\|\sqrt{\erki\xi}\Grad r_{k}\|_{L^2(\Omega)}\|\sqrt{\erki\xi}\Grad\log\vt\|_{L^2(\Omega)}
\leq \sumkN\frac{\ep}{2}\intO{\erki\xi|\Grad r_k|^2}+\sumkN\frac{\ep}{2}\intO{\erki\xi|\Grad \log\vt|^2},
\end{equation*}
and the last term  is absorbed by the $\ep$-dependent  part of the heat flux. Indeed, due to \eqref{notation_x} and \eqref{21a}
$$\sumkN\frac{\ep}{2}\intO{\erki\xi|\Grad \log\vt|^2}\leq\sumkN\frac{\ep}{2\Un{m}}\intO{m_k\erki\xi|\Grad \log\vt|^2}
\leq\frac{\ep}{2\Un{m}}\intO{\vr_\xi^n|\Grad\log\vt|^2}.$$

To control the positive part of the entropy at time $\tau$ and the negative part of it at the initial time $t=0$ we note that 
	\begin{multline*}
	\intOB{\frac{4\beta}{3}\vt^3(\tau)+{\vr\log\vt}(\tau)-\sumkN \erk r_{k}(\tau)-\delta\sumkN r_{k}\lr{\frac{r_{k}}{2}+1}(\tau)}\\
	\leq\intOB{\frac{4\beta}{3}\max\left\{\frac{4}{3},\vt(\tau)\right\}^3+[\vr\log\vt(\tau)]_++\sumlN[-\erk r_k(\tau)]_++\delta\sumkN\left[-r_{k}\lr{\frac{r_{k}}{2}+1}(\tau)\right]_+},
	\end{multline*}
where by $[\cdot]_{+}$ we denoted the positive part of a function. Hence, on account of \eqref{main_ener} we may write
	\begin{multline}\label{un:s}
	\intOB{\frac{4\beta}{3}\vt^3(\tau)+{\vr\log\vt}(\tau)-\sumkN \erk r_{k}(\tau)-\delta\sumkN r_{k}\lr{\frac{r_{k}}{2}+1}(\tau)}\\
	\leq c(\Omega,\tau)+\intOB{\beta\vt^4(\tau)+{\vr\vt}(\tau)}\leq c(\Omega,\tau)+\inttauO{\frac{\ep}{\vt^2}}.
	\end{multline}
On the other hand, we easily verify that
	$$\inttauO{\frac{\ep}{\vt^2}}+\ep\inttauO{\vt^{4}}\leq c+ \inttauO{\frac{\ep}{\vt^3}}+\ep\inttauO{\vt^{5}},$$
which appears on the l.h.s. of \eqref{eq:entr}.


Summarizing, we have shown the following estimate
	\begin{equation}\label{est:entr}
	\begin{split}
	\sup_{t\in[0,\tau]}&\intOB{{1\over2}\vr|\vu|^2(t)+\frac{\lambda}{2}|\Grad^{2s+1}\vr|^{2}(t)+\vr e(t)}\\
	+&\inttauO{\vt^{-1}\lr{2\vr^n|\D\vu|^2
	+\lambda|\lap^{s}\Grad(\vr\vu)|^2+\lambda\ep(\lap^{s+1}\vr)^{2}}}\\
	+&\sum_{k,l=1}^n\inttauO{{C_{0}(\vr,\vt)}\hat{D}_{kl}(r_\xi)\left(\Grad r_{l}+\Grad\log \vt\right)\left(\Grad r_{k}+\Grad\log \vt\right)}\\
	+&\intTO{\frac{\kappa_{\ep}(\vr,\vt,r_\xi)|\Grad\vt|^2}{\vt^2}}
	-\sumkN\inttauO{\vr^n g_k\vw_k}
	+\sumkN\inttauO{(\delta+\ep\erki{\xi})|\Grad r_{k}|^2}\\
	+&\ep\inttauO{\frac{1}{\vr\vt}\frac{\partial\pi_c}{\partial\vr}|\Grad\vr|^2}
	+\inttauO{\frac{\ep}{\vt^3}}
	+\ep\inttauO{\vt^{5}}\leq c.
	\end{split}
	\end{equation}

Taking $s$ from the density-regularizing term sufficiently large one can show that the density is separated from $0$ uniformly with respect to all approximation parameters except for $\lambda$. {This property was observed  by Bresch and Desjardins in \cite{BD, BDC} where the case of single-component heat-conducting fluid was discussed. Recalling their reasoning we may use} the Sobolev embedding $\|\vr^{-1}\|_{L^{\infty}(\Omega)}\leq c\|\vr^{-1}\|_{W^{3,1}(\Omega)}$ and
	\begin{equation*}
	\|\Grad^{3}\vr^{-1}\|_{L^{1}(\Omega)}\leq (1+\|\Grad^3\vr\|_{L^{2}(\Omega)})^{3}(1+\|\vr^{-1}\|_{L^{4}(\Omega)})^{4},
	\end{equation*}
where the last term is bounded on account of \eqref{est:entr} and the assumption that $\gamma^{-}\geq4$. So,  provided that $2s+1\geq 3$  we have

	\begin{equation}\label{below}
	\|\vr^{-1}\|_{L^{\infty}((0,\tau)\times\Omega)}\leq c(\lambda)\quad \mbox{a.e.\ in} \ (0,\tau)\times\Omega.
	\end{equation}

\subsubsection{Global in time existence of solutions.}
The uniform estimates for $\vu$ and $r_{k}$ can be summarized as follows

	\begin{equation*}\label{indnt}
	\|\sqrt{\vr}\vu\|_{L^{\infty}(0,\tau;L^2(\Omega))}+\sqrt{\lambda}\|\lap^s\Grad(\vr\vu)\|_{L^{2}(0,\tau;L^2(\Omega))}\leq c
	\end{equation*}
and
	$$\sqrt{\delta}\|r_{k}\|_{L^\infty(0,\tau;L^2(\Omega))}+\sqrt{\delta}\|r_{k}\|_{L^2(0,\tau;W^{1,2}(\Omega))}\leq c.$$

Moreover, the density $\vr$ is bounded from below by a positive constant on account of \eqref{below}. By the equivalence of norms on the finite dimensional spaces $X_N$ and $Y_N$ we can thus deduce the uniform bounds for $\vu$ and $r_{k}$ in $C([0,\tau]; X_N)$ and $C([0,\tau]; Y_N)$, $k=1,\ldots,n$, respectively.  As the global estimate does not blow up, we can return to the procedure of construction of local in time solution described in Sections \ref{ex:cont}--\ref{ex:mom}.  By a contradiction argument we get  a solution defined on $[0,T]$ for arbitrary but finite $T>0$, exactly as in \cite{MPZ1}.

	\begin{rmk}
	Finally, note that \eqref{est:entr} is global in time and independent of $\xi$. Therefore it is straightforward to let $\xi\to 0$, since $r_{k,\xi}\to r_{k}$ strongly in $C([0,T],Y_N)$. 
	\end{rmk}
This remark completes the proof of Theorem \ref{th1}. $\Box$

\section{Limit passage in the Galerkin approximation (proof of Theorem \ref{t 2})}\label{Sec:Gal}

The purpose of this section is to let $N\to \infty$ in the equations of approximate system introduced in Section \ref{Sec:app}. We start with summarizing all the estimates that are uniform with respect to $N$ derived mostly from \eqref{est:entr} and its consequences. This will be done in Subsection \ref{UB}, then in Subsection \ref{Npass} we use these estimates to extract the weekly convergent subsequences and to prove that the limit passage $N\to\infty$ can be performed.

\subsection{Estimates independent of $N$}\label{UB}
Note that the above estimates are not only uniform with respect to time but also with respect to $N$. From \eqref{un:s} and \eqref{est:entr} we get  $\vr_N s(\vt_N,\vr_N, r_N)\in L^\infty(0,T;L^1(\Omega))$, more specifically we have
	\begin{equation}\label{ess_entr}
	\|\vr_N\log\vt_N\|_{L^\infty(0,T;L^1(\Omega))}+\|\erkN r_{k,N}\|_{L^\infty(0,T;L^1(\Omega))}+\delta\|r_{k,N}^2\|_{L^\infty(0,T;L^1(\Omega))}\leq c,
	\end{equation}
also, from \eqref{main_ener}, we get that
	\begin{multline}\label{ess_ener}
	\left\|\vr_N|\vu_N|^2\right\|_{L^\infty(0,T;L^1(\Omega))}+\sqrt{\lambda}\left\|\Grad^{2s+1}\vr_N\right\|_{L^\infty(0,T;L^2(\Omega))}\\
	+\left\|\vr_N e_c(\vr_N)\right\|_{L^\infty(0,T;L^1(\Omega))}+\|\beta\vt_N^4\|_{L^\infty(0,T;L^1(\Omega))}+\|\vr_N\vt_N\|_{L^\infty(0,T;L^1(\Omega))}\leq c.
	\end{multline}
In addition, we have the estimates following from boundedness of the entropy production rate:\\
\begin{itemize}
\item[-]
the velocity estimates
	\begin{equation}\label{entr_vel}
	\left\|\sqrt{\frac{2\vr_N^n}{\vt_N}}\D\vu_N\right\|_{L^2((0,T)\times\Omega)}
	+\left\|\sqrt{\frac{\lambda}{\vt_N}}\lap^{s}\Grad(\vr_N\vu_N)\right\|_{L^2((0,T)\times\Omega)}\leq c,
	\end{equation}
\item[-]
the density estimates
	\begin{equation}\label{entr_dens}
	\sqrt{\lambda\ep}\left\|\frac{\lap^{s+1}\vr_N}{\sqrt{\vt_N}}\right\|_{L^2((0,T)\times\Omega)}+\sqrt{\ep}\left\|{\frac{1}{\sqrt{\vr_N\vt_N}}\sqrt{\frac{\partial\pi_c}{\partial\vr_N}}\Grad\vr}\right\|_{L^2((0,T)\times\Omega)}\leq c,
	\end{equation}
\item[-]
the temperature estimates
	\begin{equation}\label{entr_temp}
	\left\|\frac{\sqrt{\kappa_{\ep}(\vr_N,\vr_N^n,\vt_N)}\Grad\vt_N}{\vt_N}\right\|_{L^2((0,T)\times\Omega)}+\left\|\frac{\ep}{\vt_N^3}\right\|_{L^1((0,T)\times\Omega)}
	+\left\|\ep\vt_N^{5}\right\|_{L^1((0,T)\times\Omega)}\leq c,
	\end{equation}
\item[-]
the species densities estimates
	\begin{equation}\label{entr_part}
	\left\|\sqrt{\delta+\ep\erkN}\Grad r_{k,N}\right\|_{L^2((0,T)\times\Omega)}\leq c;
	\end{equation}
\end{itemize}
moreover, we can write
	\begin{equation}\label{estF}
	\begin{split}
	&\sum_{k,l=1}^n\intTO{{C_{0}(\vr_N,\vt_N)}\hat{D}_{kl}\left(\Grad r_{l,N}+\Grad\log \vt_N\right)\left(\Grad r_{k,N}+\Grad\log \vt_N\right)}\\
	=&\sumkN\intTO{ \frac{\vf_k(\vr_N,\vt_N, r_{N})}{m_k}\cdot\left(\Grad r_{k,N}+\Grad\log \vt\right)}\\
	=&\intTO{\frac{\pi_{m}(\vt_N,r_N)}{C_0(\vr_N,\vt_N)\vt_N}\sumkN \frac{\vf_k^2(\vr_N,\vt_N, r_{N})}{m_k \erkN}}\leq c.
	\end{split}
	\end{equation}

\bigskip

\noindent{\bf Temperature estimates.}
One of the main consequences of \eqref{est:entr} is \eqref{entr_temp} which, for $\kappa_{\ep}(\vr^n,\vr,\vt)$ satisfying \eqref{notation_x}, provides a priori estimates for the temperature
	\begin{equation}\label{main-temp-est}
	\|(1+\sqrt{\ep\vr_N^n}+\sqrt{\vr_N})\Grad\log\vt_N\|_{L^2((0,T)\times \Omega)}+  \|\sqrt{\vr_N}\Grad{\vt_N}\|_{L^2((0,T)\times \Omega)} \|\sqrt{\beta}\Grad\vt_N^{a}\|_{L^2((0,T)\times \Omega)} \leq c, 
	\end{equation}
where $a\in [0,\frac{B}{2}]$ and $B\geq8$. To control the full norm of $\vt_N^a$ in $L^2(0,T;W^{1,2}(\Omega))$ we combine the above estimates with \eqref{ess_ener}. Therefore, the Sobolev imbedding gives 
	\begin{equation}\label{maxtemp}
	\|\sqrt{\beta}\vt_N\|_{ L^{B}(0,T; L^{3B}(\Omega))} \leq c.
	\end{equation}

\bigskip

\noindent {\bf Estimates of the species densities.}
From\eqref{ess_entr} and  \eqref{entr_part} it follows that
	\begin{equation}\label{grad_est}
	|\Grad\erkN|\leq2|\Grad\sqrt{\erkN}|\sqrt{\erkN}
	\end{equation}
is bounded in ${L^{2}(0,T;L^1(\Omega))}$, thus, by the Sobolev imbedding, $\erkN$ is bounded in $L^2(0,T;L^{\frac{3}{2}}(\Omega))$. Returning to \eqref{grad_est} we get
	\begin{equation} \label{29}
	\|\Grad \erkN\|_{L^{\frac{4}{3}}(0,T; L^{\frac 65}(\Omega))}\leq c;
	\end{equation}
using once more the Sobolev imbedding theorem and the bound in $L^\infty(0,T;L^1(\Omega))$ we end up with  
	\begin{equation} \label{29a}
	\|\erkN\|_{L^{\frac{5}{3}}((0,T)\times \Omega)}\leq c(\ep).
	\end{equation}
Having this, we  return to \eqref{grad_est} to deduce
	$\|\Grad\erkN\|_{L^{\frac{5}{4}}((0,T)\times\Omega)}\leq c(\ep)$.

\noindent {\bf Kinetic energy estimate.} We  now integrate \eqref{firstep} with respect to time to get
	\begin{equation}\label{kinetic}
	\begin{split}
	&\intO{\lr{{1\over2}\vr_N|\vu_N|^2+\frac{\lambda}{2}|\Grad^{2s+1}\vr_N|^{2}+\vr_N e_c(\vr_N)}(T)}+\ep\intTO{\frac{1}{\vr}\frac{\partial\pi_c(\vr_N)}{\partial\vr_N}|\Grad\vr_N|^2}\\
	+&\intTOB{2\vr_N^n|\D(\vu_N)|^2
	+\lambda|\lap^{s}\Grad(\vr_N\vu_N)|^{2}+\lambda\ep|\lap^{s+1}\vr_N|^{2}}\\
	=&\intTO{\lr{\pi_m+\frac{\beta}{3}\vt_N^4}\Div\vu_N}+\intO{\lr{{1\over2}\vr_N|\vu_N|^2+\frac{\lambda}{2}|\Grad^{2s+1}\vr_N|^{2}+\vr_N e_c(\vr_N)}(0)}.
	\end{split}
	\end{equation}
From \eqref{kinetic} it follows that 
	\begin{multline}\label{ener_vel}
	\left\|\sqrt{{\vr_N^n}}\D\vu_N\right\|_{L^2((0,T)\times\Omega)}
	+\left\|\sqrt{\lambda}\lap^{s}\Grad(\vr_N\vu_N)\right\|_{L^2((0,T)\times\Omega)}+\sqrt{\lambda\ep}\left\|{\lap^{s+1}\vr_N}\right\|_{L^2((0,T)\times\Omega)}\\
	+\sqrt{\ep}\left\|{\frac{1}{\sqrt{\vr_N}}\sqrt{\frac{\partial\pi_c(\vr_N)}{\partial\vr_N}}\Grad\vr_N}\right\|_{L^2((0,T)\times\Omega)} \\
	\leq c+\frac{\beta}{3}\intTO{\vt_N^4|\Div\vu_N|}+c\sumkN\intTO{\sqrt{\erkN}\vt|\sqrt{\erkN}\Div\vu_N|}
	\end{multline}
and the r.h.s. is bounded. Indeed, we have
	$$\frac{\beta}{3}\|\vt_N^4\Div\vu_N\|_{L^1((0,T)\times\Omega)}\leq c\beta\|\vt_N\|^4_{L^\infty(0,T;L^4(\Omega))}\|\Grad\vu_N\|_{L^2(0,T;L^\infty(\Omega))},$$
	$$\|\Grad\vu_N\|_{L^2(0,T;L^\infty(\Omega))}\leq\|\vu_N\|_{L^2(0,T;W^{3,2}(\Omega))},$$
and the r.h.s. is bounded provided  $2s+1\geq3$. To see it, we write
	\begin{equation}\label{onlyu}
	\Grad^3\vu_N=\Grad^3(\vr_N^{-1}\vr_N\vu_N)\approx\lr{\frac{\Grad^3\vr_N}{\vr_N^2}+\frac{(\Grad\vr_N)^3}{\vr_N^4}}\vr_N\vu_N+\vr_N^{-1}\Grad^{3}(\vr_N\vu_N)
	\end{equation}
and boundedness of the r.h.s. follows from \eqref{below}, \eqref{ess_ener} and the Cauchy inequality.
To estimate the last term from the r.h.s. of  \eqref{ener_vel}, observe that we have
	$$\sqrt{\erkN}\vt_N|\sqrt{\erkN}\Div\vu_N|\leq c(\epsilon)\erkN\vt_N^2+\epsilon\erkN|\D\vu_N|^2$$
and the last term is absorbed by the l.h.s. of \eqref{kinetic}, whence for the first term we use the following estimate
	$$\|\erkN\vt_N^2\|_{L^1((0,T)\times\Omega)}\leq \|\erkN\|_{L^{\frac{5}{3}}((0,T)\times\Omega)}\|\vt_N\|^2_{L^5((0,T)\times\Omega)};$$
both terms are bounded on account of  \eqref{maxtemp} and \eqref{29a} by a constant dependent on $\ep$.

\subsection{Passage to the limit with $N$}\label{Npass}
This subsection is devoted to the limit passage $N\to\infty$. Using estimates from the previous subsection we can extract  weakly convergent subsequences, whose limits satisfy the approximate system. It should be, however, emphasized that at this level we replace the weak formulation of the thermal energy by the weak formulation of the total energy.

\subsubsection{Strong convergence of the density and passage to the limit in the continuity equation}
From \eqref{kinetic} and a procedure similar to \eqref{onlyu} we deduce that
	\begin{equation}\label{conv_vu}
	\vu_N\to\vu\quad\mbox{weakly\ in}\  L^2(0,T;W^{2s+1,2}(\Omega)),
	\end{equation}
and
	\begin{equation}\label{conv_vr}
	\vr_N\to\vr \quad\mbox{weakly\ in} \ L^2(0,T;W^{2s+2,2}(\Omega)),
	\end{equation}
at least for a suitable subsequence. In addition the r.h.s. of the linear parabolic problem
	\begin{equation*}
		\begin{array}{rcl}
		\vspace{0.1cm}
		\pt\vr_N-\ep\lap\vr_N&=&\Div(\vr_N\vu_N),\\
		\vspace{0.1cm}
		\vr_N(0,x)&=&\vr^0_{\ep,\lambda},
		\end{array}
	\end{equation*}
is uniformly bounded in $L^2(0,T;W^{2s,2}(\Omega))$ and the initial condition is sufficiently smooth, thus, applying the $L^p-L^q$ theory to this problem we conclude that $\{\pt\vr_N\}_{n=1}^\infty$ is uniformly bounded in $L^2(0,T;W^{2s,2}(\Omega))$. Hence, the standard compact embedding implies
$\vr_N\to\vr$ {a.e. in} $(0,T)\times\Omega$
and therefore passage to the limit in the approximate continuity equation is straightforward.

\subsubsection{Strong convergence of the species densities}
To show this property we take advantage of the species mass balance equation, we observe that
	\begin{multline}
	\left|{\vf_k}(\vr_N,\vt_N, r_{N})\right|=\left| \sqrt{\frac{\sumkN\erkN}{C_0(\vr_N,\vt_N)}}\frac{\vf_k(\vr_N,\vt_N, r_{N})}{\sqrt{\erkN}}\right|\left| \sqrt{\frac{C_0(\vr_N,\vt_N)\erkN}{\sumkN\erkN}}\right|\\
	\leq c\sqrt{\vr_N\vt_N}\sqrt{\frac{\pi_{m}}{C_0(\vr_N,\vt_N)}}\sumkN \frac{|\vf_k(\vr_N,\vt_N, r_{N})|}{\sqrt{m_k \vt_N\erkN}},
	\end{multline}
so, since $\|\sqrt{\vr_N\vt_N}\|_{L^{4}((0,T)\times\Omega)}\leq\|\sqrt{\vr_N}\|_{L^\infty(0,T;L^6(\Omega))}\|\sqrt{\vt_N}\|_{L^{4}(0,T;L^{12}(\Omega))}$, we have that 
\begin{equation}\label{bF}
\|\vf_k(\vr_N,\vt_N, r_{N})\|_{ L^{\frac{4}{3}}((0,T)\times\Omega)}\leq c.
\end{equation}


Having this, we can repeat our reasoning from \cite{MPZ} (since $\frac{4}{3}>\frac{5}{4}$) to prove the following lemma providing compactness with respect to time.

	\begin{lemma}
	There exists a constant $c$ depending on the initial data, $T$, and the parameter $\ep$ such that
		\begin{equation}\label{time_part}
		\delta\|\pt r_{k,N}\|_{L^{\frac{5}{4}}(0,T;W^{-1,\frac{5}{4}}(\Omega))}\leq c.
		\end{equation}
	\end{lemma}

\pf We take any $\phi\in W^{1,5}(\Omega)\subset W^{1,2}(\Omega)$ such that $\|\phi\|_{W^{1,5}(\Omega)}\leq 1$ and decompose it into $\phi=\phi_1+\phi_2$, where $\phi_1$ is an orthogonal projection of $\phi$ (with respect to the scalar product induced by the norm of the space $L^2(\Omega)$)  onto $Y_N$. Using $\phi_1$ as a test function in \eqref{rk_gal_diff}  we show that
	\begin{equation} \label{31a}
	\begin{split}
	&\intO{\pt(\delta r_{k,N}+\erkN)\phi_1}\\
	&\quad=\intO{\left(\erkN\vu-(\delta+\ep\erkN)\Grad r_{k,N}+\frac{\vf_{k,N}}{m_k}\right)\cdot\Grad\phi_1}+\intO{\frac{\vr^n_N\vt_N\vw_k}{m_k}\phi_1}\\
	&\quad\leq c\sumkN\left(\|\vu\|_{L^\infty(\Omega)}\|\erkN\|_{L^{\frac{5}{3}}(\Omega)}+\delta\|\Grad r_{k,N}\|_{L^{2}(\Omega)}+\ep\|\Grad{\erkN}\|_{L^{\frac{5}{4}}(\Omega)}\right)\|\phi_1\|_{W^{1,5}(\Omega)}~~~~~\\
	&\qquad+c\sumkN\left(\left\|\vf_{k,N}\right\|_{L^{\frac{4}{3}}(\Omega)}+\|\erkN\|_{L^{\frac{5}{3}}(\Omega)}\|\vt_N\|_{L^{5}(\Omega)}\right)\|\phi_1\|_{W^{1,5}(\Omega)}.
	\end{split}
	\end{equation}
Then we have
	\begin{multline*}
	\begin{split}
	\|\pt r_{k,N}(t,\cdot)\|_{W^{-1,\frac 54}(\Omega)} &= \sup_{\phi \in W^{1,5}(\Omega) ; \|\phi\|\leq 1} \Big|\int_\Omega \pt r_{k,N}(t,\cdot) \phi \, \dx\Big| \\
	&=   \sup_{\phi \in W^{1,5}(\Omega); \|\phi\|\leq 1} \Big|\int_\Omega \pt r_{k,N}(t,\cdot) \phi_1 \dx\Big| = \int_\Omega  \big|\pt r_{k,N}(t,\cdot) \varphi_1 \big|\dx
	\end{split}
	\end{multline*}
for some $\varphi_1 \in W^{1,5}_0(\Omega) \cap Y_N$.
Hence
	\begin{equation} \label{31b}
	\|\pt r_{k,N}(t,\cdot)\|_{W^{-1,\frac 54}(\Omega)} \leq \sup_{\phi \in W^{1,5}(\Omega) \cap Y_N ; \|\phi\|\leq 1} \frac 1\delta \Big|\int_\Omega (\delta + {\rm e}^{r_{k,N}(t,\cdot)})\pt r_{k,N}(t,\cdot) \phi \, \dx\Big| 
	\end{equation}
and due to estimate \eqref{31a} we end up with $\|\pt r_{k,N}\|_{L^{\frac 54}(0,T;W^{-1,\frac 54}(\Omega))} \leq \frac{c(\ep)}{\delta}.$ $\Box$\\

We now apply the Aubin-Lions lemma to the sequence $r_{k,N}$, we deduce from \eqref{ess_entr}, \eqref{entr_part} and \eqref{time_part} that it is possible to extract a subsequence such that
	\begin{equation}\label{est_rkN}
	\begin{split}
	&r_{k,N}\to r_{k}\quad\mbox{weakly}^*\ \mbox{in\ }L^{\infty}(0,T;L^2(\Omega)),\\
	&\Grad r_{k,N}\to\Grad r_{k}\quad\mbox{weakly\ in\ }L^2((0,T)\times\Omega),\\
	&\pt r_{k,N}\to \pt r_{k}\quad\mbox{weakly\ in\ }L^{\frac{5}{4}}(0,T;W^{-1,\frac{5}{4}}(\Omega)),\\
	&r_{k,N}\to r_{k}\quad\mbox{strongly\ in\ }L^2(0,T;L^{p}(\Omega)),~ p<6;
	\end{split}
	\end{equation}
in particular, there exists a subsequence $r_{k,N}$ which converges to $r_k$ a.e. on $(0,T)\times\Omega$. Therefore also
	\begin{equation*}
	\erkN\to\erk\quad \mbox{a.e.\ on}\ (0,T)\times\Omega.
	\end{equation*}
Moreover, we have
	\begin{equation}\label{est_erkN}
		\begin{split}
		&\Grad\erkN\to\Grad \erk\quad\mbox{weakly\ in\ }L^2(0,T;L^1(\Omega))\cap L^{\frac{5}{4}}((0,T)\times\Omega),\\
		&\erkN\to \mbox{e}^{r_{k}}\quad\mbox{strongly\ in\ }L^{q}((0,T)\times\Omega),\, q<\frac 53.
		\end{split}
		\end{equation}
The above considerations imply that the species mass balance equations will be satisfied for $N\to\infty$ if we validate that the temperature sequence converges strongly. This is the purpose of the next subsection.

\subsubsection{Strong convergence of the temperature}
For the temperature we have
	\begin{equation}\label{conv_teta}
	\vt_N\to\vt\quad\mbox{weakly\ in\ }L^2(0,T;W^{1,2}(\Omega));
	\end{equation}
note that at this level, the time-compactness can be proved directly from the internal energy equation \eqref{app_vt}. Indeed, due to the continuity equation, we have
		 \begin{equation}\label{solo_vt}
		 \begin{split}
	\ptb{\vr_N\vt_N+\beta\vt_N^4} =&-\Div(\vu_N\vr_N\vt_N+\beta\vu_N\vt_N^4)+\Div\left(\kappa_{\ep}(\vr_N,\vr_N^n,\vt_N)\Grad\vt_N\right)\\
	&-\sumkN \Div \left(\vt_N\frac{\vf_{k,N}}{m_k}-\delta \vt_N\Grad r_{k,N}-\ep\vt_N\Grad\erkN\right)\\
	&+\frac{\ep}{\vt_N^2}-\ep\vt_N^{5}-\lr{\pi_{m,N}+\frac{\beta}{3}\vt_N^4}\Div\vu_N+\ep{\frac{1}{\vr_N}\frac{\partial\pi_c(\vr_N)}{\partial\vr_N}|\Grad\vr_N|^2}\\
	&+2\vr_N^n|\D(\vu_N)|^2+\lambda|\lap^{s}\Grad(\vr_N\vu_N)|^{2}+\lambda\ep|\lap^{s+1}\vr_N|^{2}=\sum_{i=1}^{10} I_i.
	\end{split}
	\end{equation}
On account of \eqref{entr_temp} and \eqref{ener_vel} the last 7 terms are bounded in $L^1((0,T)\times\Omega)$. Then it follows from \eqref{ess_ener}, \eqref{maxtemp} and \eqref{onlyu} that $I_1$ can be estimated as
	$$\|\vu_N\vr_N\vt_N\|_{L^{{12}/{11}}((0,T)\times\Omega)}\leq c\|\sqrt{\vr_N}\vu_N\|_{L^\infty(0,T;L^2(\Omega))}\|\sqrt{\vr_N}\|_{L^\infty(0,T;L^6(\Omega))}\|\vt_N\|_{L^\infty(0,T;L^4(\Omega))}\leq c,$$
	$$\|\vu_N\vt_N^4\|_{L^{8/7}(0,T;L^{8/3}(\Omega))}\leq c\|\vu\|_{L^2(0,T;L^\infty(\Omega))}\|\vt_N^4\|_{L^{8/3}((0,T)\times\Omega)}\leq c,$$
where we used the interpolation 
	\begin{equation}\label{eqT}
	\|\vt_N\|_{L^{32/3}((0,T)\times\Omega)}\leq c\|\vt_N\|^{1/4}_{L^{\infty}(0,T;L^4(\Omega))}\|\vt_N\|^{3/4}_{L^{8}(0,T;L^{24}(\Omega))};
	\end{equation}
hence the last term is bounded provided  $B\geq 8$.

For $I_2$ recall that  we have 
	$\kappa_{\ep}(\vr_N,\vr_N^n,\vt_N)\Grad\vt_N=\lr{\kappa_0+\ep\vr_N^n+\vr_N+\vr_N\vt_N^2+\beta\vt_N^B}\Grad\vt_N$, 
therefore using  estimate \eqref{entr_temp} and \eqref{ess_ener} we verify that the most restrictive terms are bounded. Indeed,
	$$\|\vr^n_N\Grad\vt_N\|_{L^{p}((0,T)\times\Omega)}\leq c\sumkN\|\sqrt{\vr_{N}^n}\Grad\log\vt_N\|_{L^2((0,T)\times\Omega)}\|\erkN\|^{\frac{1}{2}}_{L^{\frac{5}{3}}((0,T)\times\Omega)}\|\vt_N\|_{L^{32/3}((0,T)\times\Omega)}\leq c,$$
with  $p>1$ on account of \eqref{eqT}, further
	$$\|{\vr_N}\Grad{\vt_N}\|_{L^2(0,T;L^{\frac{3}{2}}(\Omega))}\leq c\|\sqrt{\vr_N}\Grad{\vt_N}\|_{L^{2}((0,T)\times\Omega)}\|\sqrt{\vr_N}\|_{L^{\infty}(0,T;L^6(\Omega)}\leq c,$$
	\begin{equation*}
	\|\sqrt{\vr_N}\Grad{\vt_N}\sqrt{\vr_N}\vt_N^{2}\|_{L^{\frac{2B}{B+4}}(0,T;L^{\frac{3B}{2B+2}}(\Omega))}\leq \|\sqrt{\vr_N}\|_{L^\infty((0,T)\times \Omega)} \|\sqrt{\vr_N}\Grad{\vt_N}\|_{L^2(0,T;L^{\frac{3}{2}}(\Omega))}\|\vt_N\|^2_{L^B(0,T;L^{3B}(\Omega))}\leq c.
	\end{equation*}
Finally, since $B\geq8$, $\vt^{B+1}$ can be bounded using \eqref{eqT}. For $I_3$, we have that
$$\|\vt_N\vf_{k,N}\|_{L^{32/27}((0,T)\times\Omega)}\leq\|\vt_N\|_{L^{32/3}((0,T)\times\Omega)}\|\vf_{k,N}\|_{L^{4/3}((0,T)\times\Omega)}\leq c,$$
$$\|\vt_N\Grad r_{k,N}\|_{L^2(0,T;L^{\frac{4}{3}}(\Omega))}\leq\|\Grad r_{k,N}\|_{L^{2}((0,T)\times\Omega)}\|\vt_N\|_{L^\infty(0,T;L^4(\Omega))}\leq c,$$
	$$\|\vt_N\Grad\erkN\|_{L^{160/143}((0,T)\times\Omega)}\leq \|\Grad\erkN\|_{L^{\frac{5}{4}}((0,T)\times\Omega)}\|\vt_N\|_{L^{32/3}((0,T)\times\Omega)}\leq c.$$
As a conclusion we have  that
	\begin{equation}\label{ptvtN}
	\ptb{\vr_N\vt_N+\beta\vt_N^4}\in L^1(0,T;W^{-1,p}(\Omega))\cup L^p(0,T;W^{-2,q}(\Omega)),
	\end{equation}
for some $p,q>1$. On the other hand, since $\pt\vr$ is uniformly bounded in $L^2(0,T; W^{2s,2}(\Omega))$, $\vr>c(\lambda)$ and $\vt>0$, we have 
	$$\|\pt\vt_N\|_{ L^1(0,T;W^{-1,p}(\Omega))\cup L^p(0,T;W^{-2,q}(\Omega)))}\leq c\|\ptb{\vr_N\vt_N+\beta\vt_N^4}\|_{ L^1(0,T;W^{-1,p}(\Omega))\cup L^p(0,T;W^{-2,q}(\Omega)))},$$
thus an application of the Aubin-Lions lemma gives precompactness of the sequence approximating the temperature and we have
	$$\vt_N\to\vt\quad\mbox{strongly\ in\ }L^{p'}((0,T)\times \Omega)$$
for any $1\leq p'<32/3$.

\subsubsection{Passage to the limit in the momentum equation}
Having the strong convergence of the density, we  start to identify the limit for $N\to \infty$ in the nonlinear  terms of the momentum equation.

{\bf The convective term.} First, one observes that
	$$\vr_N\vu_N\to\vr\vu \quad\mbox{weakly$^*$\ in} \ L^\infty(0,T;L^{2}(\Omega)),$$
due to the uniform estimates \eqref{ess_ener} and the strong convergence of the density. 
Next, one can show that for any $\phi\in \cup_{n=1}^\infty X_N$ the family of functions $\intO{\vr_N\vu_N(t)\phi}$ is bounded and equi-continuous in $C([0,T])$, thus via the Arzel\`a-Ascoli theorem and density of smooth functions in $L^2(\Omega)$ we get that
	\begin{equation}\label{weak_mom}
	\vr_N\vu_N\to\vr\vu \quad\mbox{in} \ C([0,T];L_{\rm{weak}}^{2}(\Omega)).
	\end{equation}
Finally, by the compact embedding $L^2(\Omega)\subset W^{-1,2}(\Omega)$ and the weak convergence of $\vu_N$ (cf. \eqref{conv_vu}) we verify that
	$$\vr_N\vu_N\otimes\vu_N\to\vr\vu\otimes\vu\quad\mbox{weakly\ in} \ L^2((0,T)\times\Omega).$$

{\bf The capillarity term.} We rewrite it in the form 
	$$\intTO{\vr_N\Grad\lap^{2s+1}\vr_N\cdot\vcg{\phi}}=\intTO{\lap^s\Div\left({\vr_N\vcg{\phi}}\right)\lap^{s+1}\vr_N}.$$
Due to \eqref{conv_vr} and boundedness of the time derivative of $\vr_N$, we infer that
	\begin{equation}\label{conv_vr_strong}
	\vr_N\to\vr \quad\mbox{strongly\ in} \ L^2(0,T;W^{2s+1,2}(\Omega)),
	\end{equation}
thus
	$$
	\intTO{\lap^s\Div(\vr_N\vcg{\phi})\lap^{s+1}\vr_N} \to \intTO{\lap^s\Div(\vr \vcg{\phi})\lap^{s+1}\vr}$$
for any $\vcg{\phi} \in C^\infty((0,T)\times \overline{\Omega})$.

{\bf The momentum term.} We rewrite it in the form
	$$-\lambda\intTO{\vr_N\lap^{2s+1}(\vr_N\vu_N)\cdot\vcg{\phi}}=\lambda\intTO{\lap^{s}\Grad(\vr_N\vu_N) :\lap^{s}\Grad(\vr_N\vcg{\phi})}$$
so the convergences established in  \eqref{conv_vu} and \eqref{conv_vr_strong} are sufficient to pass to the limit here.

{\bf The molecular part of the pressure.} Passage to the limit here requires for example the weak convergence of the species densities and strong convergence of the temperature, which is guaranteed on account of the results from previous sections.\\

Strong convergence of the density, temperature together with the strong convergence of the species densities enables us to  perform the limit passage in the momentum equation \eqref{FG1} satisfied 
for any function $\vcg{\phi}\in C^1([0,T]; X_{\tilde N})$ such that $\vcg{\phi}(T)=\vc{0}$ and by the density argument we can take all such test functions from $C^1([0,T]; W^{2s+1}(\Omega))$.

\subsubsection{Passage to the limit in the species equations}

As already mentioned this passage differs from what was done in the isothermal case studied in \cite{MPZ} only due to a presence of $\vt$ in the form of diffusion fluxes $\vf_k$. However, on account of strong convergence of the temperature and species densities we can write
	\begin{equation}
	\vf_{k,N}=-\frac{C_{0}(\vt_N,\vr_N)}{\sumlN\erlN}\sumlN C_{kl}\erlN\left(\Grad r_{l,N}+\Grad\log \vt_N\right)\to\vf_k\quad\mbox{weakly\ in\ }
	L^{\frac{4}{3}}((0,T)\times\Omega).
	\end{equation}

\subsubsection{Passage to the limit in the internal energy balance equation}
Passage to the limit in the terms $2\vr^n|\D(\vu)|^2$, $\lambda|\lap^{s}\Grad(\vr\vu)|^{2}$, and $\lambda\ep|\lap^{s+1}\vr|^{2}$ requires a sort of strong convergence of these quantities. 
This will be deduced from the kinetic energy balance. For this purpose we need to show that $\vu$ can be used as a test function in the limit momentum equation. Here it is again important that we have 
the kinetic energy estimate \eqref{kinetic}. 
Indeed,  in \eqref{FG1} all terms are bounded
 due to estimates above.
Moreover, thanks to the lower bound of $\vr$ we can verify that $\vu$ is actually a continuous function with respect to time and that it is continuously differentiable. To see this it is enough to differentiate \eqref{FG1} with respect to time and use the kinetic energy balance.

Now, using $\vu$ as a test function and taking advantage of the fact that the limit continuity equation is satisfied pointwisely, we obtain
	\begin{multline}\label{kinetic_po}
	\inttOB{2\vr^n|\D(\vu)|^2
	+\lambda|\lap^{s}\Grad(\vr\vu)|^{2}+\lambda\ep|\lap^{s+1}\vr|^{2}}+\intO{\lr{{1\over2}\vr|\vu|^2+\frac{\lambda}{2}|\Grad^{2s+1}\vr|^{2}}(t)}\\
	=\inttO{{\pi}\Div\vu}+\intO{\lr{{1\over2}\vr|\vu|^2+\frac{\lambda}{2}|\Grad^{2s+1}\vr|^{2}}(0)}
	\end{multline}
for any $t\in[0,T]$.
On the other hand, due to \eqref{firstep}, we have
	\begin{equation}\label{kinetic_przed}
	\begin{split}
	&\lim_{N\to\infty}\inttOB{2\vr_{n,N}|\D(\vu_N)|^2
	+\lambda|\lap^{s}\Grad(\vr_N\vu_N)|^{2}+\lambda\ep|\lap^{s+1}\vr_N|^{2}}\\
	&+\lim_{N\to\infty}\intO{\lr{{1\over2}\vr_N|\vu_N|^2+\frac{\lambda}{2}|\Grad^{2s+1}\vr_N|^{2}}(t)}=
	\inttO{{\pi}\Div\vu}+\intO{\lr{{1\over2}\vr|\vu|^2+\frac{\lambda}{2}|\Grad^{2s+1}\vr|^{2}}(0)}.
	\end{split}
	\end{equation}
The comparison of these two expressions yields
	$$\|\sqrt{\vr_{n,N}}\D\vu_N\|^2_{L^2((0,T)\times\Omega)}\to \|\sqrt{\vr_{n}}\D\vu\|^2_{L^2((0,T)\times\Omega)},$$
	$$\|\lap^{s+1}\vr_N\|^2_{L^2((0,T)\times\Omega)}\to \|\lap^{s+1}\vr\|^2_{L^2((0,T)\times\Omega)},$$
	$$\|\sqrt{\lambda}\lap^{s}\Grad(\vr_N\vu_N)\|^2_{L^2((0,T)\times\Omega)}\to \|\sqrt{\lambda}\lap^{s}\Grad(\vr\vu)\|^2_{L^2((0,T)\times\Omega)}$$
and for all $t\in[0,T]$ we have that
	$$\|{\vr_N}|\vu_N|^2(t)\|_{L^1(\Omega)}\to \|{\vr}|\vu|^2(t)\|_{L^1(\Omega)},$$
	$$\|\Grad^{2s+1}\vr_N(t)\|^2_{L^2(\Omega)}\to \|\Grad^{2s+1}\vr(t)\|_{L^2(\Omega)}.$$
Having convergences of these norms and the relevant weakly convergent sequences we deduce the strong convergence. On account of that we are able to perform the limit passage in the internal energy equation \eqref{app_vt}
		 \begin{equation}\label{app_vtN}
		 \begin{split}
	&\intTO{\lr{\vr\vt+\beta\vt^4}\pt\phi} +\intTO{\vu\lr{\vr\vt+\vt^4}\cdot\Grad\phi}-\intTO{\kappa_{\ep}\Grad\vt\cdot\Grad\phi}\\
	+&\sumkN \intTO{ \left(\vt\frac{\vf_{k}}{m_k}-\delta \vt\Grad r_{k}-\ep\vt\Grad\erk\right)\cdot\Grad\phi} \\
	=-&\intTO{\lr{\frac{\ep}{\vt^2}-\ep\vt^{5}}\phi}+\intTO{\lr{\pi_m+\frac{\beta}{3}\vt^4}\Div\vu\phi}\\
	-&\intTO{\lr{2\vr^n|\D(\vu)|^2
	+\lambda|\lap^{s}\Grad(\vr\vu)|^{2}+\lambda\ep|\lap^{s+1}\vr|^{2}+\ep\frac{1}{\vr}\frac{\partial\pi_c(\vr)}{\partial\vr}|\Grad\vr|^2}\phi}\\
	-&\intO{\lr{\vr\vt+\beta\vt^4}(0)\phi(0)},
	\end{split}
	\end{equation}
for any smooth $\phi$ vanishing at $t=T$, where the limit of the heat flux term should be always understood in the following sense
\begin{equation}\label{weak_kappa1}
\intTO{\kappa_{\ep}\Grad\vt\cdot\Grad\phi}=\intTO{\lr{\kappa_0+\ep\vr^n+\vr+\vr\vt^2}\Grad\vt\cdot\Grad\phi}-\frac{\beta}{B+1}\intTO{\vt^{B+1}\lap\phi}.
\end{equation}

\subsubsection{Limit in the total energy balance equation}

Now we use  $\vu\phi$ as a test function in the limit momentum equation  \eqref{FG1}, using again the limit continuity equation and after integrating by parts we get 

	\begin{equation}\label{kinetic_weak}
	\begin{split}
	&\intTO{\vr\frac{|\vu|^2}{2}\pt{\phi}}+\intTO{\vr\frac{|\vu|^2}{2}\vu\cdot\Grad\phi}
	-\intTO{(2\vr^n\D(\vu)\vu-\pi\vu)\cdot\Grad\phi}\\
	=&\intTO{2\vr^n|\D(\vu)|^2\phi}+\lambda\intTO{\lap^{s}\Grad(\vr\vu) :\lap^{s}\Grad (\vr\vu\phi)}+\frac{\ep}{2} \intTO {|\vu|^2 \nabla \vr \cdot \nabla \phi}\\
	&-\intTO{\pi\Div \vu\phi}+\lambda\intTO{\lap^s\Div\left({\vr\vu\phi}\right)\lap^{s+1}\vr}-\intO{\vr\frac{|\vu|^2}{2}(0){\phi}(0)}. 
	\end{split}
	\end{equation}
We  apply to the approximate continuity equation the operator $\lap^s$ and then test it by $\lambda\Div(\Grad\lap^s\vr\phi)$ in order to obtain	
\begin{equation}\label{cont:diff}
\begin{split}
&\intTO{\frac{\lambda}{2}|\Grad\lap^s\vr|^2\pt\phi}+\lambda\intTO{\lap^s\Div(\vr\vu)\lap^{s+1}\vr\phi}+\lambda\intTO{\lap^s\Div(\vr\vu)\Grad\lap^s\vr\cdot\Grad\phi}\\
&-\lambda\ep\intTO{|\lap^{s+1}\vr|^2\phi}-\lambda\ep\intTO{\lap^{s+1}\vr\Grad\lap^s\vr\cdot\Grad\phi} + \frac{\lambda}{2} \intO{|\Grad \lap^s \vr|^2(0) \phi (0)} =0.
\end{split}
\end{equation}
Now, summing \eqref{app_vtN} with \eqref{kinetic_weak} and \eqref{cont:diff} and using the limit continuity equation to rewrite the term $\intTO{\pi_c\Div \vu\phi}$, we get the weak formulation of the total energy plus some terms which will disappear in the subsequent limit passages
	\begin{equation}\label{totEN1}
	\begin{split}
	&\intTO{\lr{\vr e+\frac{1}{2}\vr|\vu|^2+\frac{\lambda}{2}|\Grad^{2s+1}\vr|^2}\pt\phi} +\intTO{\lr{\vu\vr e+\frac{1}{2}\vr|\vu|^2\vu}\cdot\Grad\phi}\\
	-&\intTO{\kappa_{\ep}\Grad\vt\cdot\Grad\phi}+\sumkN \intTO{ \left(\vt\frac{\vf_{k}}{m_k}-\delta \vt\Grad r_{k}-\ep\vt\Grad\erk\right)\cdot\Grad\phi}\\
	+&\intTO{\pi\vu\cdot\Grad\phi}
	-\intTO{(2\vr^n\D(\vu)\vu)\cdot\Grad\phi}=-\intTO{\lr{\frac{\ep}{\vt^2}-\ep\vt^{5}}\phi}\\
	+&\intTO{R_{\ep,\lambda}(\vr,\vt,\vu,\phi)}-\intO{\lr{\vr e(0)+\frac{1}{2}\vr|\vu|^2(0)+\frac{1}{2}|\Grad^{2s+1}\vr|^2(0)}\phi(0)},
	\end{split}
	\end{equation}
with \eqref{weak_kappa1} and
\begin{equation}\label{RN1} 
\begin{split}
{R_{\ep,\lambda}(\vr,\vt,\vu,\phi)} = 
&\lambda\left[\Delta^s ( \Div (\vr \vu\phi )) \Delta^{s+1} \vr 
-\lap^s\Div(\vr\vu)\lap^{s+1}\vr\phi\right]\\
&-\lambda\lap^s\Div(\vr\vu)\Grad\lap^s\vr\cdot\Grad\phi
- \lambda\left[|\Delta^s (\nabla (\vr \vu))|^2 \phi- \Delta^s \nabla (\vr \vu): \Delta^s \nabla (\vr \vu \phi)\right] \\
&+ 
\lambda \ep\lap^{s+1}\vr\Grad \lap^{s}\vr \cdot\Grad\phi + \frac{\ep}{2} |\vu|^2 \nabla \vr \cdot \nabla\phi +\ep \nabla \vr \cdot \nabla \phi \left(e_c(\vr) + \frac{\pi_c(\vr)}{\vr}\right).
\end{split}
	\end{equation}
 Finally we  define partial densities in the following way
	$$\vr_{k}=m_k\mbox{e}^{r_{k}},\quad k=1,\ldots,n$$ which finishes the proof of Theorem \ref{t 2}. $\Box$

\section{Passage to the limit $\delta\to0$}\label{Sec:delta}

From \eqref{ess_entr} and \eqref{entr_part} we can deduce that 
	$$\vr_{k,\delta}>0\quad\mbox{a.e.\ in\ } (0,T)\times\Omega,\ k=1,\ldots,n.$$
We will concentrate here only on the strong convergence of the partial densities and on the limit passage in the species equations, since the strong convergence of $\vt_\delta,\vr_\delta$ and $\vu_\delta$ to $\vt,\vr$ and $\vu$, respectively, can be proven identically as in the previous  section. 

Passage to the limit in \eqref{rk_del} follows the same steps as in \cite{MPZ}, Chapter III.C. Indeed, repeating procedure leading to \eqref{time_part}, we can show that
	\begin{equation*}
	\|\pt \vr_{k,\delta}\|_{L^{\frac{5}{4}}(0,T;W^{-1,\frac{5}{4}}(\Omega))}
	\leq c
	\end{equation*}
and thus, the uniform estimates from Section \ref{UB} and the application of the Aubin-Lions lemma give rise to the following convergences
	\begin{equation*}
	\begin{array}{c}
	\delta \log\frac{\varrho_{k,\delta}}{m_k}\to 0\quad\mbox{strongly \ in}\ L^{\infty}(0,T;L^2(\Omega)),\\
	{\delta}\Grad \log\frac{\varrho_{k,\delta}}{m_k}\to 0\quad\mbox{strongly \ in}\ L^{2}((0,T)\times\Omega),\\
	\varrho_{k,\delta}\to \varrho_k\quad\mbox{strongly\ in\ }L^{q}((0,T)\times\Omega), q<\frac 53,\\
	{\Grad\varrho_{k,\delta}\to\Grad \varrho_k\quad\mbox{weakly\ in\ }L^2(0,T;L\log L_{\rm{weak}^*}(\Omega))\cap L^{\frac{5}{4}}((0,T)\times\Omega)},\\
	{\varrho_{k,\delta}\to \varrho_k\quad\mbox{in\ }C([0,T];L\log L_{\rm{weak}^*}(\Omega))}.
	\end{array}
	\end{equation*}
Moreover,
	\begin{equation*}
	\begin{array}{c}
	\frac{\varrho_{j,\delta}}{\sumkN\varrho_{k,\delta}}\to\frac{\varrho_j}{\sumkN \varrho_k}\quad\mbox{strongly\ in\ }L^p((0,T)\times\Omega),\ p<\infty,\\
	\varrho_k\geq0\quad\mbox{a.e.\ in\ } (0,T)\times\Omega,\ k=1,\ldots,n,
	\end{array}
	\end{equation*}
and due to strong convergence $\vt_\delta\to\vt$ 
we obtain that
	$$\vf_{k,\delta}\to\vf_k \quad\mbox{weakly\ in\ }L^{\frac{4}{3}}((0,T)\times(\Omega)),$$
where $\vf_k$ depends on the limit functions $\vr,\vt, \varrho_1,\ldots,\varrho_n$, as specified in \eqref{eq:diff}. Thus letting $\delta\to 0$ in the approximate equations \eqref{rk_del}, we verify  that:
	
	\begin{lemma}
	The limit quantities $\varrho_k$, $k=1,\ldots,n$ satisfy
	 	\begin{equation}\label{1_ep}
		\pt{\varrho_k}+\Div (\varrho_{k} \vu)-\ep\lap\varrho_k+ \Div (\vf_k(\vr,\vt,\vr_1,\ldots,\vr_n))  =  \vr^n\vt\vw_{k},\quad k=1,...,n,
		\end{equation}
	in the sense of distributions on $(0,T)\times\Omega$. \
	\end{lemma}

In addition,  summing \eqref{1_ep} with respect to $k=1,\ldots,n$, property \eqref{wform0} and Remark \ref{rem1} lead to the following equation
	 \begin{equation*}\label{2_ep}
	\pt{\vr^n}+\Div (\vr^n \vu)-\ep\lap\vr^n= 0.
	\end{equation*}
	This equation is, due to the previous lemma, satisfied in the same sense as system \eqref{1_ep}, together with the initial condition $\varrho(0,x)=\varrho^0=\sumkN\varrho^0_k(x)$ for a.a. $x\in\Omega$. Moreover, it is possible to 
identify $\vr^n$ with $\vr_\ep$ -- the unique classical  {positive} solution to the initial-value problem \eqref{apcontN_2}
constructed in Lemma \ref{kropa}, {see Ref. \cite{NS}, Sections 7.6.3--7.6.7 for more details.} In particular, we know that any solution of \eqref{1_ep} satisfies
	$$\sumkN\varrho_k=\varrho\quad \mbox{a.e.\ in\ }(0,T)\times\Omega.$$
As a corollary we have that $Y_k=\frac{\vr_k}{\vr}$ satisfies
	\begin{equation}\label{maxY}
	\|Y_k\|_{L^\infty((0,T)\times\Omega)}\leq 1,\quad k=1,\ldots,n.
	\end{equation}

	\begin{rmk}Note  that at this level we already perform the limit in the total energy balance \eqref{totEN1} rather than in the internal energy balance.  
	The latter will have from now on an auxiliary character, however after the limit passage $\delta\to0$ we will have to replace the equality  by inequality, we only have:
		 \begin{multline}\label{vt_in_delta}
		\ptb{\vr \vt+\beta\vt^4} +\Div(\vu(\vr\vt+\beta\vt^4))-\Div\left(\kappa_{\ep}(\vr,\vt)\Grad\vt\right)
		+\sumkN \Div \left(\vt\frac{\vf_{k}}{m_k}-\ep\vt\Grad\erk\right)\geq \\
		\frac{\ep}{\vt^2}-\ep\vt^{5}-\lr{\pi_m+\frac{\beta}{3}\vt^4}\Div\vu+2\vr|\D(\vu)|^2
		+\lambda|\lap^{s}\Grad(\vr\vu)|^{2}+\lambda\ep|\lap^{s+1}\vr|^{2},
		\end{multline}
	in the sense of distributions on $(0,T)\times\Omega$.
	\end{rmk}

\section{Derivation of the B-D estimate}\label{Sec:BD}
At this level we are left with only two parameters of approximation: $\ep$ and $\lambda$. From the so-far obtained a-priori estimates only the ones following from \eqref{main_ener} and \eqref{est:entr} were independent of these parameters. However having the $\ep$-dependent estimate for $\lap^{s+1}\vr$ allows us to derive a type of B-D estimate, from which it will follow that this estimate depends only on $\lambda$. As a by-product, we will derive the energy estimate independent of  $\lambda$. Note that so far in \eqref{ener_vel} we were only able to estimate the r.h.s. using the $\lambda$-dependent bounds for $\vu$. We will prove the following lemma.

%

	\begin{lemma}\label{l10} For any positive constant $r>1$, we have
		\begin{equation}\label{cacya}
		\begin{split}
				\Dt& \intOB{{1\over2}\vr|\vu+\Grad\phi(\vr)|^2+{{r-1}\over2}\vr|\vu|^{2}+\frac{r\lambda}{2}|\Grad\lap^s\vr|^{2}+r\vr e_c(\vr)}
		+\intO{\Grad\phi(\vr)\cdot\Grad {\pi}}\\
		&\quad +{1\over2}\intO{\vr|\Grad\vu-\Grad^{T}\vu|^2} +2\lambda\intO{|\lap^{s+1}\vr|^{2}}+2(r-1)\intO{\vr|\D(\vu)|^2}\\
		&\quad +r\intOB{\lambda\ep|\lap^{s+1}\vr|^{2}+\lambda|\lap^s\Grad(\vr\vu)|^2}\\
		&\leq-\ep\intO{(\Grad\vr\cdot\Grad)\vu\cdot\Grad\phi}+\ep\intO{\lap\vr\frac{|\Grad\phi|^{2}}{2}}+\ep\intO{\vr\Grad\phi(\vr)\cdot\Grad\left(\phi'(\vr)\lap\vr\right)}\\
		&\quad-\ep\intO{\Div(\vr\vu)\phi'(\vr)\lap\vr}
		+r\intO{\lr{\pi_m+\frac{\beta}{3}\vt^4}\Div\vu}-2\lambda\intO{\lap^s\Grad(\vr\vu):\lap^s\Grad^2\vr}		
%
		\end{split}
		\end{equation}
	in ${\cal D}'(0,T)$, where $\Grad\phi(\vr)=2\Grad\log\vr$, $e_c(\vr)=\int_{0}^\vr{{y^{-2}}{\pi}_{c}(y)}~\mbox{d}y\geq0$.
	\end{lemma}

\pf The basic idea of the proof is to find the explicit form of the term:
	\begin{equation}\label{calka}
	\Dt \intOB{{1\over2}\vr|\vu|^{2}+\vr\vu\cdot\Grad\phi(\vr)+{1\over 2}\vr|\Grad\phi(\vr)|^2}.
	\end{equation}
The first term can be evaluated by means of the main energy equality,  i.e.
	\begin{multline}\label{kinBD}
	\Dt\intOB{{1\over2}\vr|\vu|^2+\frac{\lambda}{2}|\Grad\lap^{s}\vr|^{2}+\vr e_c(\vr)}
	+\ep\intO{\frac{1}{\vr}\frac{\partial\pi_c}{\partial\vr}|\Grad\vr|^2}\\
	+\intTOB{2\vr|\D(\vu)|^2
	+\lambda|\lap^{s}\Grad(\vr\vu)|^{2}+\lambda\ep|\lap^{s+1}\vr|^{2}}
	=\intO{\lr{\pi_m+\frac{\beta}{3}\vt^4}\Div\vu}.
\end{multline}
To get a relevant expression for the third term in \eqref{calka}, we multiply the approximate continuity equation by $\frac{|\Grad\phi(\vr)|^2}{2}$ and we obtain the following sequence of equalities
	\begin{equation}\label{39ep}
	\begin{split}
	&\Dt \intO{{1\over2}\vr|\Grad\phi(\vr)|^2}\\
	&=\intOB{\vr\pt\frac{|\Grad\phi(\vr)|^2}{2}-\frac{|\Grad\phi(\vr)|^2}{2}\Div(\vr\vu)+\ep\frac{|\Grad\phi(\vr)|^2}{2}\lap\vr}\\
	&=\intOB{\vr\Grad\phi(\vr)\cdot\Grad\left(\phi'(\vr)\pt\vr\right)-\frac{|\Grad\phi(\vr)|^2}{2}\Div(\vr\vu)+\ep\frac{|\Grad\phi(\vr)|^2}{2}\lap\vr}.
	\end{split}
	\end{equation}
Using the approximate continuity equation, we get	
	\begin{equation}
	\begin{split}
	&\intO{\vr\Grad\phi(\vr)\cdot\Grad\left(\phi'(\vr)\pt\vr\right)}\\
	&=\intO{\ep\vr\Grad\phi(\vr)\cdot\Grad\left(\phi'(\vr)\lap\vr\right)}-\intO{\vr\Grad\vu:\Grad\phi(\vr)\otimes\Grad\phi(\vr)}\\
	&-\intO{\vr\Grad\phi(\vr)\cdot\Grad\left(\phi'(\vr)\vr\Div\vu\right)}-\intO{\vr\vu\otimes\Grad\phi(\vr):\Grad^{2}\phi(\vr)}.\\
		\end{split}
	\end{equation}
Integrating by parts the two last terms from the r.h.s.	
		\begin{equation}
	\begin{split}
	&\intO{\vr\Grad\phi(\vr)\cdot\Grad\left(\phi'(\vr)\pt\vr\right)}=\intO{\ep\vr\Grad\phi(\vr)\cdot\Grad\left(\phi'(\vr)\lap\vr\right)}-\intO{\vr\Grad\vu:\Grad\phi(\vr)\otimes\Grad\phi(\vr)}\\
	&\quad+\intO{\vr\left|\Grad\phi(\vr)\right|^2\Div\vu}+\intO{\vr^2\phi'(\vr)\lap\phi(\vr)\Div\vu}\\
	&\quad+\intO{\left|\Grad\phi(\vr)\right|^2\Div(\vu\vr)}+\intO{\vr\vu\cdot\Grad(\Grad\phi(\vr))\cdot\Grad\phi(\vr)}.
	\end{split}
	\end{equation}
Combining the three previous equalities we finally obtain
	
		\begin{equation}
	\begin{split}
	&\Dt \intO{{1\over2}\vr|\Grad\phi(\vr)|^2}\\
	&=\intO{\ep\vr\Grad\phi(\vr)\cdot\Grad\left(\phi'(\vr)\lap\vr\right)}-\intO{\vr\Grad\vu:\Grad\phi(\vr)\otimes\Grad\phi(\vr)}+\intO{\ep\frac{|\Grad\phi(\vr)|^2}{2}\lap\vr}\\
	&\quad+\intO{\vr^2\phi'(\vr)\lap\phi(\vr)\Div\vu}+\intO{\vr\left|\Grad\phi(\vr)\right|^2\Div\vu}.
	\end{split}
	\end{equation}	
In the above series of equalities, each one holds pointwisely with respect to time due to the regularity of $\vr$ and $\Grad\phi$. This is not the case of the middle integrant of \eqref{calka}, for which one should really think of weak in time formulation. Denote 
$$V=W^{2s+1,2}(\Omega),\quad
\mbox{and}\quad \vc{v}=\vr\vu,\quad \vc{h}=\Grad\phi.$$ 
We know that $\vc{v}\in L^2(0,T;V)$ and its weak derivative with respect to time variable $\vc{v}'\in L^2(0,T;V^*)$, where $V^*$ denotes the dual space to $V$. Moreover, $\vc{h}\in L^2(0,T; V)$, $\vc{h}'\in L^2(0,T;W^{2s-1,2}(\Omega))$. Now, let $\vc{v}_m$, $\vc{h}_m$ denote the standard mollifications in time of $\vc{v}$ and $\vc{h}$ respectively. By the properties of mollifiers we know that
$$\vc{v}_m, \vc{v}'_m\in C^\infty(0,T;V), \quad \vc{h}_m,\vc{h}'_m\in C^\infty(0,T;V),$$
and
\begin{equation}\label{univh}
\begin{array}{ll}
\vspace{0.2cm}
\vc{v}_m\to\vc{v}\quad \mbox{in} \ L^2(0,T;V), &\vc{h}_m\to\vc{h}\quad \mbox{in} \ L^2(0,T;V),\\
\vc{v}'_m\to\vc{v}'\quad \mbox{in} \ L^2(0,T;V^*),\quad &\vc{h}'_m\to\vc{h}'\quad \mbox{in} \ L^2(0,T;V^*).
\end{array}
\end{equation}
For these regularized sequences we may write
\begin{equation}\label{mvh}
\Dt\intO{\vc{v}_m\cdot\vc{h}_m}=\Dt(\vc{v}_m,\vc{h}_m)_{V}=(\vc{v}'_m,\vc{h}_m)_{V}+(\vc{v}_m,\vc{h}'_m)_{V}.
\end{equation}
Using the Riesz representation theorem we verify that $\vc{v}'_m\in C^\infty(0,T;V)$ uniquely determines the functional $\Phi_{\vc{v}'_m}\in V^*$ such that $(\vc{v}'_m,\psi)_{V}= \langle{\Phi_{\vc{v}'_m},\psi\rangle}_{V^*,V}=\intO{\vc{v}'_m\cdot \psi},$ $\forall\psi\in V$; for the second term from the r.h.s. of  \eqref{mvh} we can simply replace $V=L^2(\Omega)$
%
and thus we obtain
$$-\intT{(\vc{v}_m,\vc{h}_m)_{V}\psi'}=\intT{\langle{\vc{v}'_m,\vc{h}_m\rangle}_{V^*,V}\psi}+\intT{(\vc{v}_m,\vc{h}'_m)_{L^2(\Omega)}\psi} \quad\forall\psi\in {\cal D}(0,T).$$
Observe that both integrands from the r.h.s. are uniformly bounded in $L^1(0,T)$, thus, using \eqref{univh}, we let $m\to \infty$ to obtain
\begin{equation*}
\Dt(\vc{v},\vc{h})_{V}=\langle{\vc{v}',\vc{h}\rangle}_{V^*,V}+(\vc{v}\cdot\vc{h}')_{L^2(\Omega)}\quad \mbox{in}\ {\cal D}'(0,T).
\end{equation*}
Coming back to our original notation, this means that the operation 
\begin{equation}\label{crossep1}
\Dt\intO{\vr\vu\cdot\Grad\phi(\vr)}=\langle{\pt(\vr\vu),\Grad\phi\rangle}_{V^*,V}+\intO{\vr\vu\cdot\pt\Grad\phi}
\end{equation}
is well defined and is nothing but equality between two scalar distributions.
By the fact that $\pt\Grad{\phi}$ exists a.e. in $(0,T)\times\Omega$ we may use approximate continuity equation to write
\begin{equation}\label{crossep}
\intO{\vr\vu\cdot\pt\Grad\phi}=\intO{(\Div(\vr\vu))^2 \phi'(\vr)}-\ep\intO{\Div(\vr\vu)\phi'(\vr)\lap\vr},
\end{equation}
whence the first term on the r.h.s. of \eqref{crossep1} may be evaluated by testing the approximate momentum equation by $\Grad\phi(\vr)$
	\begin{equation}\label{secondaep}
	\begin{split}
	&\langle{\pt(\vr\vu),\Grad\phi\rangle}_{V^*,V}\\
	&\quad=-\intO{2\vr\lap\phi(\vr)\Div\vu}+2\intO{\Grad\vu:\Grad\phi(\vr)\otimes\Grad \vr}
	-2\intO{\Grad\phi(\vr)\cdot\Grad \vr\Div\vu}\\
	&\qquad-\intO{\Grad\phi(\vr)\cdot\Grad {\pi}}-\lambda\intO{\lap^{s+1}\vr\lap^s\Div(\vr\Grad\phi(\vr))}
	-\lambda\intO{\lap^{s}\Grad(\vr\vu) :\lap^{s}\Grad (\vr\Grad{\phi})}\\
	&\qquad-\intO{\Grad\phi(\vr)\cdot\Div(\vr\vu\otimes\vu)}
 	-\ep\intO{(\Grad\vr\cdot\Grad)\vu\cdot\Grad\phi(\vr)}.
	\end{split}
	\end{equation}
Recalling the form of $\phi(\vr)$ it can be deduced that the combination of (\ref{39ep}) with (\ref{crossep1}--\ref{secondaep}) yields
	\begin{multline}\label{secondep}
	\begin{split}
	\Dt&\intOB{\vr\vu\cdot\Grad\phi(\vr)+{1\over2}\vr|\Grad\phi(\vr)|^2}+\intO{\Grad \pi\cdot\Grad\phi(\vr)}+2\lambda\intO{|\lap^{s+1}\vr|^{2}}\\
	&=-\intO{\Grad\phi(\vr)\Div(\vr\vu\otimes\vu)}+\intO{(\Div(\vr\vu))^2 \phi'(\vr)}-2\lambda\intO{\lap^s\Grad(\vr\vu):\lap^s\Grad^2\vr}\\
	&\quad-\ep\intO{\Div(\vr\vu)\phi'(\vr)\lap\vr}
	+\ep\intO{\frac{|\Grad\phi(\vr)|^2}{2}\lap\vr}\\
	&\quad-\ep\intO{(\Grad\vr\cdot\Grad)\vu\cdot\Grad\phi(\vr)}+\ep\intO{\vr\Grad\phi(\vr)\cdot\Grad\left(\phi'(\vr)\lap\vr\right)}.
	\end{split}
	\end{multline}
The first two terms from the r.h.s of (\ref{secondep}) can be transformed into
	\begin{equation*}
	\begin{split}
	&\intO{\left[(\Div(\vr\vu))^2 \phi'(\vr)-\Grad\phi(\vr)\Div(\vr\vu\otimes\vu)\right]}\\
	=&\intOB{\vr^2\phi'(\vr)(\Div\vu)^2+\vr\phi'(\vr)\vu\cdot\Grad\vr\Div\vu-\vr\phi'(\vr)\Grad\vr(\vu\cdot\Grad\vu)}\\
	=&2\intOB{\vr(\Div\vu)^2-\vr(\Div\vu)^2-\vr\vu\cdot\Grad\Div\vu+\vr\partial_iu_j\partial_j u_i+\vr\vu\Grad\Div\vu}\\
	=&2\intO{\vr\partial_iu_j\partial_j u_i}=2\intO{\vr|\D\vu|^2}-2\intO{\vr\lr{\frac{\partial_iu_j-\partial_j u_i}{2}}^2}
	\end{split}
	\end{equation*}
and thus, the assertion of Lemma \ref{l10} follows by adding \eqref{kinBD} multiplied by $r$ to \eqref{secondep}. 
$\Box$

In order to deduce the uniform estimates from \eqref{cacya} we need to control all the non-positive contributions to the l.h.s. as well as the terms from the r.h.s. The $\ep$-dependent terms can be bounded similarly as in \cite{MPZ}, so we focus only on the new aspects. To this purpose we first derive the uniform bounds for partial pressures. Denoting
	\begin{equation}\label{denot}
	\vC\Grad_{x_{i}}{p}=(\Grad_{x_{i}}{p})^{I},
	\end{equation}
where 
	\begin{equation*}
	{p}=\left(\begin{array}{c}
		p_1\\
		\vdots\\
		p_n
		\end{array}\right)\quad\mbox{and} \quad
	\Grad{p}=\left(\begin{array}{c}
		\Grad p_1\\
		\vdots\\
		\Grad p_n
		\end{array}\right),
	\end{equation*}
we obtain, for every $k$-th coordinate $k\in\{1,\ldots,n\}$ and every $i$-th space coordinate $i\in\{1,2,3\}$, the following decomposition
    	\begin{equation}\label{uu}
        (\Grad_{x_{i}}  p)_{k}=(\Grad_{x_{i}} {p})^{I}_{k}+\alpha_{i} Y_k.
   	 \end{equation}
Next, multiplying the above expression by $m_{k}$ and summing over $k\in\{1,\ldots,n\}$ one gets
    	$$\alpha_{i}=\frac{\Grad_{x_{i}}(\vr\vt)}{\sum_{k=1}^{n}m_{k} Y_k}-\frac{\sum_{k=1}^{n}m_{k}(\Grad_{x_{i}}{p})^{I}_{k}}{\sum_{k=1}^{n}m_{k} Y_k }.$$
Returning to \eqref{uu} we can express the full gradients of partial pressures in terms of gradients of temperature, density and the gradient of ''known'' part of the pressure
    	\begin{equation}\label{decomp}
        \Grad{p} = 
        ( \Grad{p})^{I} + \left(\frac{\Grad(\vr\vt)}{\sum_{k=1}^{n}m_{k} Y_k}-\frac{\sum_{k=1}^{n}m_{k}(\Grad{p})^{I}_{k}}{\sum_{k=1}^{n}m_{k} Y_k }\right)  Y.
    	\end{equation}
As was announced, we will use the above expression in order to control the molecular part of the pressure from the l.h.s. of  \eqref{calka}.

\bigskip

 {\bf {Estimate of $\Grad \pi\cdot\Grad\phi$.} }
Since $\Grad\phi=2\Grad\log\vr$ and due to \eqref{decomp} we obtain
   	\begin{equation}\label{rs}
        \Grad\phi(\vr)\cdot\Grad \pi=2\pi'_{c}(\vr)\frac{|\Grad\vr|^2}{\vr}+
       2 \frac{\Grad\vr\cdot\Grad\pi_{m}}{\vr}+\frac{2\beta}{3}\frac{\Grad\vt^4\cdot\Grad\vr}{\vr}.
    	\end{equation}
The first  term is non-negative due to \eqref{coldp}, so it can be considered on the l.h.s. of \eqref{cacya} and we only need to estimate the second and the third one. Since 
$\Grad\pi_{m}=\sumkN(\Grad p)_{k}$ {and} $\sumkN ({Y})_k=1,$
we may use \eqref{decomp} to write
    	\begin{multline}\label{pepi0}
        \intO{\frac{\Grad\vr\cdot\Grad\pi_{m}}{\vr}}
        =\intO{\frac{\Grad\vr\cdot\sumkN(\Grad p)^{I}_{k}}{\vr}}+\intO{\frac{|\Grad\vr|^2 \vt}{\sumkN\vr_{k}m_{k}}}\\
        +\intO{\frac{\Grad\vr\cdot\Grad\vt \vr}{\sumkN\vr_{k}m_{k}}}-
        \intO{\frac{\Grad\vr\cdot\sumkN m_{k}(\Grad p)^{I}_{k}}{\sumkN\vr_{k}m_{k}}}=\sum_{i=1}^4 I_i.
    	\end{multline}
Note that $I_2$ is non-negative, so we can put it to the l.h.s. of \eqref{cacya} as well.\\
Next, $I_1$ and $I_4$ can be estimated in a similar way,  we have

	\begin{equation}\label{ltop}
	\intO{\frac{|\Grad\vr||\sumkN( \Grad{p})_k^{I}|}{\vr}} \leq \epsilon\intO{\frac{|\Grad\vr |^{2}\vt}{\vr}}+ c(\epsilon)\intO{\frac{|\sumkN(\Grad  p)^{I}_k|^{2}}{\vt\vr}},
	\end{equation}
so for $\ep$ sufficiently small, the first term can be controlled  by $I_2$ since $0\leq\vr_k\leq\vr$, for all $k=1,\ldots,n$.
Concerning the second integral, from (\ref{estF}) we have
	\begin{equation}\label{reduced}
	\intTO{\sumkN\frac{\pi_m\vf_{k}^{2}}{\vC_0\vt\vr_{k}}}\leq c.
	\end{equation}
Using \eqref{difp}, the integral may be transformed as follows
	\begin{equation}\label{reduced2}
	\intTO{\sumkN\frac{\vC_0(\vC\Grad  p)_k^{2}}{\pi_m\vt\vr_{k}}}\leq c,
	\end{equation}
thus, due to \eqref{denot} and \eqref{assC}, the integral over time of the r.h.s. of \eqref{ltop} is bounded.

For  $I_3$ we verify that
    $$\left|\Grad\vr\cdot\Grad\vt\frac{\vr}{\sumkN\vr_{k}m_{k}}\right|\leq c(\epsilon)\vr{|\Grad\vt|^{2}}
    +\epsilon{|\Grad\sqrt{\vr}|^{2}},$$
and the first term is bounded in view of \eqref{main-temp-est} whereas boundedness of the second one follows from the Gronwall inequality applied to (\ref{cacya}). 

Finally
	$$\frac{2\beta}{3}\frac{\Grad\vt^4\cdot\Grad\vr}{\vr}\leq c(\epsilon)\beta\|\Grad\vt^4\|^2_{L^2((0,T)\times\Omega)}+\epsilon\|\Grad\log\vr\|^2_{L^2((0,T)\times\Omega)}$$
and the first term is bounded for $B\geq8$ while the second one can be estimated differently in two cases:
	\begin{itemize}
		\item[(i)] $\vr\geq 1$, then $\vr^{-1}\leq1$ and $\vr^{-2}|\Grad\vr|^2\leq\vr^{-1}|\Grad\vr|^2$ which is then bounded by the Gronwall inequality applied to (\ref{cacya})
		
		\item[(ii)] $\vr<1$, then $\vr^{-\gamma}\geq1$ and $\ep\vr^{-2}|\Grad\vr|^2\leq\vr^{-2-\gamma}|\Grad\vr|^2\leq \ep\pi'_c(\vr)\vr^{-1}|\Grad\vr|^2$
which is absorbed by the analogous term from the l.h.s. of \eqref{cacya} (the first term from \eqref{rs}).\\
\end{itemize}

{\bf Estimate of $\lr{\pi_m+\frac{\beta}{3}\vt^4}\Div\vu$.} Here we proceed in a little different way than in case of kinetic energy balance. The problem is that we want to have estimates uniform with respect to $\lambda$, so we cannot use the bound for $\Grad\vu(t)$ in $L^\infty(\Omega)$ any more.
For the molecular part of the pressure we have
	\begin{multline}\label{lasttop}
	\left|\intO{\vr_k\vt\Div\vu}\right|  \leq \epsilon\|\sqrt{\vr}\Div\vu\|^2_{L^{2}(\Omega)}+c(\epsilon)\left\|\frac{\vr\vt}{\sqrt{\vr}}\right\|^2_{L^{2}(\Omega)}\\
	\leq
	\epsilon\|\sqrt{\vr}\Div\vu\|^2_{L^{2}(\Omega)}+c(\epsilon)\|\vr \|^{\frac{1}{2}}_{L^{\frac{3}{2}}(\Omega)}\|\vt\|_{L^{6}(\Omega)}.
    	\end{multline}
On account of \eqref{maxtemp},  $\vt\in L^2(0,T; L^6(\Omega))$. Moreover, the term $\epsilon\|\sqrt{\vr}\Div\vu\|^2_{L^{2}(\Omega)}$ is then absorbed by the l.h.s. of \eqref{kinBD}. Therefore, since the Sobolev imbedding theorem implies that  $\|{\vr}\|_{L^{\frac{p}{2}}(\Omega)}\leq c\left\|\frac{\Grad\vr}{\sqrt{\vr}}\right\|_{L^{2}(\Omega)}$ {for} $1\leq p\leq 6,$ the Gronwall inequality can applied to \eqref{cacya} from which boundedness of \eqref{lasttop} follows.

The radiative term is slightly more difficult, however, we still can write
	\begin{multline*}
	\intTO{\vt^4|\Div\vu|}=\intTO{\vt^4\vr^{-1/2}|\sqrt{\vr}\Div\vu|}\\
	\leq
	\|\vt\|^4_{L^p(0,T;L^q(\Omega))}\|\vr^{-1/2}\|_{L^{2\gamma^{-}}(0,T;L^{6\gamma^{-}}(\Omega))}\|\sqrt\vr\Div\vu\|_{L^2((0,T)\times\Omega)},
	\end{multline*}
where $p=\frac{8\gamma^{-}}{\gamma^{-}-1}$, $q=\frac{24\gamma^{-}}{3\gamma^{-}-1}$. By the interpolation $\|\vt\|_{L^p(0,T;L^q(\Omega))}\leq\|\vt\|_{L^\infty(0,T;L^4(\Omega))}^{1-a}\|\vt\|_{L^G(0,T;L^{3G}(\Omega))}^a$
for $a=\frac{2}{3}$ and $G=\frac{16\gamma^-}{3(\gamma^--1)}$, where $G\leq B$ provided $\gamma^{-}\geq3$. Thus, we can estimate
	\begin{multline}\label{vt4vu}
	\intTO{\vt^4|\Div\vu|}
	\leq 
	c(\epsilon)\lr{\|\vt\|_{L^\infty(0,T;L^4(\Omega))}^{1/3}\|\vt\|_{L^B(0,T;L^{3B}(\Omega))}^{2/3}}^{2(\gamma^--1)/\gamma^-}\\
	+\epsilon\|\vr^{-\gamma^{-}/2}\|^2_{L^{2}((0,T)\times\Omega)}
	+\epsilon\|\sqrt\vr\Div\vu\|^2_{L^2((0,T)\times\Omega)},
	\end{multline}
and the two last terms are estimated by the l.h.s. of \eqref{cacya} and \eqref{kinBD}, while the boundedness of the first one follows from \eqref{ess_ener} and \eqref{maxtemp}.

\bigskip

{\bf Estimate of $\lambda{\lap^s\Grad(\vr\vu):\lap^s\Grad^2\vr}$.}  We have
	$$2\lambda\intO{|\lap^s\Grad(\vr\vu):\lap^s\Grad^2\vr|}
	\leq c\lambda\|\lap^s\Grad(\vr\vu)\|^2_{L^2((0,T)\times\Omega)}+\lambda \|\lap^{s+1}\vr\|^2_{L^2((0,T)\times\Omega)},$$
therefore for $r$ sufficiently large, such that $r\lambda^{-1}>c$ both terms are bounded by the r.h.s. of \eqref{cacya}.\\

\section{Estimates independent of $\ep,\lambda$, passage to the limit $\ep,\lambda\to0$}\label{Sec:ep_lambda}
In this section we first present the new uniform bounds arising from the estimate of B-D entropy, performed in Section \ref{Sec:BD}, and then we let the last two approximation parameters to 0. Note that the limit passage $\lambda\to0$ and $\ep\to 0$ could be done in  a single step, however, for the sake of transparency of this proof we do it separately.

We  complete the set of uniform bounds by the following ones
	\begin{equation}\label{negdens}
	\sqrt{\lambda}\left\|\lap^{s+1}\vr\right\|_{L^{2}((0,T)\times\Omega)}+\left\|\sqrt{{\vt}{\vr}^{-1}}\Grad\vr\right\|_{L^{2}((0,T)\times\Omega)}+\left\|\sqrt{\pi'_c(\vr)\vr^{-1}}\Grad\vr\right\|_{L^{2}((0,T)\times\Omega)}\leq c,
	\end{equation}
moreover
	\begin{equation}\label{posdens}
	\sqrt{\lambda}\left\|{\Grad^{2s+1}{\vr}}\right\|_{L^{\infty}(0,T;L^2(\Omega))}+\left\|{\Grad\sqrt{\vr}}\right\|_{L^{\infty}(0,T;L^2(\Omega))}\leq c.
	\end{equation}
The uniform estimates for the velocity vector field are the following ones
    	\begin{equation}\label{lambda_vu}
	\sqrt{\lambda}\|\lap^s\Grad(\vr\vu)\|_{L^{2}((0,T)\times\Omega)}+\|\sqrt{\vr}\Grad\vu\|_{L^{2}((0,T)\times\Omega)}+\left\|{\sqrt{\vr\vt^{-1}}}\Grad\vu\right\|_{L^{2}((0,T)\times\Omega)}\leq c
   	 \end{equation}
and the constants from the r.h.s. are independent of $\ep$ and $\lambda$.\\

{\bf Estimates of species densities.} Finally, we  take advantage of the entropy estimate \eqref{est:entr} which together with \eqref{posdens} may be used to deduce boundedness of gradients of all species densities. 
	
	\begin{lemma}
	We have
		\begin{equation}\label{gradcomp}
		\left\|\sqrt{{1+\vt}}\Grad
		\sqrt{\vr_{k}}\right\|_{L^{2}((0,T)\times\Omega)}\leq c.
		\end{equation}
	\end{lemma}

\pf
First, using the form of matrix $C$ we may write
	$${\frac{\pi_m\vf_k^2}{\vC_0\vt\vr_k}}=\frac{\vC_0|\Grad p_k|^2}{\pi_m\vr_k\vt}-2\frac{Y_k\vC_0\Grad p_k\cdot\Grad\pi_m}{\pi_m\vr_k\vt}+\frac{Y_k^2\vC_0|\Grad \pi_m|^2}{\pi_m\vr_k\vt},$$
which is bounded in $L^1((0,T)\times\Omega)$ on account of \eqref{reduced}. Clearly,
	\begin{equation}\label{almostY}
	\intTO{\frac{\vC_0|\Grad p_k|^2}{\pi_m\vr_k\vt}}\leq c\left(1+\intTO{\frac{Y_k^2\vC_0|\Grad \pi_m|^2}{\pi_m\vr_k\vt}}\right).
	\end{equation}
The r.h.s. of above can be, due to \eqref{decomp}, estimated as follows
	\begin{multline*}
	\intTO{\frac{Y_k^2\vC_0|\Grad \pi_m|^2}{\pi_m\vr_k\vt}}=\intTO{\frac{Y_k\vC_0|\sumkN(\Grad  p_k)|^2}{\pi_m\vr\vt}}\\
	\leq c\intTO{\frac{\vC_0}{\pi_m\vr\vt}\left({\left|\sumkN(\vC\Grad  p)_k\right|^2}+\frac{|\Grad(\vr\vt)|^2}{\left(\sum_{k=1}^{n}m_{k} Y_k\right)^2}+\frac{|\sum_{k=1}^{n}m_{k}(\vC\Grad{p})_{k}|^2}{\left(\sum_{k=1}^{n}m_{k} Y_k \right)^2}\right)},
	\end{multline*}
which is bounded thanks to \eqref{main-temp-est}, \eqref{reduced2} and \eqref{negdens}. In consequence, \eqref{almostY} is bounded. 
Recalling assumptions imposed on $C_0$ \eqref{assC} and the form of molecular pressure $\pi_m$, we deduce that
	$$\intTO{\frac{\Un{\vC_0}(1+\vt)|\Grad\vr_k|^2}{\vr_k}}\leq c\left(1+\intTO{\frac{(1+\vt)\vr_k|\Grad\vt|^2}{\vt^2}}\right)$$
and the r.h.s. is bounded, again by  \eqref{maxY} and \eqref{main-temp-est}. $\Box$

{We now present several additional estimates of $\vr$ and $\vu$ based on imbeddings of Sobolev spaces and simple interpolation inequalities. Once the B-D estimate is obtained, these estimates can be proven exactly as in the paper of Bresch and Desjardins devoted  to the Navier-Stokes-Fourier system \cite{BD}. However, for the sake of completeness, we recall them below.}

\bigskip

{\bf Further estimates of $\vr$.} From \eqref{coldp} and \eqref{negdens}  we deduce that there exist functions
$\xi_1(\vr)=\vr$ for $\vr<(1-h)$, $\xi_1(\vr)=0$ for $\vr>1$  and $\xi_2(\vr)=0$ for $\vr<1$, $\xi_2(\vr)=\vr$ for $\vr>(1+h)$ and small parameter $h>0$,
such that
	$$\|\Grad\xi_1^{-\frac{\gamma^-}{2}}\|_{L^{2}((0,T)\times\Omega)},\  \|\Grad\xi_2^{\frac{\gamma^+}{2}}\|_{L^{2}((0,T)\times\Omega)}\leq c,$$
additionally in accordance to \eqref{ess_ener} we are allowed to use the Sobolev imbeddings, thus
	\begin{equation}\label{xipos}
	\|\xi_1^{-\frac{\gamma^-}{2}}\|_{L^{2}(0,T;L^6(\Omega))},\ \|\xi_2^{\frac{\gamma^+}{2}}\|_{L^{2}(0,T;L^6(\Omega))}\leq c.
	\end{equation}

	\begin{rmk}
	Note in particular that the first of these estimate implies that
		\begin{equation*}\label{!!!}
		\vr(t,x)>0\quad \mbox{a.e.\ on\ }(0,T)\times\Omega.
		\end{equation*}
	Moreover, by interpolation between \eqref{xipos} and \eqref{ess_ener} we can check that
	\begin{equation}\label{pic}
		\|\vr e_c(\vr)\|_{L^{5/3}((0,T)\times\Omega)}\leq c.
		\end{equation}
	
	\end{rmk}
Similarly, combination of \eqref{posdens} with \eqref{ess_ener} leads to $\|\vr^{\frac{1}{2}}\|_{L^{6}(\Omega)}\leq c\left\|{\Grad\sqrt{\vr}}\right\|_{L^{2}(\Omega)},$
and therefore
	\begin{equation}\label{maxdens}
	\|\vr\|_{L^\infty(0,T;L^3(\Omega))}\leq c.
	\end{equation}

\bigskip

{\bf Estimate of $\vu$.} We use the
H${\rm\ddot{o}}$lder inequality  to write
	\begin{equation}\label{gradu}
	\|\Grad\vu\|_{L^{p}(0,T;L^{q}(\Omega))}\\
	\leq c\left(1+\|\xi_1(\vr)^{-1/2}\|_{L^{2\gamma^-}(0,T;L^{6\gamma^-}(\Omega))}\right)\|\sqrt{\vr}\Grad\vu\|_{{L^{2}((0,T)\times\Omega)}},
	\end{equation}
where $p=\frac{2\gamma^-}{\gamma^-+1}$, $q=\frac{6\gamma^-}{3\gamma^-+1}.$
Therefore, the Korn inequality together with the Sobolev imbedding imply
	\begin{equation}\label{estuI}
	\|\vu\|_{L^{\frac{2\gamma^-}{\gamma^-+1}}(0,T;L^{\frac{6\gamma^-}{\gamma^-+1}}(\Omega))}\leq c.
	\end{equation}
Next, by a similar argument 
	\begin{equation}\label{moduI}
	\|\vu\|_{L^{p'}(0,T;L^{q'}(\Omega))}
	\leq c\left(1+\|\xi_1(\vr)^{-1/2}\|_{L^{2\gamma^-}(0,T;L^{6\gamma^-}(\Omega))}\right)\|\sqrt{\vr}\vu\|_{{L^{\infty}(0,T;L^2(\Omega))}},
	\end{equation}
with $p'=2\gamma^-$, $q'=\frac{6\gamma^-}{3\gamma^-+1}$. By a simple interpolation between \eqref{estuI} and \eqref{moduI}, we obtain 
	\begin{equation}\label{optuI}
	\|\vu\|_{ L^{\frac{10\gamma^-}{3\gamma^-+3}}(0,T;L^{\frac{10\gamma^-}{3\gamma^-+3}}(\Omega))}\leq c,
	\end{equation}
and since $\gamma^->3$, we see in particular that $\|\vu\|_{ L^{{5}/{2}}(0,T;L^{{5}/{2}}(\Omega))}\leq c$ uniformly with respect to $\ep$ and $\lambda$.\\

{\bf Strict positivity of the absolute temperature.}
We now give the proof of uniform, with respect to $\ep$, positivity of $\vt$. Note that so far this was clear on account of the bound for $\ep\vt^{-3}$ in $L^1((0,T)\times\Omega)$ following from \eqref{est:entr}.
	
	\begin{lemma}
	We have, uniformly with respect to $\ep$ and $\lambda$:
		\begin{equation*}\label{???}
		\vt_\ep(t,x)> 0\quad \mbox{a.e.\ on\ }(0,T)\times\Omega.
		\end{equation*}
	\end{lemma}

\pf The above statement is a consequence of the following estimate 
	\begin{equation}\label{tempbound}
	\intTOB{|\log\vt_\ep|^{2}+|\Grad\log\vt_\ep|^{2}}\leq c,
	\end{equation}
which can be obtained by the application of generalized Korn inequality provided that we control the $L^1((0,T)\times\Omega)$ norm of $\vr|\log\vt|$.
By  \eqref{est:entr} we have
	$$\intO{(\vr_\ep s_\ep)^{0}}\leq \intO{\vr_\ep s_\ep (T)},$$
thus substituting the form of $\vr s$ from \eqref{entropy} we obtain
	$$- \intO{\vr_\ep\log{{\vt_\ep}}(T)}\leq-\sumkN\intO{ \frac{\vr_{k,\ep}}{m_{k}}\log{\frac{\vr_{k,\ep}}{m_{k}}}(T)}-\intO{(\vr_\ep s_\ep)^{0}}$$
and the r.h.s. is bounded on account of  \eqref{maxY}, \eqref{maxdens} and the initial condition. On the other hand, the positive part of the integrand $\vr_\ep\log\vt_\ep$ is bounded from above by $\vr_\ep\vt_\ep$ which belongs to $L^\infty(0,T;L^1(\Omega))$ due to \eqref{main_ener}, so we  end up with
	\begin{equation*}\label{sinf}
	{\rm{ess}}\sup_{t\in(0,T)}\intO{\left|\vr_\ep\log\vt_\ep(t)\right|}\leq c,
	\end{equation*}
which completes the proof of \eqref{tempbound}. $\Box$

\subsection{Passage to the limit with $\ep\to 0$ }
With the B-D estimate at hand, especially with the bound on $\lap^{s+1}\vr_\ep$ in $L^2((0,T)\times\Omega)$, which is now uniform with respect to $\ep$, we may perform the limit passage similarly as in previous step. Indeed, the uniform estimates allow us to extract subsequences, such that
	$$\ep\lap^s\Grad\vu_{\ep},~\ep\Grad\vr_{\ep},~\ep\lap^{s+1}\vr_\ep\rightarrow 0\ {\rm strongly\ in\ } L^{2}((0,T)\times\Omega),$$
therefore
	$$\ep\Grad\vr_{\ep}\Grad\vu_{\ep}\rightarrow 0\ {\rm strongly\ in\ } L^{1}((0,T)\times\Omega).$$
The strong convergence of the density as well as the velocity (since $\vr_\ep>c(\lambda)$) can be obtained identically as  in the previous step. Therefore we focus only on the strong convergence of the temperature and the limit passage in the total energy balance.

Recall that from \eqref{est:entr}, \eqref{main-temp-est} it follows that
	$$\vt_\ep\to\vt\quad\mbox{weakly\ in\ }L^2(0,T;W^{1,2}(\Omega))$$
and
	$$\ep\vt_\ep^{-2},\ep\vt_\ep^5\to 0\quad\mbox{strongly\ in\ }L^1((0,T)\times\Omega).$$
The pointwise convergence of the temperature is to be deduced from the version of the Aubin-Lions lemma, see \cite{EF1}:

	\begin{lemma}\label{LLions}
	Let $v_{\ep}$ be sequence of functions bounded in $L^{2}(0,T;L^{q}(\Omega))$ and  in  $L^{\infty}(0,T;L^{1}(\Omega))$, where $q>\frac{6}{5}$. Furthermore, assume that
		\begin{equation}\label{cond0}
		\pt v_\ep\geq g_\ep\quad in\  {\cal D}'((0,T)\times\Omega),
		\end{equation}
	where
		\begin{equation}\label{condI}
		g_{\ep}\ is\ bounded\ in\ L^{1}(0,T;W^{-m,r}(\Omega))\ for\ some\ m\geq0,\ r>1\ independently \ of\ \ep.
		\end{equation}
	Then there exists a subsequence $v_{\ep}$ which converges to $v$ strongly in $L^2(0,T;W^{-1,2}(\Omega))$.
	\end{lemma}

We will apply this lemma to $v_\ep=\vr_\ep \vt{\ep}+{\beta}\vt_\ep^4$.
Then, on account of \eqref{vt_in_delta} we can repeat the estimates following \eqref{solo_vt} that led to \eqref{ptvtN} to check that
	\begin{multline*}
	\pt v_\ep\geq g_\ep=-\Div(\vu_\ep(\vr_\ep \vt_{\ep}+\beta\vt^4_\ep))+\Div\left(\kappa_{\ep}(\vr,\vt)\Grad\vt_\ep\right)
	-\sumkN \Div \left(\vt_\ep\frac{\vf_{k,\ep}}{m_k}-\ep\vt_\ep\Grad\erki{\ep}\right)\\
	+\frac{\ep}{\vt^2}-\ep\vt^{5}
	+\ep{\frac{1}{\vr_\ep}\frac{\partial\pi_c(\vr_\ep)}{\partial\vr_\ep}|\Grad\vr_\ep|^2}-\lr{\pi_{m,\ep}+\frac{\beta}{3}\vt^4_\ep}\Div\vu_\ep+2\vr_\ep|\D(\vu_\ep)|^2
	+\lambda|\lap^{s}\Grad(\vr_\ep\vu_\ep)|^{2}+\lambda\ep|\lap^{s+1}\vr_\ep|^{2}.
	\end{multline*}
Moreover, the r.h.s. is bounded in $L^{1}(0,T;W^{-1,p}(\Omega))\cup L^{1}(0,T;W^{-2,q}(\Omega))$ for some $p,q>1$. 
Therefore, the above lemma and the strong convergence  of $\vr_\ep$ imply in particular that
	$$\vt_\ep^4\to\Ov{\vt^4}\quad \mbox{strongly\ in\ }L^2(0,T;W^{-1,2}(\Omega)).$$
On the other hand, we know also that $\vt_\ep\to\vt$ weakly in $L^2(0,T;W^{1,2}(\Omega))$, therefore a simple argument based on the monotonicity of $f(x)=x^4$ implies strong convergence of $\vt_\ep$ in $L^q(0,T; L^3q(\Omega))$ for any $q<B$.\\

Let us finish this subsection with the list of the limit equations:
\begin{itemize}
\item[-]
the continuity equation
	\begin{equation*}\label{apcontL}
	\pt\vr+\Div (\vr \vu)= 0
	\end{equation*}
is satisfied pointwisely on $[0,T]\times\Omega$;

\item[-]
the momentum equation
	\begin{equation}\label{momL}
	\begin{split}
	-&\intTO{\vr\vu\cdot\pt\vcg{\phi}}-\intO{\vc{m}^{0}\cdot\vcg{\phi}(0)}
	+\intTO{\lambda\lap^{s}\Grad(\vr\vu) :\lap^{s}\Grad (\vr\vcg{\phi})}\\
	-&\intTO{(\vr{\vu}\otimes{\vu}):\Grad \vcg{\phi}}
	+\intTO{2\vr\D\vu:\D \vcg{\phi}}
	-\intTO{\pi\Div \vcg{\phi}}\\
	+&\lambda\intTO{\lap^s\Div\left({\vr\vcg{\phi}}\right)\lap^{s+1}\vr}=0
	\end{split}
	\end{equation}
holds for any test function $\vcg{\phi}\in L^2(0,T;W^{2s+1}(\Omega))\cap W^{1,2}(0,T;W^{1,2}(\Omega))$ such that $\vcg{\phi}(\cdot,T)=\vc{0}$.

\item[-]
the species equations
	\begin{multline*}\label{rk_del}
 	 \intTO{\vr_k\pt\phi}+\intO{\vr_k^0\phi(0)}\\+\intTO{ \vr_{k} \vu\cdot\Grad\phi}+\intTO{\vf_k\cdot\Grad\phi}  =  -\intTO{\vr\vt\vw_{k}\phi},\quad k\in\{1,...,n\}
	\end{multline*}
are fulfilled for any smooth function $\phi$ such that $\phi(\cdot,T)=0$;

\item[-]
the total energy equation
	\begin{equation}\label{totL}
	\begin{split}
	&\intTO{\lr{\vr e+\frac{1}{2}\vr|\vu|^2+\frac{\lambda}{2}|\Grad^{2s+1}\vr|^2}\pt\phi} +\intO{\lr{\vr e+\frac{1}{2}\vr|\vu|^2+\frac{\lambda}{2}|\Grad^{2s+1}\vr|^2}(0)\phi(0)}\\
	&+\intTO{\lr{\vr e\vu+\frac{1}{2}\vr|\vu|^2\vu}\cdot\Grad\phi}-\intTO{\kappa\Grad\vt\cdot\Grad\phi}+\sumkN \intTO{ \vt\frac{\vf_{k}}{m_k}\cdot\Grad\phi}\\
	&+\intTO{\pi\vu\cdot\Grad\phi}-\intTO{(2\vr\D(\vu)\vu)\cdot\D\phi}= \intTO{R_\lambda(\vr,\vt,\vu,\phi)},
	\end{split}
	\end{equation}
holds for any smooth function $\phi$ such that $\phi(\cdot,T)=0$;
where the heat flux term is to be understood as in \eqref{weak_kappa1}, and
\begin{equation}\label{RL} 
\begin{split}
{R_{\lambda}(\vr,\vt,\vu,\phi)} = &\lambda\left[\Delta^s ( \Div (\vr \vu \phi)) \Delta^{s+1} \vr 
-\lap^s\Div(\vr\vu)\lap^{s+1}\vr\phi\right]\\
&-\lambda\lap^s\Div(\vr\vu)\Grad\lap^s\vr\cdot\Grad\phi
- \lambda\left[|\Delta^s (\nabla (\vr \vu))|^2 \phi- \Delta^s \nabla (\vr \vu): \Delta^s \nabla (\vr \vu \phi)\right] .
\end{split}
\end{equation}
\end{itemize}
Moreover, using the lower weak semicontinuity of norm and passing to the limit in \eqref{vt_in_delta}	 
\begin{equation}\label{intL}
\begin{split}
		\ptb{\vr\vt+\beta\vt^4} +\Div(\vu(\vr\vt+\beta\vt^4))-\Div\left(\kappa\Grad\vt\right)
		+\sumkN \Div \left(\vt\frac{\vf_{k}}{m_k}\right)&\geq \\
	-\lr{\pi_m +\frac{\beta}{3}\vt^4}\Div\vu &+2\vr|\D(\vu)|^2
		+\lambda|\lap^{s}\Grad(\vr\vu)|^{2},
		\end{split}
		\end{equation}
satisfied in the sense of distributions on $(0,T)\times\Omega$.

\subsection{Passage with $\lambda\to 0$ }
In this section we present the argument for the convergence of a sequence $(\vr_\lambda,\vu_\lambda,\vt_\lambda, \vr_{1,\lambda},\ldots,\vr_{n,\lambda})$ to a solution $(\vr,\vu,\vt,\vr_1,\ldots,\vr_n)$ as specified at the beginning. Some of the arguments here are repetitions from our previous works \cite{EZ2,EZ3,MPZ1}, and so we only recall their formulations.\\

{\bf Strong convergence of the density.} 
The strong convergence of a sequence $\vr_\lambda$ is guaranteed by the following lemma
\begin{lemma}\label{Lemma:vr}
There exists a subsequence $\vr_\lambda$ such that
\begin{equation*}\label{pointvr}
\sqrt{\vr_\lambda}\rightarrow\sqrt{\vr}\quad a.e.\ and\  strongly\ in\ L^2((0,T)\times\Omega).
\end{equation*}
Moreover $\vr_\lambda\rightarrow\vr$ strongly in $C([0,T];L^p(\Omega))$, $p<3$.
\end{lemma}
For the proof see \cite{EZ2}, Lemma 7.\\

{\bf Strong convergence of the species densities.}
Analogously we show the strong convergence of species densities. We have
%
\begin{lemma}\label{Lemma_vrk}
Up to a subsequence the partial densities $\vr_{k,\lambda}$, $k=1,\ldots,n$ converge strongly in $L^p(0,T;L^{q}(\Omega))$, $1\leq p<\infty, 1\leq q<3$ to $\vr_k$. In particular
\begin{equation*}\label{pointvrk}
\vr_{k,\lambda}\to\vr_k\quad a.e.\ in\  (0,T)\times \Omega.
\end{equation*}
Moreover $\vr_{k,\lambda}\to\vr_k$  in $C([0,T];L^3_{\rm{weak}}(\Omega))$.
\end{lemma}
For the proof see \cite{EZ3}, Lemma 17.\\

{\bf Convergence of the convective term.} {To prove this, we will again take advantage of the special form of the cold component of the pressure close to vacuum \eqref{coldp}. First we show under what condition on $\gamma^-$,  the convective term $\vr_\lambda|\vu_\lambda|^3$ is uniformly bounded in $L^p((0,T)\times\Omega)$ for some $p>1$. Following the proof from \cite{BD}, we may write}
	$$\vr_\lambda^{1/3}|\vu_\lambda|=\vr_\lambda^{1/3-\alpha}\vr^\alpha|\vu_\lambda|^{2\alpha}|\vu_\lambda|^{1-2\alpha}$$
and we will use the interpolation inequality for $\vr_\lambda\in L^{\infty}(0,T;L^{\gamma^+}(\Omega))$, $\vr_\lambda|\vu_\lambda|^2\in L^\infty(0,T;L^1(\Omega))$ and $\vu_\lambda\in L^{\frac{2\gamma^-}{\gamma^-+1}}(0,T;L^{\frac{6\gamma^-}{\gamma^-+1}}(\Omega))$. So, $\vr_\lambda^{1/3}\vu_\lambda\in L^p(0,T;L^q(\Omega))$ with $p,q>3$ if
	$$\lr{1-2\alpha}\frac{\gamma^-+1}{2\gamma^-}<\frac{1}{3},$$
	$$\lr{{1\over3}-\alpha}\frac{1}{\gamma^+}+\alpha+\lr{1-2\alpha}\frac{\gamma^-+1}{6\gamma^-}<\frac{1}{3},$$
meaning that $\gamma^-$ and $\gamma^+$ must satisfy the following relation
$$\gamma^->\frac{5\gamma^+-3}{\gamma^+-3}.$$
\begin{rmk}
Improvement of this condition would require some faster growth at infinity of the viscosity coefficiet $\mu(\vr)$ which is equal $\vr$ in our case. Indeed, then the above interpolation procedure could use only the additional bound of $\vr$ following from the B-D estimate, without involving the bound in $L^\infty(0,T; L^{\gamma^+}(\Omega))$, similarly as in \cite{BD}. However, in the case of chemically reacting mixtures, modification of the viscosity coefficient would lead to a problem with closing the B-D estimate, see \eqref{pepi0}.
\end{rmk}
The above estimate implies that provided $\gamma^-$, $\gamma^+$ fulfill conditions specified above, the convective term $\vr_\lambda\vu_\lambda^3$ converges weakly to $\Ov{\vr\vu^3}$ in $L^r((0,T)\times\Omega)$ for some $r>1$. To identify the limit, we prove the following lemma.
\begin{lemma}\label{lemma:vu}
We have, up to a choice of subsequence
\begin{eqnarray}
\vr_\lambda\vu_\lambda&\rightarrow& \vr\vu\quad   in\   C([0,T];L^\frac{3}{2}_{\rm{weak}}(\Omega))\label{Cpweak},\\
\vr_\lambda^{1/3}\vu_\lambda&\rightarrow& \vr^{1/3}\vu\quad  strongly\  in\   L^{p}((0,T)\times\Omega),\ for \ some \ p>3,\label{convconv}\\
\vr_\lambda\D\vu_\lambda\vu_\lambda&\rightarrow& \vr\D\vu\vu\quad  weakly\  in\   L^{1}((0,T)\times\Omega)\label{convD}.
\end{eqnarray}
\end{lemma}
\pf  We already know that $\vr_\lambda$ converges  to $\vr$ a.e. on $(0,T)\times\Omega$. Moreover, due to \eqref{estuI}, up to extracting a subsequence, $\vu_\lambda$ converges weakly to $\vu$ in  $L^p(0,T;L^q(\Omega))$ for some $p>\frac{3}{2},\ q>\frac{9}{2}$. Since $\sqrt{\vr_\lambda}\vu_\lambda$ is uniformly bounded in $L^\infty(0,T;L^2(\Omega))$ and $\sqrt{\vr_\lambda}$ is uniformly bounded in $L^\infty(0,T;L^6(\Omega))$,  the sequence $\vr_\lambda\vu_\lambda$ is uniformly bounded in  $L^{\infty}(0,T; L^{\frac{3}{2}}(\Omega))$. All together it implies that 
$$\vr_\lambda\vu_\lambda \rightarrow \vr\vu\quad \mbox{ weakly}^*\mbox{\  in}\   L^{\infty}(0,T;L^{\frac{3}{2}}(\Omega)).$$
Now, we are aimed at improving the time compactness of this sequence. Using the differential form of \eqref{momL}
and the uniform estimates \eqref{negdens}, \eqref{posdens}, \eqref{gradcomp}, \eqref{xipos}, \eqref{gradu} we verify that the sequence of functions
$t\mapsto\intO{\vr_\lambda\vu_\lambda\cdot\vcg{\phi} }$
is uniformly bounded and equicontinuous in $C([0,T])$, where $\vcg{\phi}\in C^\infty_c(\Omega)$. But since the smooth functions are dense in $L^3(\Omega)$, applying the Arzel\`a-Ascoli theorem, we show \eqref{Cpweak}.

On the other hand, due to \eqref{gradu}, $\vu_\lambda$ is uniformly bounded in $L^p(0,T;W^{1,q}(\Omega))$ for some $p>\frac{3}{2},\ q>\frac{9}{5}$, so it converges to $\vu$ weakly in this space. Since $W^{1,q}(\Omega)$ with $q>\frac{9}{5}$ is compactly embedded into $L^3(\Omega)$, by \eqref{Cpweak}, we deduce that $\vr_\lambda|\vu_\lambda|^2$ converges to $\vr|\vu|^2$ strongly in $L^1((0,T)\times\Omega)$. 
We can use this fact and the strong convergence of the density in order to show that $\vr_\lambda^{1/3}\vu_\lambda=\vr_\lambda^{1/3}\vu_{\lambda}\vc{1}_{\{\vr_\lambda<\eta\}}+\vr^{1/2}_\lambda\vu_\lambda\vr^{-1/6}\vc{1}_{\{\vr_\lambda>\eta\}}$ converges to $\vr^{1/3}\vu$ strongly in $L^1((0,T)\times\Omega)$. This in turn, when combined with the uniform bound on $\vr_\lambda^{1/3}\vu_\lambda$ in $L^p((0,T)\times\Omega)$, $p>3$, yields \eqref{convconv}. 
Finally, rewriting $\vr_\lambda(\D\vu_\lambda)\vu_\lambda=\sqrt{\vr_\lambda}\D\vu_\lambda\vr_\lambda^{1/3}\vu_\lambda\vr_\lambda^{1/6}$ and using \eqref{convconv}, strong convergence of the density and a weak convergence of $\sqrt{\vr_\lambda}(\D\vu_\lambda)$ to $\sqrt{\vr}(\D\vu)$ in $L^2((0,T)\times\Omega)$ we prove \eqref{convD}.
$\Box$\\

{\bf Strong convergence of the temperature.} The difference with respect to previous chapter is that we cannot use the higher order estimates either for the velocity or for the density in order to deduce the boundedness of the time derivative of temperature in an appropriate space. However, the idea of proving compactness of the temperature is, as previously, to apply Lemma \ref{LLions} with $\ep=\lambda$, $v_\lambda=\vr_\lambda \vt_\lambda+\beta\vt_\lambda^4$. Therefore, our next aim is to check that its assumptions are satisfied uniformly with respect to $\lambda$.

First, let us note that $v_\lambda$ is bounded in $L^{2}(0,T;L^{q}(\Omega))$ and  in  $L^{\infty}(0,T;L^{1}(\Omega))$, where $q>\frac{6}{5}$, uniformly with respect to $\lambda$. Indeed, it follows directly from \eqref{ess_ener} and \eqref{maxtemp}. Further, from \eqref{intL} one deduces  that $\pt v_\lambda\geq g_\lambda$, where $g_\lambda$  has the following form
\begin{equation}\label{pgL}
\begin{split}
g_\lambda=&-\Div\lr{\vu_\lambda\lr{\vr_\lambda \vt_\lambda+\beta\vt_\lambda^4}}+\Div\left(\kappa(\vr_\lambda,\vt_\lambda)\Grad\vt_\lambda\right)
		-\sumkN \Div \left(\vt_\lambda\frac{\vf_{k,\lambda}}{m_k}\right)\\
	&-\lr{\pi_{m,\lambda}+\frac{\beta}{3}\vt_\lambda^4}\Div\vu_\lambda+2\vr_\lambda|\D(\vu_\lambda)|^2
		+\lambda|\lap^{s}\Grad(\vr_\lambda\vu_\lambda)|^{2}
\end{split}
\end{equation}
and is bounded\ in\ $L^{1}(0,T;W^{-m,r}(\Omega))$ for some $m\geq0,\ r>1$ independently of $\lambda$. Indeed,  this can be estimated, similarly to (\ref{solo_vt}--\ref{ptvtN}) except for the terms that contains velocity. For them we may write
	$$\|\vu_\lambda\vr_\lambda\vt_\lambda\|_{L^{{12}/{11}}((0,T)\times\Omega)}\leq c\|\sqrt{\vr_\lambda}\vu_\lambda\|_{L^\infty(0,T;L^2(\Omega))}\|\sqrt{\vr_\lambda}\|_{L^\infty(0,T;L^6(\Omega))}\|\vt_\lambda\|_{L^\infty(0,T;L^4(\Omega))}\leq c,$$
on account of \eqref{eqT} and \eqref{optuI}, further 
	$$\|\vu_\lambda\vt_\lambda^4\|_{L^{40/31}((0,T)\times\Omega)}\leq\|\vu_\lambda\|_{L^{5/2}((0,T)\times\Omega)}\|\vt_\lambda^4\|_{L^{8/3}((0,T)\times\Omega)}\leq c.$$
For the internal pressure we have
$$\|\vr_{k,_\lambda}\vt_\lambda\Div\vu_\lambda\|_{L^{{12}/{11}}((0,T)\times\Omega)}\leq c\|\sqrt{\vr_\lambda}\Div\vu_\lambda\|_{L^2((0,T)\times\Omega)}\|\sqrt{\vr_\lambda}\|_{L^\infty(0,T;L^6(\Omega))}\|\vt_\lambda\|_{L^\infty(0,T;L^4(\Omega))}\leq c$$
and the other term $\vt_\lambda^4\Div\vu_\lambda$ is bounded in $L^{1}((0,T)\times\Omega)$ as  was shown above in \eqref{vt4vu}.
Since the two last terms in \eqref{pgL} are also uniformly bounded in $L^1((0,T)\times\Omega)$, the assumptions of Lemma \ref{LLions} are satisfied with $m=1$, $r>1$. Therefore, there exists a subsequence, $v_\lambda$ converging to $v$ strongly in $L^2(0,T; W^{-1,2}(\Omega))$, which can be used to show the strong convergence of $\vt_\lambda$ exactly as in the previous section.\\

{\bf Passage to the limit in the nonlinear terms.} The last step in the limit passage $\lambda\to0$ is verification of convergence in the nonlinear terms of the system. The most demanding of them are in the energy equation \eqref{totL} and \eqref{RL}, and we will justify the limit passage only there. The correction of energy ${\lambda}\Grad^{2s+1}\vr_\lambda\to0$ strongly in $L^2((0,T)\times\Omega)$, therefore the energy $E_{\lambda}=\vr_\lambda e_c(\vr_\lambda)+\vr_\lambda\vt_\lambda+\beta\vt_\lambda+\frac{1}{2}\vr_\lambda|\vu_\lambda|^2+\frac{\lambda}{2}|\Grad^{2s+1}\vr_\lambda|^2$ converges to $E$ due to strong convergence of $\vr_\lambda$, $\vt_\lambda$ and Lemma \ref{lemma:vu}. Similarly $\vu_\lambda\vr_\lambda e_{\lambda}$, $\vr_\lambda\vu_\lambda^3$ and $\pi(\vr_\lambda,\vt_\lambda,\vr_{k,\lambda})\vu_\lambda$ converge weakly to $\vu\vr e$, $\vr\vu^3$ and $\pi\vu$ respectively, due to uniform bounds in $L^p((0,T)\times\Omega)$ for $p>1$ from above, estimate \eqref{pic}, the strong convergence of $\vr_\lambda$, $\vt_\lambda$, $\vr_{k,\lambda}$ and Lemma \ref{lemma:vu}.

Limit passage in the heat flux term $\kappa(\vr_\lambda,\vt_\lambda)\Grad\vt_\lambda$ can be performed within formulation analogous to \eqref{weak_kappa1}, since it involves only the sequences $\vr_\lambda$, $\vt_\lambda$ which are strongly convergent and a sequence $\Grad\vt_\lambda$ which converges to $\Grad\vt$ weakly in $L^2((0,T)\times\Omega)$. Similarly, we note that $\vt_\lambda\vf_k(\vr_\lambda,\vt_\lambda,\vr_{k,\lambda})$ is affine with respect to weakly convergent sequences $\Grad\vr_{k,\lambda}$, $\Grad\vt_\lambda$ with coefficients that are strongly convergent as $\lambda\to 0$. Thus, using the uniform bounds from above we verify that up to a subsequence $\vt_\lambda\vf_k(\vr_\lambda,\vt_\lambda,\vr_{k,\lambda})\to \vt\vf_k\lr{\vr,\vt,\vr_k}$ weakly in $L^p((0,T)\times\Omega)$ for some $p>1$. Passage to the limit in the last term from the r.h.s. of \eqref{totL} was proven in Lemma \ref{lemma:vu}.

We are now ready to prove that the corrector term $R_\lambda$ converges to 0 strongly in $L^1((0,T)\times\Omega)$ as $\lambda\to0$, or rather that the most demanding terms listed in \eqref{RL} vanish when $\lambda\to0$. 

First of all observe that due to \eqref{negdens} \eqref{lambda_vu} we have
\begin{equation*}
\begin{split}
&\intTO{\left|R_\lambda(\vr_\lambda,\vu_\lambda,\phi)\right|}\\
=&\lambda\intTO{\left|\lr{
\lap^{s}\Grad(\vr_\lambda\vu_\lambda)\lap^{s}(\vr_\lambda\vu_\lambda)
+\lap^s(\vr_\lambda\vu_\lambda)\lap^{s+1}\vr_\lambda+\lap^s\Div(\vr_\lambda\vu_\lambda)\Grad\lap^s\vr_\lambda}\cdot\Grad\phi\right|}\\
\leq &
c\lr{\|\Grad\phi\|_{L^\infty((0,T)\times\Omega)}}\lambda\Big[\|\vr_\lambda\vu_\lambda\|_{L^2(0,T;W^{2s+1,2}(\Omega))}\|\vr_\lambda\vu_\lambda\|_{L^2(0,T;W^{2s,2}(\Omega))}\\
&+\|\vr_\lambda\vu_\lambda\|_{L^2(0,T;W^{2s,2}(\Omega))}\|\vr_\lambda\|_{L^2(0,T;W^{2s+2,2}(\Omega))}+\|\vr_\lambda\vu_\lambda\|_{L^2(0,T;W^{2s+1,2}(\Omega))}\|\vr_\lambda\|_{L^2(0,T;W^{2s,2}(\Omega))}\Big]\\
\leq &c\lr{\|\Grad\phi\|_{L^\infty((0,T)\times\Omega)}}\sqrt{\lambda}\Big[\|\vr_\lambda\vu_\lambda\|_{L^2(0,T;W^{2s,2}(\Omega))}+\|\vr_\lambda\|_{L^2(0,T;W^{2s,2}(\Omega))}\Big],
\end{split}
	\end{equation*}
thus the task is to show that the term after the last inequality symbol converges to $0$. But this is evident, since one can use the Gagliardo-Nirenberg interpolation inequality and uniform bounds for $\vr_\lambda\vu_\lambda$ in $L^\infty(0,T; L^{3/2}(\Omega))$ and for $\vr_\lambda$ in $L^\infty(0,T;L^3(\Omega))$ together with uniform bounds for $\sqrt\lambda\vr_\lambda\vu_\lambda$  and for $\sqrt{\lambda }\vr_\lambda$ in $L^2(0,T;W^{2s+1,2}(\Omega))$ and $L^2(0,T;W^{2s+2,2}(\Omega))$, respectively. This finished the proof of the main theorem. $\Box$

\medskip

\noindent{\footnotesize {\bf Acknowlegdement.} The work of P.B.M. and E.Z. has been supported by Polish MN grant IdP2011 000661 and by the International Ph.D.
Projects Programme of Foundation for Polish Science operated within the
Innovative Economy Operational Programme 2007-2013 funded by European
Regional Development Fund (Ph.D. Programme: Mathematical
Methods in Natural Sciences)}.

%

\footnotesize

\begin{thebibliography}{10}

\bibitem{BE}
J.~Bebernes and D.~Eberly.
\newblock {\em Mathematical problems from combustion theory}, volume~83 of {\em
  Applied Mathematical Sciences}.
\newblock Springer-Verlag, New York, 1989.

\bibitem{B2010}
D.~Bothe.
\newblock On the {M}axwell-{S}tefan approach to multicomponent diffusion.
\newblock In {\em Parabolic problems}, volume~80 of {\em Progr. Nonlinear
  Differential Equations Appl.}, pages 81--93. Birkh{\"a}user/Springer Basel
  AG, Basel, 2011.

\bibitem{BGS12}
L.~Boudin, B.~Grec, and F.~Salvarani.
\newblock A mathematical and numerical analysis of the maxwell-stefan diffusion
  equations.
\newblock {\em Discrete Contin. Dyn. Syst. Ser. B}, 17(5):1427--1440, 2012.

\bibitem{BD2}
D.~Bresch and B.~Desjardins.
\newblock Existence of global weak solutions for a 2{D} viscous shallow water
  equations and convergence to the quasi-geostrophic model.
\newblock {\em Comm. Math. Phys.}, 238(1-2):211--223, 2003.

\bibitem{BDC} 
\newblock D. Bresch and B. Desjardins,
\newblock {On the construction of approximate solutions for the 2{D} viscous shallow water model and for compressible {N}avier-{S}tokes models},
\newblock \emph{J. Math. Pures Appl. (9)}, {86} (2006), 362--368.

\bibitem{BD}
D.~Bresch and B.~Desjardins.
\newblock On the existence of global weak solutions to the {N}avier-{S}tokes
  equations for viscous compressible and heat conducting fluids.
\newblock {\em J. Math. Pures Appl. (9)}, 87(1):57--90, 2007.

\bibitem{BDG}
D.~Bresch, B.~Desjardins, and D.~G{{\'e}}rard-Varet.
\newblock On compressible {N}avier-{S}tokes equations with density dependent
  viscosities in bounded domains.
\newblock {\em J. Math. Pures Appl. (9)}, 87(2):227--235, 2007.

\bibitem{BD3}
D.~Bresch, B.~Desjardins, and C.-K. Lin.
\newblock On some compressible fluid models: {K}orteweg, lubrication, and
  shallow water systems.
\newblock {\em Comm. Partial Differential Equations}, 28(3-4):843--868, 2003.

\bibitem{CHT}
G.-Q. Chen, D.~Hoff, and K.~Trivisa.
\newblock Global solutions to a model for exothermically reacting, compressible
  flows with large discontinuous initial data.
\newblock {\em Arch. Ration. Mech. Anal.}, 166(4):321--358, 2003.

\bibitem{CJ13}
X.~Chen and A.~J{\"u}ngel.
\newblock Analysis of an incompressible {N}avier-{S}tokes-{M}axwell-{S}tefan
  system.
\newblock {\em arXiv:1310.3376}, 10 2013.

\bibitem{DT2}
D.~Donatelli and K.~Trivisa.
\newblock On the motion of a viscous compressible radiative-reacting gas.
\newblock {\em Comm. Math. Phys.}, 265(2):463--491, 2006.

\bibitem{DT1}
D.~Donatelli and K.~Trivisa.
\newblock A multidimensional model for the combustion of compressible fluids.
\newblock {\em Arch. Ration. Mech. Anal.}, 185(3):379--408, 2007.

\bibitem{EF1}
E.~Feireisl.
\newblock {\em Dynamics of viscous compressible fluids}, volume~26 of {\em
  Oxford Lecture Series in Mathematics and its Applications}.
\newblock Oxford University Press, Oxford, 2004.

\bibitem{FN}
E.~Feireisl and A.~Novotn{{\'y}}.
\newblock {\em Singular limits in thermodynamics of viscous fluids}.
\newblock Advances in Mathematical Fluid Mechanics. Birkh{\"a}user Verlag,
  Basel, 2009.

\bibitem{FTP}
E.~Feireisl, H.~Petzeltov{{\'a}}, and K.~Trivisa.
\newblock Multicomponent reactive flows: global-in-time existence for large
  data.
\newblock {\em Commun. Pure Appl. Anal.}, 7(5):1017--1047, 2008.

\bibitem{FGM05a}
J.~Frehse, S.~Goj, and J.~M{{\'a}}lek.
\newblock On a {S}tokes-like system for mixtures of fluids.
\newblock {\em SIAM J. Math. Anal.}, 36(4):1259--1281 (electronic), 2005.

\bibitem{FGM05b}
J.~Frehse, S.~Goj, and J.~M{{\'a}}lek.
\newblock A uniqueness result for a model for mixtures in the absence of
  external forces and interaction momentum.
\newblock {\em Appl. Math.}, 50(6):527--541, 2005.

\bibitem{Gio90}
V.~Giovangigli.
\newblock Mass conservation and singular multicomonent diffusion algorithms.
\newblock {\em IMPACT Comput. Sci. Eng.}, 2(2):73--97, 1990.

\bibitem{Gio91}
V.~Giovangigli.
\newblock Convergent iterative methods for multicomponent diffusion.
\newblock {\em Impact Comput. Sci. Engrg.}, 3(3):244--276, 1991.

\bibitem{VG}
V.~Giovangigli.
\newblock {\em Multicomponent flow modeling}.
\newblock Modeling and Simulation in Science, Engineering and Technology.
  Birkh{\"a}user Boston Inc., Boston, MA, 1999.

\bibitem{HMPW13}
M.~Herberg, M.~Meyries, J.~Pr{\"u}ss, and M.~Wilke.
\newblock Reaction-diffusion systems of maxwell-stefan type with reversible
  mass-action kinetics.
\newblock {\em arXiv:1310.4723}, 2013.

\bibitem{HS01}
D.~Hoff and J.~Smoller.
\newblock Non--formation of vacuum states for compressible {N}avier-{S}tokes
  equations.
\newblock {\em Comm. Math. Phys.}, 2016(2):255--276, 2001.

\bibitem{JS13}
A.~J{{\"u}}ngel and I.V. Stelzer.
\newblock Existence analysis of {M}axwell-{S}tefan systems for multicomponent
  mixtures.
\newblock {\em SIAM J. Math. Anal.}, 45(4):2421--2440, 2013.

\bibitem{KMP12}
N.~A. Kucher, A.~E. Mamontov, and D.~A. Prokudin.
\newblock Stationary solutions to the equations of the dynamics of mixtures of
  viscous compressible heat-conducting fluids.
\newblock {\em Sibirsk. Mat. Zh. Sib. Math. J.}, 53(6):1075--1088, 2012.

\bibitem{KT}
Y.-S. Kwon and K.~Trivisa.
\newblock Stability and large time behaviour for multicomponent reactive flows.
\newblock {\em Nonlinearity}, 22(10):2443--2471, 2009.

\bibitem{LSU}
O.~A. Lady{z}henskaya, V.~A. Solonnikov, and N.~N. Uralceva.
\newblock {\em Linear and quasilinear equations of parabolic type}.
\newblock Translated from the Russian by S. Smith. Translations of Mathematical
  Monographs, Vol. 23. American Mathematical Society, Providence, R.I., 1967.

\bibitem{LM}
M.~Lewicka and P.~B. Mucha.
\newblock On temporal asymptotics for the {$p$}th power viscous reactive gas.
\newblock {\em Nonlinear Anal.}, 57(7-8):951--969, 2004.

\bibitem{PLL}
P.-L. Lions.
\newblock {\em Mathematical topics in fluid mechanics. {V}ol. 2}, volume~10 of
  {\em Oxford Lecture Series in Mathematics and its Applications}.
\newblock The Clarendon Press Oxford University Press, New York, 1998.
\newblock Compressible models, Oxford Science Publications.

\bibitem{MT13}
M.~Marion and R.~Temam.
\newblock Global existence for fully nonlinear reaction-diffusion systems
  describing multicomponent reactive flows.
\newblock {\em arXiv:1310.2624}, 10 2013.

\bibitem{MV07}
A.~Mellet and A.~Vasseur.
\newblock On the barotropic compressible {N}avier-{S}tokes equations.
\newblock {\em Comm. Partial Differential Equations}, 32(1-3):431--452, 2007.

\bibitem{MV08}
A.~Mellet and A.~Vasseur.
\newblock Existence and uniqueness of global strong solutions for
  one-dimensional compressible {N}avier-{S}tokes equations.
\newblock {\em SIAM J. Math. Anal.}, 39(4):1344--1365, 2007/08.

\bibitem{MP3}
P.~B. Mucha and M.~Pokorn{{\'y}}.
\newblock Weak solutions to equations of steady compressible heat conducting
  fluids.
\newblock {\em Math. Models Methods Appl. Sci.}, 20(5):785--813, 2010.

\bibitem{MPZ}
P.~B. Mucha, M.~Pokorn{{\'y}}, and E.~Zatorska.
\newblock Approximate solutions to model of two-component reactive flow.
\newblock {\em Discrete Contin. Dyn. Syst. Ser. S}, (in press), 2013.

\bibitem{MPZ1}
P.~B. Mucha, M.~Pokorn{{\'y}}, and E.~Zatorska.
\newblock Chemically reacting mixtures in terms of degenerated parabolic
  setting.
\newblock {\em J. Math. Phys.}, 54(071501), 2013.

\bibitem{M1D}
P.B. Mucha.
\newblock Compressible {N}avier-{S}tokes system in 1-d.
\newblock {\em Math. Methods Appl. Sci.}, 24(9):607--622, 2001.

\bibitem{NP1}
A.~Novotn{{\'y}} and M.~Pokorn{{\'y}}.
\newblock Steady compressible {N}avier-{S}tokes-{F}ourier system for monoatomic
  gas and its generalizations.
\newblock {\em J. Differential Equations}, 251(2):270--315, 2011.

\bibitem{NP2}
A.~Novotn{{\'y}} and M.~Pokorn{{\'y}}.
\newblock Weak and variational solutions to steady equations for compressible
  heat conducting fluids.
\newblock {\em SIAM J. Math. Anal.}, 43(3):1158--1188, 2011.

\bibitem{NS}
A.~Novotn{{\'y}} and I.~Stra{\v{s}}kraba.
\newblock {\em Introduction to the mathematical theory of compressible flow},
  volume~27 of {\em Oxford Lecture Series in Mathematics and its Applications}.
\newblock Oxford University Press, Oxford, 2004.

\bibitem{RT}
K.~R. Rajagopal and L.~Tao.
\newblock {\em Mechanics of mixtures}, volume~35 of {\em Series on Advances in
  Mathematics for Applied Sciences}.
\newblock World Scientific Publishing Co. Inc., River Edge, NJ, 1995.

\bibitem{T2000}
S.~R. Turns.
\newblock An introduction to combustion: Concepts and applications.
\newblock {\em McGraw-Hill Series in Mechanical Engineering}, 2000.

\bibitem{Wal62}
L.~Waldmann and E.~Tr{{\"u}}benbacher.
\newblock Formale kinetische {T}heorie von {G}asgemischen aus anregbaren
  {M}olek{\"u}len.
\newblock {\em Z. Naturforsch}, 17a:363--376, 1962.

\bibitem{EZ}
E.~Zatorska.
\newblock On a steady flow of multicomponent, compressible, chemically reacting
  gas.
\newblock {\em Nonlinearity}, 24:3267--3278, 2011.

\bibitem{EZ2}
E.~Zatorska.
\newblock On the flow of chemically reacting gaseous mixture.
\newblock {\em J. Differential Equations}, 253(12):3471--3500, 2012.

\bibitem{EZ3}
E.~Zatorska.
\newblock Mixtures: sequential stability of variational entropy solutions.
\newblock {\em arXiv:1312.4391}, 2013.

\end{thebibliography}
%
%
%
%
%

\end{document}